\documentclass[a4paper, 10pt, final]{amsart}

\usepackage{amsthm}

\theoremstyle{plain}

\usepackage[latin1]{inputenc}
\usepackage[T1]{fontenc}
\usepackage[english]{babel}
\usepackage{fullpage}
\usepackage{amsmath}
\usepackage{amsfonts}
\usepackage{amssymb}
\usepackage{amscd}
\usepackage{mathrsfs}

\usepackage[OT2,T1]{fontenc}

\usepackage{amsmath, amsthm, amsfonts}

\usepackage{lmodern}
\usepackage{graphicx}
\usepackage{tikz}
\usetikzlibrary{calc}
\usepackage{pgfplots}
\usepgfplotslibrary{groupplots}
\usepgfplotslibrary{units}

\usetikzlibrary{arrows}
\usetikzlibrary{graphs,shapes.geometric,arrows,fit,matrix,positioning}

\DeclareMathOperator{\weight}{weight}
\DeclareMathOperator{\id}{id}

\DeclareMathOperator{\Sign}{Sign}

\DeclareMathOperator{\depth}{depth}

\DeclareMathOperator{\Ad}{Ad}
\DeclareMathOperator{\DR}{DR}
\DeclareMathOperator{\un}{un}

\DeclareMathOperator{\Wd}{Wd}

\DeclareMathOperator{\Spec}{Spec}

\DeclareMathOperator{\Har}{Har}
\DeclareMathOperator{\loc}{loc}

\DeclareMathOperator{\har}{har}

\DeclareMathOperator{\SubWd}{SubWd}

\DeclareMathOperator{\Li}{Li}

\DeclareMathOperator{\comp}{comp}

\DeclareMathOperator{\KZ}{KZ}
\DeclareMathOperator{\an}{an}

\DeclareMathOperator{\shft}{shft}
\DeclareMathOperator{\deloc}{deloc}

\theoremstyle{plain}

\newtheorem{Theoreme}{Theoreme}[subsection]
\newtheorem{Proposition}[Theoreme]{Proposition}

\newtheorem{Proposition-Definition}[Theoreme]{Proposition-Definition}
\newtheorem{Lemma-Notation}[Theoreme]{Lemma-Notation}
\newtheorem{Lemma-Definition}[Theoreme]{Lemma-Definition}

\newtheorem{Nota Bene}[Theoreme]{Nota Bene}
\newtheorem{Corollary}[Theoreme]{Corollary}
\theoremstyle{definition}
\newtheorem{Definition-Notation}[Theoreme]{Definition-Notation}

\newtheorem{Lemma}[Theoreme]{Lemma}
\newtheorem{Definition}[Theoreme]{Definition}
\newtheorem{Remark}[Theoreme]{Remark}

\newtheorem{Example}[Theoreme]{Example}
\newtheorem{Examples}[Theoreme]{Examples}

\newtheorem{Notation}[Theoreme]{Notation}

\input cyracc.def
\DeclareFontFamily{U}{russian}{}
\DeclareFontShape{U}{russian}{m}{n}
        { <5><6> wncyr5
        <7><8><9> wncyr7
        <10><10.95><12><14.4><17.28><20.74><24.88> wncyr10 }{}
\DeclareSymbolFont{Russian}{U}{russian}{m}{n}
\DeclareSymbolFontAlphabet{\mathcyr}{Russian}
\makeatletter
\let\@math@cyr\mathcyr
\renewcommand{\mathcyr}[1]{\@math@cyr{\cyracc #1}}
\makeatother
\newcommand{\sh}{\mathcyr{sh}} 

\numberwithin{equation}{subsection}

\makeatletter
\newcounter{subsubsubsection}[subsubsection]
\renewcommand\thesubsubsubsection{\thesubsubsection .\@alph\c@subsubsubsection}
\newcommand\subsubsubsection{\@startsection{subsubsubsection}{4}{\z@}%
                                     {-3.25ex\@plus -1ex \@minus -.2ex}%
                                     {1.5ex \@plus .2ex}%
                                     {\normalfont\normalsize\bfseries}}
\newcommand*\l@subsubsubsection{\@dottedtocline{3}{10.0em}{4.1em}}
\newcommand*{\subsubsubsectionmark}[1]{}
\makeatother

\usepackage[top=2.5cm, bottom=2.5cm, left=2.5cm, right=2.5cm]{geometry}

\author{David Jarossay}

\address{De Vinci Higher Education, De Vinci Research Center, Paris, France}
\email{david.jarossay@devinci.fr}

\begin{document}

\title{Pro-unipotent harmonic actions and a computation of $p$-adic cyclotomic multiple zeta values}

\maketitle
\noindent

\begin{abstract}
We obtain formulas relating $p$-adic cyclotomic multiple zeta values and cyclotomic multiple harmonic sums. In particular, we obtain a series formula for $p$-adic cyclotomic multiple zeta values, and conversely a formula for certain cyclotomic multiple harmonic sums in terms of $p$-adic cyclotomic multiple zeta values. Our formulas are related to the motivic framework via a new notion which we call pro-unipotent harmonic actions, which are ad hoc $p$-adic byproducts of the Ihara action.

As an application, we prove a conjecture of Akagi, Hirose and Yasuda on the relation between $p$-adic multiple zeta values and multiple harmonic sums, and we generalize it to the cyclotomic case. We also deduce bounds on the dimension of the spaces of finite cyclotomic multiple zeta values.

This is Part I-2 of \emph{$p$-adic cyclotomic multiple zeta values and $p$-adic pro-unipotent harmonic actions}.
\end{abstract}

\tableofcontents

\numberwithin{equation}{subsection}

\section{Introduction}

\subsection{Complex and $p$-adic cyclotomic multiple zeta values}

Cyclotomic multiple zeta values (CMZVs) are the following complex numbers. Let $N \in \mathbb{N}_{\geq 1}$. For any $d \in \mathbb{N}_{\geq 1}$,  $n_{i} \in \mathbb{N}_{\geq 1}$, ($1 \leqslant i \leqslant d$) and $\xi_{i}$ $N$-th roots of unity ($1 \leqslant i \leqslant d$) such that $(\xi_{d},n_{d}) \not= (1,1)$ :
\begin{equation}
\label{eq:multizetas series} \zeta\big((n_{i})_{d}; (\xi_{i})_{d}\big) =
\sum_{0<m_{1}<\ldots<m_{d}} \frac{\big( \frac{\xi_{2}}{\xi_{1}}\big)^{m_{1}} \ldots \big( \frac{1}{\tilde{\xi_{d}}}\big)^{m_{d}}}{m_{1}^{n_{1}} \ldots m_{d}^{n_{d}}}.
\end{equation}
Here $n = n_{d}+\ldots+n_{1}$ resp. $d$ is called the weight, resp. the depth of $\big((n_{i})_{d};(\xi_{i})_{d}\big) = \big(n_{1},\ldots,n_{d}; \xi_{1},\ldots,\xi_{d} \big)$.
\newline One has the following integral formula, where we denote the sequence $(\underbrace{0,\ldots,0}_{n_{d}-1},\xi_{d},\ldots,\underbrace{0,\ldots,0}_{n_{1}-1},\xi_{1})$ by $(\epsilon_{n},\ldots,\epsilon_{1})$ :
\begin{equation} \label{eq:multizetas integrals}
\zeta \big((n_{i})_{d}; (\xi_{i})_{d}\big) = (-1)^{d}\int_{t_{n}=0}^{1} \frac{dt_{n}}{t_{n} - \epsilon_{n}} \int_{t_{n-1}=0}^{t_{n}} \ldots \int_{t_{1}=0}^{t_{2}} \frac{dt_{1}}{t_{1} - \epsilon_{1}}.
\end{equation}
\indent Equation (\ref{eq:multizetas integrals}) shows that cyclotomic multiple zeta values are Betti-De Rham periods of the pro-unipotent fundamental groupoid ($\pi_{1}^{\un}$) of $\mathbb{P}^{1} \setminus \{0,\mu_{N},\infty\}$ (\cite{Deligne Goncharov}, \S5.16).
\newline\indent Let $p$ be a prime number which does not divide $N$, and let $K$ be the extension of $\mathbb{Q}_{p}$ generated by a primitive $N$-th root of unity. $p$-adic cyclotomic multiple zeta values ($p$CMZVs) are numbers $\zeta_{p,\alpha}\big((n_{i})_{d};(\xi_{i})_{d}\big)\in K $ defined as $p$-adic analogues of the integrals (\ref{eq:multizetas integrals}), where $\alpha$ is any non-zero integer (\cite{Deligne Goncharov}, \cite{Yamashita}, \cite{U1}, \cite{U2}, \cite{I-1} Definition 1.2.2). A different but equivalent notion (\cite{Furusho 1}, \cite{Furusho 2}, \cite{Yamashita}) defines $p$-adic cyclotomic multiple zeta values $\zeta_{p}^{\KZ}\big((n_{i})_{d};(\xi_{i})_{d}\big)\in K$ as the Coleman integrals analogous to (\ref{eq:multizetas integrals}). These definitions do not come with an explicit formula. Both notions of $p$-adic integrals refer to the Frobenius structure of the KZ differential equation on $\mathbb{P}^{1} \setminus \{0,\mu_{N},\infty\}$ (\ref{eq:KZ equation}) which is the connection canonically associated with $\pi_{1}^{\un,\DR}(\mathbb{P}^{1} \setminus \{0,\mu_{N},\infty\})$ in the sense of \cite{Deligne}. That Frobenius structure characterizes the crystalline resp. rigid pro-unipotent fundamental groupoid of $\mathbb{P}^{1} \setminus \{0,\mu_{N},\infty\}$ in the sense of \cite{Deligne} \S13.6, \cite{Shiho 1} \cite{Shiho 2}, resp. \cite{CLS}.
\newline\indent By \cite{Yamashita}, $p$CMZVs are reductions of the periods associated with the crystalline pro-unipotent fundamental groupoid of $\mathbb{P}^{1} \setminus \{0,\mu_{N},\infty\}$.

\subsection{The question of computing $p$-adic cyclotomic multiple zeta values}

The question of finding a convenient $p$-adic analogue of (\ref{eq:multizetas series}), which would be an explicit formula, has been raised first by Deligne in 2002 in the $N=1$ case and appears in \cite{Deligne Goncharov}, \S5.28. In the complex case, equation (\ref{eq:multizetas series}) can be viewed as the value at $z =1$ of the power series expansion at $0$ of multiple polylogarithms, which are solutions to the KZ equation : for any positive integers
$d$ and $n_{i}$ ($1 \leqslant i \leqslant d$) and for any roots of unity $\xi_{i}$ ($1 \leqslant i \leqslant d$), for $z \in \mathbb{C}$ such that $|z| < 1$,
\begin{equation} \label{eq:Li} \Li\big((n_{i})_{d};(\xi_{i})_{d} \big)(z) = \sum_{0<m_{1}<\ldots<m_{d}} \frac{\big( \frac{\xi_{2}}{\xi_{1}} \big)^{m_{1}} \ldots \big( \frac{z}{\xi_{d}} \big)^{m_{d}}}{m_{1}^{n_{1}} \ldots m_{d}^{n_{d}}}.
\end{equation}
\indent In the $p$-adic case, one has $p$-adic multiple polylogarithms, defined by Coleman integration (\cite{Furusho 1} \cite{Furusho 2} \cite{Yamashita}), solutions to the KZ equation and also admitting (\ref{eq:Li}) as a power series expansion at $0$. Thus, a $p$-adic analogue of (\ref{eq:multizetas series}) would mean a formula for $p$CMZVs in terms of the coefficients of the power series expansion (\ref{eq:Li}) ; for us, this will mean in terms of the weighted multiple harmonic sums (let $m$ be a positive integer) :
\begin{equation} \label{eq:multiple harmonic sum}
\har_{m} \big((n_{i})_{d}; (\xi_{i})_{d+1}\big) = m^{n_{d}+\ldots+n_{1}} \sum_{0<m_{1}<\ldots<m_{d}<m} 
\frac{ \big( \frac{\xi_{2}}{\xi_{1}} \big)^{m_{1}} \ldots \big( \frac{\xi_{d+1}}{\xi_{d}} \big)^{m_{d}}  
\big(\frac{1}{\xi_{d+1}} \big)^{m}}{m_{1}^{n_{1}}\ldots m_{d}^{n_{d}}}.
\end{equation}
\indent The power series expansion in (\ref{eq:Li}) converges for $z \in \mathbb{C}_{p}$ such that $|z|_{p}<1$. However, $\mathbb{C}_{p}$ is totally disconnected and one cannot take the limit of (\ref{eq:Li}) when $z \rightarrow 1$ in $\mathbb{C}_{p}$. This is what makes not immediate to find a $p$-adic analogue of (\ref{eq:multizetas series}).
\newline\indent The overconvergence of the Frobenius of the KZ equation provides a certain substitute to the operation $\underset{z \rightarrow 1}{\lim}$ in (\ref{eq:Li}), and gives a sort of $p$-adic analogue of (\ref{eq:multizetas series}). This has been used in \cite{U1}, \cite{U2}, \cite{U3}, \cite{U4}. However, the formulas obtained are very recursive and it seems difficult to read them and use them, because it requires to deal with the differential equation satisfied by the Frobenius (equation (\ref{eq:first})), which is complicated combinatorially.
\newline\indent A hope for the existence of simpler formulas for $p$CMZVs is provided by Kaneko-Zagier's work on finite multiple zeta values and the conjecture of Akagi-Hirose-Yasuda inspired by that work (see \S0.3 and \S0.4 for details) ; it is also motivated by a question asked by Deligne and Goncharov (\cite{Deligne Goncharov}, \S5.28). We propose in this paper a realization of this hope. This opens an explicit theory of $p$CMZVs.
\newline\indent The starting point is the observation that the equation satisfied by the Frobenius (\ref{eq:first}) is surprisingly constant with respect to one of its parameters when reformulated in a convenient way (Proposition \ref{previous formula}). Combined with the bound on the norm of the Frobenius given by the main result of \cite{I-1}, this leads us to a big simplification of this equation (Lemma \ref{cor tend to 0}), which will allow us to compute $p$CMZVs in a different way.

\subsection{Summary of the paper}

We are going to work not in terms of each $p$CMZV, but in terms of their non-commutative generating series $\Phi_{p,\alpha}$ (equation (\ref{eq:Phipalpha})) which is a $K$-point of the affine scheme $\Pi_{1,0}=\pi_{1}^{\un,\DR}(\mathbb{P}^{1} \setminus \{0,\mu_{N},\infty\},-\vec{1}_{1},\vec{1}_{0})$, and we will encode our computation by a new structure on $\pi_{1}^{\un,\DR}(\mathbb{P}^{1} \setminus \{0,\mu_{N},\infty\})$. This structure will keep track of the motivic Galois theory. In the complex setting, CMZVs are periods and their non-commutative generating series represents a point of a torsor under the motivic Galois group of a Tannakian category of mixed Tate motives. In the $p$-adic setting, the situation is different : the generating series of $p$CMZVs represents a point of that motivic Galois group. Our computation will keep track of this feature as follows. 
\newline\indent With the simplification evoked in \S0.3, we will replace the Frobenius by a simpler variant which we will call the \emph{harmonic Frobenius of integrals} (Definition \ref{def har Frob}), which we will view as an operation on the space of the weighted multiple harmonic sums (\ref{eq:multiple harmonic sum}). The passage from the Frobenius to the harmonic Frobenius will be lifted to a passage from the Ihara action, which is a byproduct of the motivic Galois action useful to express the Frobenius (equation (\ref{eq:Ihara})) to a new object, the \emph{pro-unipotent harmonic action of integrals} $\circ_{\har}^{\smallint}$ (Definition \ref{dR-rt harmonic Ihara action}). We will construct a torsor under $\circ_{\har}^{\smallint}$ and we will characterize $p$CMZVs in terms of the action $\circ_{\har}^{\smallint}$ on this torsor.
\newline\indent The definition of $\circ_{\har}^{\smallint}$, which is not an algebraic operation but involves infinite summations, will be prepared by \S1, where we will make out of $\pi_{1}^{\un,\DR}(X)$, which is a groupoid in pro-affine schemes over $X$, a groupoid in complete ultrametric $K$-algebras (Proposition \ref{corollary le corollary}), whose topologies are induced by certain norms (Definition \ref{def norms}), and which includes a notion of ``summable elements'' (Definition \ref{def summable}).
\newline\indent In \S4 we construct a \emph{harmonic Frobenius of series} on weighted multiple harmonic sums (Definition \ref{def harmonic Frobenius of series}) in an elementary way using explicit sums of series (Definition \ref{def harmonic Frobenius of series}). It involves to define a \emph{pro-unipotent harmonic action of series} $\circ_{\har}^{\Sigma}$ (Proposition-Definition \ref{series harmonic action}). The construction involves a notion called localized multiple harmonic sums, which is introduced and studied as a preliminary in \S3 (Definition \ref{def localized mult har sums}).
\newline\indent Having obtained two expressions of the harmonic Frobenius, it remains to say that they are equal. This is the purpose of \S5. We define maps of comparison between series and integrals, $\comp^{\smallint\Sigma}$ and $\comp^{\Sigma\smallint}$ (Definition \ref{def comp from series to integrals}, Definition \ref{def comp from integrals to series}). We show that they enable to relate $\circ_{\har}^{\smallint}$ and $\circ_{\har}^{\Sigma}$.
\newline\indent Below, the $e_{x}$'s where $x \in \{0\} \cup \mu_{N}(K)$ are generators of the Lie algebra of $\pi_{1}^{\un,\DR}(\mathbb{P}^{1} \setminus \{0,\mu_{N},\infty\},\vec{1}_{0})$ and $\har_{\mathbb{N}}$, $\har_{p^{\alpha}\mathbb{N}}$, $\har_{p^{\alpha}}$, $\har_{\mathbb{N}}^{(p^{\alpha})}$ are non-commutative generating series of weighted multiple harmonic sums (Definition \ref{def non comm gen seri wmhs}). 

\subsection{Main result and applications}

The main result is the following :
\newline
\newline \textbf{Theorem.}
\noindent\newline (i) \emph{(integrals) The $p$-adic pro-unipotent harmonic action of integrals is a continuous group action, there exists a torsor for $\circ_{\har}^{\smallint}$ containing $\har_{\mathbb{N}}^{(p^{\alpha})}$, and we have
\begin{equation} \label{eq: property of ac 01} \har_{p^{\alpha}\mathbb{N}} = \big(\Ad_{\Phi_{p,\alpha}^{(\xi)}}(e_{\xi})\big)_{\xi \in \mu_{N}(K)} 
\text{ } \circ_{\har}^{\smallint} \text{ }
\har_{\mathbb{N}}^{(p^{\alpha})}.
\end{equation}
(ii) (series) The $p$-adic pro-unipotent harmonic action of series is continuous and we have
\begin{equation} \label{eq: property ac sigma} \har_{p^{\alpha}\mathbb{N}} = \har_{p^{\alpha}}
\text{ } \circ_{\har}^{\Sigma} \text{ }
\har_{\mathbb{N}}^{(p^{\alpha})} 
\end{equation}
\noindent (iii) (comparison between integrals and series) The maps of comparison satisfy, for $h$ in the orbit of $\har_{\mathbb{N}}^{(p^{\alpha})}$, and for any $g$,
\begin{equation} \label{eq:comparison 1} g \circ_{\har}^{\smallint} \comp^{\smallint\Sigma} h = g \circ_{\har}^{\Sigma} h,
\end{equation}
and
\begin{equation} \label{eq:comparison 2}\comp^{\Sigma\smallint} \circ \comp^{\smallint\Sigma} = \id ,
\end{equation}
\begin{equation} \label{eq:series expansion} \big( \Ad_{\Phi_{p,\alpha}^{(\xi)}}(e_{\xi})\big)_{\xi \in \mu_{N}(K)}
= \comp^{\smallint\Sigma} (\har_{p^{\alpha}}),
\end{equation}
\begin{equation} \label{eq:inversion of series expansion} \har_{p^{\alpha}} = \comp^{\Sigma\smallint} \big( \Ad_{\Phi_{p,\alpha}^{(\xi)}}(e_{\xi})\big)_{\xi \in \mu_{N}(K)}.
\end{equation}}

The simplest terms of equations (\ref{eq: property of ac 01}) and (\ref{eq: property ac sigma}) (depth $\leq 2$, in $N=1$ case) are written explicitly without the combinatorial tools used in this paper, respectively by Example \ref{example circ har int} and Example \ref{example theorem part Sigma}. Equation (\ref{eq:series expansion}) is an expression of $p$CMZVs in terms of prime weighted cyclotomic multiple harmonic sums, i.e. the numbers (\ref{eq:multiple harmonic sum}) with $m=p^{\alpha}$. Its explicit version can be obtained by combining the formulas of Proposition-Definition \ref{loc series harmonic action} and Proposition \ref{formula for loc}. Equation (\ref{eq:inversion of series expansion}) is an expression of prime weighted cyclotomic multiple harmonic sums in terms of $p$CMZVs. 
Equation (\ref{eq:inversion of series expansion}) is actually a particular case of equation (\ref{eq: property of ac 01}), via the fact that all multiple harmonic sums $\har_{1}(w)$, being an iterated sums as in (\ref{eq:multiple harmonic sum}) on an empty domain of summation, vanish ; it is also obtained by joining (\ref{eq:comparison 2}) and (\ref{eq:series expansion}).
\newline\indent The explicit version of equation (\ref{eq:inversion of series expansion}) is the following (the notation $\zeta_{p,\alpha}^{(\xi)}$ is introduced in \S1.1.3 ; see also Notation \ref{notation words}) :
\begin{multline} \label{eq:formula for n=1}
\har_{p^{\alpha}} \big((n_{i})_{d}; (\xi_{i})_{d+1}\big) = (-1)^{d} \sum_{\xi \in \mu_{N}(K)} \xi^{-p^{\alpha}} \big( {\Phi^{(\xi)}_{p,\alpha}}^{-1}e_{\xi}\Phi^{(\xi)}_{p,\alpha}\big)\bigg[\frac{1}{1-e_{0}}e_{\xi_{d+1}}e_{0}^{n_{d}-1}e_{\xi_{d}}\ldots e_{0}^{n_{1}-1}e_{\xi_{1}} \bigg]
\\ = \sum_{d'=0}^{d} \sum_{l_{d'+1},\ldots,l_{d} = 0}^{\infty} \xi_{d-d'+1}^{p^{\alpha}} \bigg( \prod_{i=d'+1}^{d} (-1)^{n_{i}} {-n_{i} \choose l_{i}}
\bigg) \zeta^{(\xi_{d'+1})}_{p,\alpha} 
\big((n_{d'+i}+l_{d'+i})_{d-d'} ; (\xi_{d'+1+i})_{d-d'} \big)
\text{ }\zeta^{(\xi_{d'+1})}_{p,\alpha} \big((n_{i})_{d'} ; (\xi_{i})_{d'}\big),
\end{multline}
\noindent in particular, in the case of $\mathbb{P}^{1} \setminus \{0,1,\infty\}$,
\begin{multline}
\label{eq:explicit inversion of series expansion N=1} \har_{p^{\alpha}}(n_{1},\ldots,n_{d}) = (-1)^{d} (\Phi_{p,\alpha}^{-1}e_{1}\Phi_{p,\alpha})\bigg[\frac{1}{1-e_{0}}e_{1}e_{0}^{n_{d}-1}e_{1}\ldots e_{0}^{n_{1}-1}e_{1} \bigg] \\ = \sum_{d'=0}^{d} \sum_{l_{d'+1},\ldots,l_{d} =0}^{\infty} \bigg(\prod_{i=d'}^{d} (-1)^{n_{i}} {-n_{i} \choose l_{i}} \bigg) \zeta_{p,\alpha}(n_{d}+l_{d},\ldots,n_{d'+1}+l_{d'+1})\zeta_{p,\alpha}(n_{1},\ldots,n_{d'}).
\end{multline}

The case $\alpha=1$, $N=1$ and depth $1$ of equation (\ref{eq: property ac sigma}) was known by Boyd (\cite{Boyd}, Theorem 5.2). The case $\alpha=1$ and $d=1$ of equation (\ref{eq:explicit inversion of series expansion N=1}) was known by a result of Washington (\cite{Washington}, Theorem 1 (a)) combined to a result of Coleman \cite{Coleman} (equation 4 p. 173). Akagi, Hirose and Yasuda had conjectured the $\alpha=1$ case of equation (\ref{eq:explicit inversion of series expansion N=1}) and M. Hirose had proved it for $\alpha=1$ and $d=2$ \cite{Yasuda}.

We also give an application to finite cylcotomic multiple zeta values in \S6. This generalizes an application due to Akagi-Hirose-Yasuda in the $N=1$ case \cite{Yas2}. 

The formulas of this paper keep track of the motivic Galois action by the pro-unipotent harmonic actions. We will find an algebraic and motivic background behind these formulas in next papers \cite{II-1,II-2,II-3}.

\emph{Acknowledgements.} This work has been done at Universit\'{e} Paris Diderot with the support of ERC grant 257638, then has been extended and revised at Universit\'{e} de Strasbourg with the support of Labex IRMIA and at Universit\'{e} de Gen\`{e}ve with the support of NCCR SwissMAP, and at De Vinci Research Center in Paris. I thank Seidai Yasuda and Francis Brown for having transmitted to me \cite{Yasuda} which contained the statement of the conjecture of Akagi, Hirose and Yasuda mentioned in \S0.3. I also thank Pierre Cartier and Ahmed Abbes for encouragements.

\section{Setting for the pro-unipotent harmonic action of integrals}

This section is a prerequisite for \S2. We review the combinatorics of some operations on $\pi_{1}^{\un}(\mathbb{P}^{1} \setminus \{0,\mu_{N},\infty\})$, and we define a few operations and a topological structure on  $\pi_{1}^{\un,\DR}(\mathbb{P}^{1} \setminus \{0,\mu_{N},\infty\})(K)$.
In all this paper we denote by $\mathbb{N}$ resp. $\mathbb{N}^{\ast}$ the set of nonnegative resp. positive integers.

\subsection{Review on $\pi_{1}^{\un}(\mathbb{P}^{1} \setminus \{0,\mu_{N},\infty\})$, and an adjoint Ihara action}

\subsubsection{The De Rham unipotent fundamental groupoid of $\mathbb{P}^{1} \setminus \{0,\mu_{N},\infty\}$}

Let $X$ be $\mathbb{P}^{1} \setminus \{0,\mu_{N},\infty\}$ over a field $K$ of characteristic $0$ containing a primitive $N$-th root of unity. Let $\pi_{1}^{\un,\DR}(X)$ be the De Rham realization of the unipotent fundamental groupoid of $X$ (\cite{Deligne}, \S10.27, \S10.30,(ii)). It is a groupoid in pro-affine schemes on $X$. Its base-points are the points of $X$, the tangential base-points of $X$ i.e. the non-zero tangent vectors $\vec{v}_{x}$ at a point $x\in \{0\}\cup \mu_{N}(K)\cup \{\infty\}$, (\cite{Deligne}, \S15) and the canonical base-point $\omega_{\DR}$ (\cite{Deligne}, (12.4.1)).
\newline\indent Let $e_{0 \cup \mu_{N}}$ be the alphabet $\{e_{x}\text{ }|\text{ }x \in \{0\} \cup \mu_{N}(K)\}$. Let $\mathcal{O}^{\sh}$ be the shuffle Hopf algebra on $e_{0 \cup \mu_{N}}$. It is a Hopf algebra over $\mathbb{Q}$ whose underlying vector space admits as a basis the set of words on $e_{0 \cup \mu_{N}}$, including the empty word, and whose product is the shuffle product of words on $e_{0 \cup \mu_{N}}$, denoted by $\sh$. The weight of a word on $e_{0 \cup \mu_{N}}$ is its number of letters. We usually write a word on $e_{0\cup\mu_{N}}$ in the form $e_{0}^{n_{d}-1}e_{\xi_{d}} \ldots e_{0}^{n_{1}-1}e_{\xi_{1}}e_{0}^{n_{0}-1}$ where $d$ and the $n_{i}$'s ($0\leqslant i \leqslant d$) are positive integers and the $\xi_{i}$'s ($1\leqslant i \leqslant d$) are $N$-th roots of unity. For most computations it is sufficient to consider the words such that $n_{0}=1$. The depth of a word $w$ on $e_{0 \cup \mu_{N}}$ is its number $d$ of letters distinct from $e_{0}$.
\newline\indent The pro-unipotent affine group scheme $\pi_{1}^{\un,\DR}(X,\omega_{\DR})$ is canonically isomorphic to $\Spec(\mathcal{O}^{\sh,e_{0 \cup \mu_{N}}})$ (by \cite{Deligne}, \S12.9). Let $K \langle\langle e_{0 \cup \mu_{N}} \rangle\rangle$ be the non-commutative $K$-algebra of formal power series with variables the letters of $e_{0 \cup \mu_{N}}$ and coefficients in $K$.
\begin{Notation} \label{notation words}An element $f$ of $K \langle \langle e_{0 \cup \mu_{N}} \rangle\rangle$ can be written in a unique way as
$f = \sum\limits_{w\text{ word on }e_{0\cup \mu_{N}}} f[w]w$ i.e. for any word $w$ on $e_{0 \cup \mu_{N}}$, we denote by $f[w] \in K$ the coefficient of $w$ in $f$.
\end{Notation}
We have a canonical inclusion $\pi_{1}^{\un,\DR}(X,\omega_{\DR})(K) \subset K \langle \langle e_{0 \cup \mu_{N}} \rangle\rangle$, whose image is the group of formal power series $f$ satisfying $f[\emptyset]=1$ and the shuffle equation, i.e. $f[w\text{ }\sh\text{ }w'] = f[w] f[w']$ for all words $w,w'$ on $e_{0 \cup \mu_{N}}$.
\newline\indent For any base-points $x,y$, the scheme $\pi_{1}^{\un,\DR}(X,y,x)$ is canonically isomorphic to $\pi_{1}^{\un,\DR}(X,\omega_{\DR})$ and these isomorphisms are compatible with the groupoid maps $\pi_{1}^{\un,\DR}(X,z,y)\times \pi_{1}^{\un,\DR}(X,y,x) \rightarrow \pi_{1}^{\un,\DR}(X,z,x)$ (\cite{Deligne} \S12). The image of $1 \in \pi_{1}^{\un,\DR}(X,\omega_{\DR})(K)$ in $\pi_{1}^{\un,\DR}(X,y,x)(K)$ is denoted by ${}_y 1_{x}$, and called the canonical path from $x$ to $y$.
\newline\indent The KZ connection on $\mathbb{P}^{1} \setminus \{0,\mu_{N},\infty\}$ is the connection on $\pi_{1}^{\un,\DR}(X,\omega_{\DR}) \times X$ defined as follows (\cite{Deligne}, \S12.4) :
\begin{equation} \label{eq:KZ equation} \nabla_{\KZ} : f \mapsto df - \bigg(  \frac{dz}{z} e_{0} + \sum_{\xi \in \mu_{N}(K)} \frac{dz}{z - \xi} e_{\xi} \bigg) f .
\end{equation}

\begin{Notation} \label{notation DG 5} (i) (\cite{Deligne Goncharov}, \S5)
For all base-points $x,y$, let $\Pi_{y,x} =\pi_{1}^{\un,\DR}(\mathbb{P}^{1} \setminus \{0,\mu_{N},\infty\},y,x)$.	
For $x,y \in \{0\} \cup \mu_{N}(K)$, let $\Pi_{y,x}=\Pi_{\vec{1}_{y},\vec{1}_{x}}$. Let $\Pi = \pi_{1}^{\un,\DR}(\mathbb{P}^{1} \setminus \{0,\mu_{N},\infty\},\omega_{\DR})$.
\newline (ii) For any point $g \in \Pi_{1,0}(K)$, let $g^{(\xi)}$ be the element of $\Pi_{\xi,0}(K)$ obtained from $g$ by functoriality of $\pi_{1}^{\un,\DR}$ with respect to the automorphism $x \mapsto \xi x$ of $X$.
We will sometimes identify $g$ and the sequence $(g^{(\xi)})_{\xi \in \mu_{N}(K)}$.
\end{Notation}

\subsubsection{Some byproducts of the motivic Galois action on $\pi_{1}^{\un,\DR}(\mathbb{P}^{1} \setminus \{0,\mu_{N},\infty\})$}

The operations reviewed below will be used to express the pro-unipotent harmonic actions.
\newline\indent Let $G_{\omega}$ be the motivic Galois group defined as the Tannakian group associated with the category of mixed Tate motives over the  $N$-th cyclotomic field which are unramified at primes $p$ prime to $N$ (\cite{Deligne Goncharov}, \S1.6) and the canonical fiber functor $\omega$ (\cite{Deligne Goncharov}, \S1.1). We have $ G_{\omega} = \mathbb{G}_{m} \ltimes U_{\omega}$, where $U_{\omega}$ is a pro-unipotent algebraic group (\cite{Deligne Goncharov}, \S2.1.2).
\newline\indent By \cite{Deligne Goncharov} \S5, $G_{\omega}$ acts on $\Pi_{1,0}$, and this action encodes the motivic Galois theory of CMZVs. By \cite{Yamashita}, it also encodes the motivic Galois theory of $p$CMZVs, with the only difference that the $p$-adic analogue of $\zeta(2)$ is zero. This action is described as follows. The action of $\mathbb{G}_{m}$ on $\Pi_{1,0}$, and more generally, on any $\pi_{1}^{\un,\DR}(X,y,x)$, is given by
\begin{equation} \label{eq:tau} \tau : \begin{array}{cc} \mathbb{G}_{m} \times \pi_{1}^{\un,\DR}(X,y,x) \rightarrow \pi_{1}^{\un,\DR}(X,y,x)
\\ \big( \lambda,f(e_{0},(e_{\xi})_{\xi \in \mu_{N}(K)})\big) \mapsto f(\lambda e_{0},(\lambda e_{\xi})_{\xi \in \mu_{N}(K)})
\end{array}
\end{equation}
i.e. applying $\tau(\lambda)$ multiplies the terms of weight $n$ of an element $f$ by $\lambda^{n}$, for all $n$. We will also denote by $\tau$ the action on $K \langle\langle e_{0 \cup \mu_{N}} \rangle\rangle$ defined in the same way. The action of $U_{\omega}$ on $\Pi_{1,0}$ makes $\Pi_{1,0}$ into a torsor under a quotient $V_{\omega}$ of $U_{\omega}$ (\cite{Deligne Goncharov}, \S5.12), in such a way that the isomorphism of schemes $V_{\omega} \simeq \Pi_{1,0}$ obtained by choosing the canonical path ${}_{\vec{1}_{1}} 1_{\vec{1}_{0}}$ of $\Pi_{1,0}$ (in the sense reviewed in \S1.1.1) identifies the action of $V_{\omega}$ with 

\begin{equation} \label{eq:Ihara}
\circ^{\smallint_{1,0}} :
\begin{array}{cc} \Pi_{1,0} \times \Pi_{1,0} \rightarrow \Pi_{1,0}
\\
(g,f) \mapsto g \circ^{\smallint_{1,0}} f = g(e_{0},(e_{\xi})_{\xi \in \mu_{N}(K)}) \times f\big(e_{0},(\Ad_{g^{(\xi)}}(e_{\xi}))_{\xi \in \mu_{N}(K)}\big)
\end{array}
\end{equation}

where, because of our convention of reading the multiplication of the groupoid $\pi_{1}^{\un,\DR}(X)$ from the right to the left, we take the convention that $\Ad(e_{\xi})$ is $f \mapsto f^{-1}e_{\xi}f$. The group law $\circ^{\smallint_{1,0}}$  is sometimes called the twisted Magnus product, or the Ihara product or the Ihara action, at the base-points $(\vec{1}_{1},\vec{1}_{0})$ (our notation $\circ^{\smallint_{1,0}}$ is not standard). It induces a group law $\circ^{\xi,0}$ on $\Pi_{\smallint_{\xi,0}}$ for all $\xi \in \mu_{N}(K)$, by functoriality of $\pi_{1}^{\un,\DR}$. By the same isomorphism $V_{\omega} \simeq \Pi_{1,0}$, the action of $V_{\omega}$ on $\Pi_{0,0}$ induced by the motivic Galois action on $\Pi_{0,0}$ is identified with 

\begin{equation} \label{eq:Ihara bis}
\circ^{\smallint_{0,0}} :
\begin{array}{cc} \Pi_{1,0} \times \Pi_{0,0} \rightarrow \Pi_{0,0}
\\
(g,f) \mapsto g \circ^{\smallint_{0,0}} f = f\big(e_{0},(\Ad_{g^{(\xi)}}(e_{\xi}))_{\xi \in \mu_{N}(K)}\big)
\end{array}
\end{equation}

which we call the Ihara action at the base-point $\vec{1}_{0}$.
\newline\indent We now introduce a push-forward by $\Ad(e_{1})$ of the actions $\circ^{\smallint_{1,0}}$ and $\circ^{\smallint_{0,0}}$ :

\begin{Definition} \label{def Ad action 1,0}
(i) Let the adjoint Ihara action at the base-points $(\vec{1}_{1},\vec{1}_{0})$ be the map $\circ^{\smallint_{1,0}}_{\Ad} : \Ad_{\Pi_{1,0}}(e_{1}) \times \Ad_{\Pi_{1,0}}(e_{1}) \rightarrow \Ad_{\Pi_{1,0}}(e_{1})$, 
$(g,f) \mapsto f(e_{0},(g^{(\xi)})_{\xi \in \mu_{N}(K)})$.
\newline (ii) \label{def Ad action 0,0} Let the adjoint Ihara action at the base-point $\vec{1}_{0}$ be the map $\circ_{\Ad}^{\smallint_{0,0}} : \Ad_{\Pi_{1,0}}(e_{1}) \times \Pi_{0,0}\rightarrow \Pi_{0,0}$, $(h,f) \mapsto f(e_{0},(h^{(\xi)})_{\xi \in \mu_{N}(K)})$. We will also denote by $\circ_{\Ad}^{\smallint_{0,0}}$ the map $ K \langle \langle e_{0 \cup \mu_{N}} \rangle\rangle \times K \langle \langle e_{0 \cup \mu_{N}} \rangle\rangle \rightarrow  K \langle \langle e_{0 \cup \mu_{N}} \rangle\rangle$ defined by the same formula.
\end{Definition}

\begin{Proposition} \label{adjoint morphism}(i) $(\Ad_{\Pi_{1,0}}(e_{1}),\circ_{\Ad}^{\smallint_{1,0}})$ is a group scheme such that $\Ad(e_{1})$ is an morphism of group schemes $(\Pi_{1,0},\circ^{\smallint_{1,0}}) \mapsto (\Ad_{\Pi_{1,0}}(e_{1}),\circ_{\Ad}^{\smallint_{1,0}})$.
\newline (ii) $\circ^{\smallint_{0,0}}_{\Ad}$ is an algebraic group action of $(\Ad_{\Pi_{1,0}}(e_{1}),\circ_{\Ad}^{\smallint_{1,0}})$ such that $\Ad(e_{1})$ induces a morphism $\circ^{\smallint_{0,0}} \mapsto \circ_{\Ad}^{\smallint_{0,0}}$ of algebraic group actions.
\end{Proposition}

\begin{proof} Follows directly from the formulas and from the fact that the composition of non-commutative formal power series is associative.
\end{proof}

\subsubsection{The Frobenius of $\pi_{1}^{\un,\DR}(X)$}

Let us now assume that $K$ is the extension of $\mathbb{Q}_{p}$ generated by a primitive $N$-th root of unity, where $p$ is a prime number which does not divide $N$. Let $\phi$ be the crystalline Frobenius of $\pi_{1}^{\un,\DR}(X)$ in the sense of \cite{Deligne}, \S13. Let $\sigma$ be the Frobenius automorphism of $K$, which generates the Galois group of $K/\mathbb{Q}_{p}$. It induces an automorphism of $K\langle \langle e_{0 \cup \mu_{N}} \rangle\rangle$ which we also denote by $\sigma$. Let a positive integer $\alpha$. The map $\tau(p^{\alpha}) \circ \phi^{\alpha}$ at base-points $(\vec{1}_{0},\vec{1}_{1})$ is of the form
\begin{equation} \label{eq: Frobenius at 10} \tau(p^{\alpha})\circ\phi^{\alpha} : \begin{array}{cc} \Pi^{(p^{\alpha})}_{1,0}(K) \rightarrow \Pi_{1,0}(K) 
\\ \text{ }\text{ }\text{ }\text{ }\text{ }\text{ }\text{ }\text{ }\text{ }f \mapsto \Phi_{p,\alpha} \circ^{\smallint_{1,0}} \sigma^{\alpha}(f)
\end{array}
\end{equation}
where $\Pi^{(p^{\alpha})}_{1,0}$ is the pull-back of $\Pi_{1,0}$ by $\sigma^{\alpha}$, and
\begin{equation} \label{eq:Phipalpha}
\Phi_{p,\alpha} = \tau(p^{\alpha})\circ \phi^{\alpha}(_{\vec{1}_{1}} 1^{(p^{\alpha})} _{\vec{1}_{0}}) \in \Pi^{(p^{\alpha})}_{1,0}(K).
\end{equation}
\newline\indent The numbers 
\begin{equation} \zeta_{p,\alpha} \big((n_{i})_{d}:(\xi_{i})_{d}\big) = (-1)^{d} \Phi_{p,\alpha}[e_{0}^{n_{d}-1}e_{\xi_{d}} \ldots e_{0}^{n_{1}-1}e_{\xi_{1}} ] \in K 
\end{equation}
with $d$ and the $n_{i}$'s ($1\leqslant i \leqslant d$) positive integers and the $\xi_{i}$'s ($1\leqslant i \leqslant d$) $N$-th roots of unity, are $p$-adic cyclotomic multiple zeta values. If $N=1$, they are $p$-adic multiple zeta values. We also denote, for any $\xi \in \mu_{N}(K)$ and word $w$ on $e_{0\cup \mu_{N}}$, by $\zeta_{p,\alpha}^{(\xi)}(w) = \Phi_{p,\alpha}^{(\xi)}[w]$, where $\Phi_{p,\alpha}^{(\xi)}$ is in the sense of Notation \ref{notation DG 5} (ii).
\newline\indent Let the affinoid rigid analytic space $U^{\an} = \mathbb{P}^{1,\an} \setminus \underset{\xi \in \mu_{N}(K)}{\cup} B(\xi,1)$ over $K$, where $B(\xi,1)$ is the open ball of center $\xi$ and radius $1$. Let, on that space, $\Li_{p,\alpha}^{\dagger}(z) = \tau(p^{\alpha})\phi^{\alpha}(_{z} 1 _{\vec{1}_{0}})$. The coefficients  $\Li_{p,\alpha}^{\dagger}[w]$ are overconvergent analytic functions on $U^{\an}$ called overconvergent $p$-adic multiple polylogarithms.
Let $X^{(p^{\alpha})}$ be the pull-back of $X$ by $\sigma^{\alpha}$. Let $\log_{p}$ be any determination of the $p$-adic logarithm. Let $\Li_{p,X}^{\KZ}$, resp. $\Li_{p,X^{(p^{\alpha})}}^{\KZ}$ (\cite{Furusho 1} for $N=1$, \cite{Yamashita} for any $N$) be the non-commutative generating series of Coleman functions on $X$, resp. $X^{(p^{\alpha})}$, which is a horizontal section of $\nabla_{\KZ}$ (\ref{eq:KZ equation}), resp. of the pull-back of $\nabla_{\KZ}$ by $\sigma^{\alpha}$, with the asymptotics $\Li_{p,X}^{\KZ}(z) \underset{z \rightarrow 0}{\sim} e^{ e_{0} \log_{p}(z)}$, resp. $\Li_{p,X^{(p^{\alpha})}}^{\KZ}(z) \underset{z \rightarrow 0}{\sim} e^{ e_{0} \log_{p}(z)}$. We have the following equation on $U^{\an}$ (\cite{I-1}, Proposition 2.2.21) :
\begin{multline} \label{eq:first} \Li_{p,\alpha}^{\dagger}(z)\big(e_{0},(e_{\xi})_{\xi \in \mu_{N}(K)}\big)\text{ }\Li_{p,X^{(p^{\alpha})}}^{\KZ}(z^{p^{\alpha}})\big(e_{0},
(\Ad_{{\Phi^{(\xi)}_{p,\alpha}}}(e_{\xi}))_{\xi \in \mu_{N}(K)} \big) = \Li_{p,X}^{\KZ}(z)\big(p^{\alpha}e_{0},(p^{\alpha}e_{\xi})_{\xi \in \mu_{N}(K)}\big),
\end{multline}
which is equivalent to a differential equation satisfied by $\Li_{p,\alpha}^{\dagger}$ with  $\Li_{p,\alpha}^{\dagger}(0)=1$ (\cite{I-1}, Proposition 2.1.3) and which characterizes the Frobenius.

\subsection{Duals of some usual operations on $\Pi_{1,0}$}

We discuss the combinatorics of some usual operations which appeared above, in particular with respect to the depth filtration.

\begin{Definition} For any $\xi \in \mu_{N}(K)$, let $\tilde{\Pi}_{\xi,0}$ be the subscheme of $\Pi_{\xi,0}$ defined by the equations $f[e_{0}] = f[e_{\xi}] = 0$. 
\end{Definition}

By the shuffle equation, the points of $\tilde{\Pi}_{\xi,0}$ satisfy more generally $f[e_{0}^{n}] = f[e_{\xi}^{n}] = 0$ for any $n>0$. It follows from the definitions that $\tilde{\Pi}_{\xi,0}$ is a sub-group scheme of $\Pi_{\xi,0}$ for the usual group scheme structure on $\Pi_{\xi,0}$ and that $\tilde{\Pi}_{\xi,0}$ is the image of $\tilde{\Pi}_{1,0}$ by the automorphism $(x \mapsto \xi x)_{\ast}$ of $\pi_{1}^{\un,\DR}(\mathbb{P}^{1} \setminus \{0,\mu_{N},\infty\})$. 
\newline\indent We have $\Phi^{(\xi)}_{p,\alpha} \in \tilde{\Pi}_{\xi,0}(K)$ (\cite{U4}, equation (4.1.3) and Proposition 4.3.1 in the $\alpha=-1$ case ; the same proof works for any $\alpha$).
\newline\indent For any $\xi \in \mu_{N}(K)$, one has the implication $f e_{\xi} = e_{\xi}f \Rightarrow f \in  K \langle \langle e_{\xi} \rangle\rangle$, for $f \in K \langle \langle e_{0 \cup \mu_{N}} \rangle\rangle$ ; it follows that $\Ad(e_{\xi})$ restricted to $\tilde{\Pi}_{\xi,0}(K)$ is injective.

\begin{Definition} \label{def subword} Let $w$ be a word on $e_{0 \cup \mu_{N}}$.
\newline (i) Let $\SubWd(w)$ be the set of subwords of $w$ that contain all the letters of $w$ that are not $e_{0}$.
\newline (ii) Let $sw \in \SubWd(w)$. A connected partition $(sw_{j})_{j \in J}$ of $sw$ is a partition of $sw$, viewed as the set of its letters, in subwords as $sw = \amalg_{j\in J} sw_{j}$, such that the letters of each $sw_{j}$ are consecutive in $sw$ (we will say that each $sw_{j}$ is connected in $sw$), and such that at least one letter of each $sw_{j}$ is not $e_{0}$.
\newline (iii) We say that a subword $sw \in \SubWd(w)$ is maximally at the left of $w$ if it contains the  first letter different from $e_{0}$ in $w$  (where words over $e_{0 \cup \mu_{N}}$ are read from the right to the left).
\newline (iv) \label{def coloring}A coloring of a connected partition $(sw_{i})_{i \in I}$ of an element $sw$ of $\SubWd(w)$ is a map $I \mapsto \mu_{N}(K)$, which we will denote by $i \mapsto \xi_{j(i)}$.
\newline (v) Let $sw \in \SubWd(w)$, $(sw_{i})_{i \in I}$ a connected partition of $sw$ and $C=(\xi_{j(i)})_{i\in I}$ be a coloring of $(sw_{i})_{i \in I}$. We call the quotient of $w$ by the partitioned subword $sw = \amalg_{i \in I} sw_{i}$ colored in $C$ the word obtained by replacing, in $w$, each subword $sw_{i}$ by the letter $e_{\xi_{j(i)}}$ ; we denote it by $\frac{w}{((sw_{i})_{i \in I},C)}$.
\end{Definition}

Let $\Wd(e_{0 \cup \mu_{N}})$ be the set of words on $e_{0 \cup \mu_{N}}$. For any non-negative integers $n,d$, let $\Wd_{n}(e_{0 \cup \mu_{N}})$, resp. $\Wd_{\ast,d}(e_{0 \cup \mu_{N}})$, resp. 
$\Wd_{n,d}(e_{0 \cup \mu_{N}}) = \Wd_{n}(e_{0 \cup \mu_{N}}) \cap \Wd_{\ast,d}(e_{0 \cup \mu_{N}})$ be the subset of $\Wd(e_{0 \cup \mu_{N}})$ consisting of the words of weight $n$, resp. of depth $d$, resp. of weight $n$ and depth $d$, on $e_{0 \cup \mu_{N}}$ ; let $\mathcal{O}^{\sh}_{n}$, $\mathcal{O}^{\sh}_{\ast,d}$, $\mathcal{O}^{\sh}_{n,d}$ be the vector subspaces of $\mathcal{O}^{\sh}$ generated respectively by these sets, and let ${}_{\mathbb{Z}} \mathcal{O}^{\sh}_{n}$, 
${}_{\mathbb{Z}} \mathcal{O}^{\sh}_{\ast,d}$, 
${}_{\mathbb{Z}} \mathcal{O}^{\sh}_{n,d}$ be their restriction of scalars to $\mathbb{Z}$.

\begin{Proposition} \label{prop dual composition} \label{cor dual ihara} Let $(n,d) \in \mathbb{N}^{2}$ with $d \leqslant n$.
\newline (i) For any $\xi \in \mu_{N}(K)$, the dual of $\circ^{\smallint_{\xi,0}}$ restricted to $\tilde{\Pi}_{\xi,0}$, resp. of  $\circ_{\Ad}^{\smallint_{\xi,0}}$ restricted to 
$\Ad_{\tilde{\Pi}_{\xi,0}}(e_{\xi})$,
sends :
${}_{\mathbb{Z}} \mathcal{O}^{\sh}_{n,d}
\rightarrow \underset{\substack{n_{1}+n_{2}-1=n \\ d_{1}+d_{2}-1 = d}}{\bigoplus} {}_{\mathbb{Z}} \mathcal{O}^{\sh}_{n_{1},d_{1}} \otimes {}_{\mathbb{Z}} \mathcal{O}^{\sh}_{n_{2},d_{2}}$.
\newline (ii) For any $\xi \in \mu_{N}(K)$, the dual of the map of inversion for the product $\circ^{\smallint_{\xi,0}}$, resp. $\circ_{\Ad}^{\smallint_{\xi,0}}$, sends ${}_{\mathbb{Z}} \mathcal{O}^{\sh}_{n,d} \mapsto {}_{\mathbb{Z}} \mathcal{O}^{\sh}_{n,d}$.
\newline (iii) The action $\circ^{\smallint_{0,0}}$, resp. $\circ_{\Ad}^{\smallint_{0,0}}$, restricted to an action of $\tilde{\Pi}_{1,0}$, resp. of $\Ad_{\tilde{\Pi}_{1,0}}(e_{1})$, sends 
${}_{\mathbb{Z}} \mathcal{O}^{\sh}_{n,d}
\rightarrow \underset{\substack{n_{1}+n_{2}-1=n \\ d_{1}+d_{2}-1 = d}}{\bigoplus} \mathbb{Z}.{}_{\mathbb{Z}} \mathcal{O}^{\sh}_{n_{1},d_{1}} \otimes {}_{\mathbb{Z}} \mathcal{O}^{\sh}_{n_{2},d_{2}}$.
\end{Proposition}

\begin{proof} The result follows from the facts below and the formulas of equation (\ref{eq:Ihara}), (\ref{eq:Ihara bis}) and of Definition \ref{def Ad action 1,0}.
\newline\indent (a) the product $(g,f) \mapsto fg$, whose dual sends ${}_{\mathbb{Z}} \mathcal{O}^{\sh}_{n,d} \mapsto \underset{\substack{n_{1}+n_{2}=n\\ d_{1}+d_{2}=d}}{\bigoplus} {}_{\mathbb{Z}} \mathcal{O}^{\sh}_{n_{1},d_{1}} \otimes {}_{\mathbb{Z}} \mathcal{O}^{\sh}_{n_{2},d_{2}}$.
\newline\indent (b) The isomorphism $\Pi_{1,0} \mapsto \Pi_{\xi,0}$ $f \mapsto f^{(\xi)}$, whose dual sends ${}_{\mathbb{Z}} \mathcal{O}^{\sh}_{n,d} \mapsto {}_{\mathbb{Z}} \mathcal{O}^{\sh}_{n,d}$ and commutes with the Ihara product and the adjoint action.
\newline\indent (c) The composition of non-commutative formal power series ; let $f$ in $K \langle\langle e_{0 \cup \mu_{N}} \rangle\rangle$, and $(h_{\xi})_{\xi \in \mu_{N}(K)}$ in $K \langle\langle e_{0 \cup \mu_{N}} \rangle\rangle^{N}$, such that $h_{\xi}[\emptyset] =0$, and $f[e_{0}^{n}] = h_{\xi}[e_{0}^{n}] =0$ for any $\xi \in \mu_{N}(K)$ and any positive integer $n$ ; we have 
\begin{multline*} f \big(e_{0},(h_{\xi})_{\xi \in \mu_{N}(K)} \big) = f[\emptyset]+
\\ \displaystyle \sum_{\substack{d \in \mathbb{N}^{\ast}
\\ (n_{0},\ldots,n_{d}) \in (\mathbb{N}^{\ast})^{d+1}
\\ \xi_{1},\ldots,\xi_{d} \in \mu_{N}(K)^{d}}}  f[e_{0}^{n_{d}-1}e_{\xi_{d}} \ldots e_{\xi_{1}}e_{0}^{n_{0}-1}]
e_{0}^{n_{d}-1}
\bigg( \sum_{w_{d} \in \Wd(e_{0 \cup \mu_{N}})} h_{\xi_{d}}[w_{d}]w_{d} \bigg) 
\ldots
\bigg( \sum_{w_{1} \in \Wd(e_{0 \cup \mu_{N}})} h_{\xi_{1}}[w_{1}]w_{1} \bigg)
e_{0}^{n_{0}-1} 
\end{multline*}
and rewriting the right-hand side in that equation as a sum indexed by the words on $e_{0 \cup \mu_{N}}$, we obtain, for any word $w$,
\begin{equation} \label{eq:part c of proof} f(e_{0},(h_{\xi})_{\xi \in \mu_{N}(K)})[w]= f[\emptyset] + \sum_{sw \in \SubWd(w)}
\sum_{\substack{(sw_{i})_{i\in I} \\ \text{connected } \\ \text{ partition} \\ \text{ of }sw}} \sum_{\substack{C =(\xi_{j(i)})_{i \in I} \\ \text{ coloring of } \\ (sw_{j})_{i \in I}}}
\Big( \prod_{i \in I} h_{\xi_{i(j)}}[sw_{i}] \Big) f \Big[\frac{w}{((sw_{i})_{i \in I},C)} \Big].
\end{equation}
\noindent We check that, for any $w \in \Wd(e_{0 \cup \mu_{N}})$, $sw \in \SubWd(w)$, $(sw_{i})_{i \in I}$ connected partition of $sw$, and $C$ coloring of $(sw_{i})_{i \in I}$, we have $\depth(\frac{w}{((sw_{i})_{I \in I},C)}) = \depth(w) - \sum\limits_{i \in I} (\depth(w_{i})-1)$, and $\weight(\frac{w}{((sw_{i})_{I \in I},C)}) = \weight(w) - \sum\limits_{i \in I} (\weight(w_{i})-1)$. Let us assume that there exists $g \in \tilde{\Pi}_{1,0}(K)$ such that $h_{\xi} = \Ad_{g^{(\xi)}}(e_{\xi})$ for any $\xi \in \mu_{N}(K)$. Using (b), the shuffle equation for $g^{(\xi)}$, and the fact that the antipode of $\mathcal{O}^{\sh}$ is given by $e_{x_{l}} \ldots e_{x_{1}} \mapsto (-1)^{l}e_{x_{1}}\ldots e_{x_{l}}$ which gives a description of the coefficients of ${g^{(\xi)}}^{-1}$, we deduce that the dual of the map
$(g,f) \mapsto f\big(e_{0},(\Ad_{g^{(\xi)}}(e_{\xi}))_{\xi \in \mu_{N}(K)}\big)$ sends
${}_{\mathbb{Z}} \mathcal{O}^{\sh}_{n,d} \mapsto
\underset{\substack{d',d'',n',n'' \geqslant 0 \\ d'+d''-1 = d \\ n' + n''-1=n}}{\bigoplus} {}_{\mathbb{Z}} \mathcal{O}^{\sh}_{n',d'} \otimes {}_{\mathbb{Z}} \mathcal{O}^{\sh}_{n'',d''}$.
\end{proof}

\subsection{Groupoids in ultrametric complete groups associated with $\pi_{1}^{\un,\DR}(\mathbb{P}^{1} \setminus \{0,\mu_{N},\infty\})$\label{topology}}

The following definitions will enable us to define the $p$-adic pro-unipotent harmonic action of integrals as a continuous action of a complete topological group (Definition \ref{dR-rt harmonic Ihara action}), and will play a central role in \cite{I-3}.

\subsubsection{$K \langle \langle e_{0 \cup \mu_{N}} \rangle\rangle$ as a ultrametric complete normed algebra}

Let us consider formal variables $U_{1},\ldots,U_{\nu}$, where $\nu \in \mathbb{N}^{\ast}$, and let us equip the set $\mathbb{R}_{+}[[U_{1},\ldots,U_{\nu}]]$ with the product topology associated with the real topology on $\mathbb{R}_{+}$ and the natural identification $\mathbb{R}_{+}[[U_{1},\ldots,U_{\nu}]] \simeq \mathbb{R}_{+}^{\mathbb{N}^{\nu}}$. Let us define a partial order on $\mathbb{R}_{+}[[U_{1},\ldots,U_{\nu}]]$ by declaring that 
$\sum\limits_{(n_{1},\ldots,n_{\nu})\in \mathbb{N}^{\nu}} s_{n_{1},\ldots,n_{\nu}} U_{1}^{n_{1}}\ldots,U_{\nu}^{n_{\nu}} \leqslant \sum\limits_{(n_{1},\ldots,n_{\nu})\in \mathbb{N}^{\nu}} s'_{n_{1},\ldots,n_{\nu}} U_{1}^{n_{1}} \ldots U_{\nu}^{n_{\nu}}$ if, for all $(n_{1},\ldots,n_{\nu}) \in \mathbb{N}^{\nu}$, we have $s_{n_{1},\ldots,n_{\nu}} \leqslant s'_{n_{1},\ldots,n_{\nu}}$. If $S \leqslant S'$ in the sense above, then we have $SR \leqslant S'R$ for all $R \in \mathbb{R}_{+}[[U_{1},\ldots,U_{\nu}]]$. 
Let the maximum of two elements $\sum\limits_{(n_{1},\ldots,n_{\nu})\in \mathbb{N}^{m}} s_{n_{1},\ldots,n_{\nu}} U_{1}^{n_{1}}\ldots,U_{\nu}^{n_{\nu}}$ and $\sum\limits_{(n_{1},\ldots,n_{\nu})\in \mathbb{N}^{\nu}} s'_{n_{1},\ldots,n_{\nu}} U_{1}^{n_{1}}\ldots,U_{\nu}^{n_{\nu}}$ be
$\sum\limits_{(n_{1},\ldots,n_{\nu})\in \mathbb{N}^{m}} \max(s_{n_{1},\ldots,n_{\nu}},s'_{n_{1},\ldots,n_{\nu}}) U_{1}^{n_{1}}\ldots,U_{\nu}^{n_{\nu}}$.
\newline\indent Let $\mathcal{C}$ be a $K$-algebra equipped with a map $\mathcal{N} : \mathcal{C} \rightarrow \mathbb{R}_{+}[[U_{1},\ldots,U_{\nu}]]$ satisfying the axioms of an (ultrametric) algebra norm, adapted to maps having target $ \mathbb{R}_{+}[[U_{1},\ldots,U_{\nu}]]$ with the notions of order (and maximum) on $\mathbb{R}_{+}[[U_{1},\ldots,U_{\nu}]]$ defined above, and satisfying $\mathcal{N}(1_{\mathcal{C}}) = 1$. Then we say that $\mathcal{C}$ is a (ultrametric) normed $K$-algebra with norm $\mathcal{N}$. Any (ultrametric) normed $K$-algebra in the sense of this definition is in particular a (ultra)metric space where the distance is defined by the norm.
 
\begin{Definition} (i) \label{def bounded} Let $K \langle \langle e_{0 \cup \mu_{N}} \rangle \rangle_{<\infty}$, be the subset of
$K \langle \langle e_{0 \cup \mu_{N}} \rangle \rangle$ of the elements $f$ such that, for each $d \in \mathbb{N}^{\ast}$, we have
$\underset{w \in \Wd_{\ast,d}(e_{0 \cup \mu_{N}})}{\sup} |f[w]|_{p} < \infty$. We say that the elements of $K \langle \langle e_{0 \cup \mu_{N}} \rangle \rangle_{<\infty}$ are the bounded elements of $K \langle \langle e_{0 \cup \mu_{N}} \rangle \rangle$. 
\newline (ii) \label{def summable}Let $K \langle \langle e_{0 \cup \mu_{N}} \rangle \rangle_{o(1)}$, be the subset of $K \langle \langle e_{0 \cup \mu_{N}} \rangle \rangle$ consisting of elements $f$ such that, for all $d \in \mathbb{N}^{\ast}$, we have : $\underset{w \in \Wd_{n,d}(e_{0 \cup \mu_{N}})}{\sup} \big|f[w] \big|_{p}  \underset{n\rightarrow \infty}{\longrightarrow} 0$, i.e. $\sum\limits_{l \in \mathbb{N}} |f[w_{l}]|_{p} < +\infty$ for all sequences $(w_{l})_{l \in\mathbb{N}}$ of words over $e_{0 \cup \mu_{N}}$ such that $\displaystyle\weight(w_{l})\underset{l \rightarrow \infty}{\longrightarrow} \infty$ and $\displaystyle\underset{l \rightarrow \infty}{\limsup}\depth(w_{l}) < \infty$. We say that the elements of $K \langle \langle e_{0 \cup \mu_{N}} \rangle \rangle_{o(1)}$ are the summable elements of $K \langle \langle e_{0 \cup \mu_{N}} \rangle \rangle$.
\end{Definition}

\begin{Definition} \label{def norms} (i)
$\mathcal{N}_{\Lambda} : K \langle \langle e_{0 \cup \mu_{N}} \rangle \rangle \longrightarrow \mathbb{R}_{+} [[\Lambda]]$, $ f \mapsto \sum\limits_{n \in \mathbb{N}} 
\underset{w \in\Wd_{n}(e_{0 \cup \mu_{N}})}{\max} \big|f[w]\big|_{p} \Lambda^{n}$.
\newline (ii) Let $\mathcal{N}_{\Lambda,D} : K \langle \langle e_{0 \cup \mu_{N}} \rangle \rangle \rightarrow \mathbb{R}_{+}[[\Lambda,D]]$, 
$ f \mapsto 
\sum\limits_{(n,d) \in \mathbb{N}^{2}} 
\underset{w \in\Wd_{n,d}(e_{0 \cup \mu_{N}})}{\max} \big|f[w]\big|_{p} \Lambda^{n}D^{d} $.
\newline(iii) Let $\mathcal{N}_{D} : K \langle \langle e_{0 \cup \mu_{N}} \rangle \rangle_{<\infty} \rightarrow \mathbb{R}_{+}[[D]]$, 
$ f \mapsto \sum\limits_{d \in \mathbb{N}} 
\bigg( \underset{w \in \Wd_{\ast,d}(e_{0 \cup \mu_{N}})}{\sup}|f[w]|_{p} \bigg) D^{d}$.
\end{Definition}

The topology induced by $\mathcal{N}_{\Lambda,D}$ resp. $\mathcal{N}_{\Lambda}$ is the topology of pointwise convergence on $K \langle \langle e_{0 \cup \mu_{N}} \rangle \rangle$ viewed as the set of maps $\Wd(e_{0 \cup \mu_{N}})\rightarrow K$, and we will use only $\mathcal{N}_{\Lambda,D}$ in the rest of this text. The topology induced by $\mathcal{N}_{D}$ is the topology on $K \langle\langle e_{0 \cup \mu_{N}} \rangle\rangle_{<\infty}$, viewed as a set of maps $\Wd(e_{0 \cup \mu_{N}})\rightarrow K$, of uniform convergence on all the subsets $\Wd_{\ast,d}(e_{0 \cup \mu_{N}})$, $d \in \mathbb{N}^{\ast}$, i.e. the topology of uniform convergence in bounded depth. The topology defined by $\mathcal{N}_{D}$ will be natural when we deal with the sums of $p$-adic series arising from the study of $p$CMZVs ; our computation of $p$CMZVs will be compatible with the depth filtration.
\newline\indent $\mathcal{N}_{D}$ and $\mathcal{N}_{\Lambda}$ can be factorized by $\mathcal{N}_{\Lambda,D}$ from which it follows the implications  $\mathcal{N}_{\Lambda,D}(f) \leqslant \mathcal{N}_{\Lambda,D}(g) \Rightarrow \mathcal{N}_{\Lambda}(f) \leqslant \mathcal{N}_{\Lambda}(g)$ and $\mathcal{N}_{\Lambda,D}(f) \leqslant \mathcal{N}_{\Lambda,D}(g) \Rightarrow \mathcal{N}_{D}(f) \leqslant \mathcal{N}_{D}(f)$, which prove that $\mathcal{N}_{\Lambda}$ and $\mathcal{N}_{D}$ inherit of most of the properties of $\mathcal{N}_{\Lambda,D}$.

\begin{Proposition} \label{axioms} (i) $K \langle \langle e_{0 \cup \mu_{N}} \rangle \rangle$ equipped with $\mathcal{N}_{\Lambda,D}$ is an ultrametric complete normed $K$-algebra.
\newline (ii) $K\langle\langle e_{0 \cup \mu_{N}} \rangle\rangle_{<\infty}$ and $K\langle\langle e_{0 \cup \mu_{N}} \rangle\rangle_{o(1)}$ equipped with $\mathcal{N}_{D}$ are complete ultrametric normed $K$-algebras.
\end{Proposition}

\begin{proof} (i) It is clear that $\mathcal{N}_{\Lambda,D}$ satisfies the separation and homogeneity properties of norms ; moreover, for any $f,g \in K \langle\langle e_{0 \cup \mu_{N}} \rangle \rangle$ we have 
$\mathcal{N}_{\Lambda,D}(f+g)\leqslant \max(\mathcal{N}_{\Lambda,D}(f),\mathcal{N}_{\Lambda,D}(g))$ and $\mathcal{N}_{\Lambda,D}(gf) \leqslant \mathcal{N}_{\Lambda,D}(g)\mathcal{N}_{\Lambda,D}(f)$ : the first inequality is clear and the second is obtained by writing, for any word $w \in \Wd_{n,d}(e_{0\cup \mu_{N}})$,  $|(gf)[w]|_{p} = \big|\sum\limits_{w_{1}w_{2}=w} g[w_{1}]f[w_{2}]\big|_{p} \leqslant 
\sum\limits_{w_{1}w_{2}=w} |g[w_{1}]|_{p}|f[w_{2}]|_{p} \leqslant \sum\limits_{\substack{n_{1}+n_{2}=n\\ d_{1}+d_{2}=d}} \underset{w_{1} \in \Wd_{n_{1},d_{1}}(e_{0\cup \mu_{N}})}{\sup |f[w_{1}]|_{p}} \underset{w_{2} \in \Wd_{n_{2},d_{2}}(e_{0\cup \mu_{N}})}{\sup |g[w_{2}]|_{p}}$. This proves that $K \langle\langle e_{0 \cup \mu_{N}} \rangle \rangle$ equipped with $\mathcal{N}_{\Lambda,D}$ is a normed $K$-algebra. Its completeness follows from the fact that $K$ is complete.
\newline (ii) Let
$\mathbb{R}_{+}[[\Lambda,D]]_{<\infty}$ resp.  $\mathbb{R}_{+}[[\Lambda,D]]_{o(1)}$ be the set of elements $S=\sum\limits_{n,d=0}^{\infty} s_{n,d}D^{d}\Lambda^{n}$ such that, for all $d$, $\underset{n\geqslant 0}{\sup}\text{ }s_{n,d} < \infty$, resp. $s_{n,d} \underset{n\rightarrow \infty}{\rightarrow} 0$.
One can check easily that if $S$ and $S'$ are in
$\mathbb{R}_{+}[[\Lambda,D]]_{<\infty}$ 
resp. $\mathbb{R}_{+}[[\Lambda,D]]_{o(1)}$, then  $\max(S,S')$ and $S \times S'$ satisfy the same property. This shows that $K \langle \langle e_{0 \cup \mu_{N}} \rangle\rangle_{<\infty}$ and $K \langle \langle e_{0 \cup \mu_{N}} \rangle\rangle_{o(1)}$ are subalgebras of $K \langle \langle e_{0 \cup \mu_{N}} \rangle\rangle$. The axioms of an algebra norm for $\mathcal{N}_{D}$ on $K \langle\langle e_{0 \cup \mu_{N}} \rangle \rangle_{<\infty}$ are checked as in (i). Thus $K \langle \langle e_{0 \cup \mu_{N}} \rangle\rangle_{<\infty}$ and $K \langle \langle e_{0 \cup \mu_{N}} \rangle\rangle_{o(1)}$ are normed algebras with $\mathcal{N}_{D}$. Their completeness follows from the fact that the spaces of sequences $\ell^{\infty}(K)$ and $c_{0}(K)$ equipped with the norm $||.||_{\infty}$ are complete.
\end{proof}

\subsubsection{Groupoids in complete ultrametric groups associated with $\pi_{1}^{\un,\DR}(\mathbb{P}^{1} \setminus \{0,\mu_{N},\infty\}(K)$}

From now on, $x,y$ are any two base-points of $\pi_{1}^{\un,\DR}(\mathbb{P}^{1} \setminus \{0,\mu_{N},\infty\})$ ; $\Pi_{y,x}(K)$ is identified to the subset of $K\langle\langle e_{0\cup \mu_{N}}\rangle\rangle$ of elements satisfying the shuffle equation and having constant coefficient equal to $1$ (\S1.1.1).

\begin{Definition} \label{summable points}Let $\Pi_{y,x}(K)_{< \infty} = \Pi_{y,x}(K) \cap K \langle\langle e_{0 \cup \mu_{N}}\rangle\rangle_{<\infty}$ and $\Pi_{y,x}(K)_{o(1)} = \Pi_{y,x}(K) \cap K \langle\langle e_{0 \cup \mu_{N}}\rangle\rangle_{o(1)}$.
\end{Definition}

\begin{Proposition} \label{corollary le corollary}
(i) $\Pi_{y,x}(K)$ is a complete topological group for the $\mathcal{N}_{\Lambda,D}$-topology ; $\Pi_{y,x}(K)_{<\infty}$ and $\Pi_{y,x}(K)_{o(1)}$ are complete topological groups for the $\mathcal{N}_{D}$-topology.
\newline (ii) The groupoid law on $\pi_{1}^{\un,\DR}(\mathbb{P}^{1} \setminus \{0,\mu_{N},\infty\}(K)$, resp. on its subgroups of bounded, resp. summable points is continuous for the $\mathcal{N}_{\Lambda,D}$-topology, resp. for the $\mathcal{N}_{D}$-topology.
\end{Proposition}

\begin{proof} (i) We know that $\Pi_{y,x}(K)$ resp. $K \langle\langle e_{0 \cup \mu_{N}}\rangle\rangle_{<\infty}$ and $K \langle\langle e_{0 \cup \mu_{N}}\rangle\rangle_{o(1)}$ is resp. are stable by multiplication (respectively by \S1.1.1 and Proposition \ref{axioms}), so $\Pi_{y,x}(K)_{<\infty}$ and $\Pi_{y,x}(K)_{o(1)}$ are stable by multiplication. 
\newline\indent For $f \in \Pi_{y,x}(K)$, and $l \in \mathbb{N}^{\ast}$, we have  $\mathcal{N}_{\Lambda,D}(f^{l}) = \mathcal{N}_{\Lambda,D}(f)$ ; indeed, this amounts to, for all $n,d \in \mathbb{N}^{\ast}$, $\underset{\substack{(w_{1},\ldots,w_{l})\in \Wd(e_{0 \cup \mu_{N}})^{l} \\ \text{s.t. }w_{1}\ldots w_{l} \in \Wd_{n,d}(e_{0 \cup \mu_{N}})}}{\max} \big| \prod_{i=1}^{l} f[w_{i}] \big|_{p} = \underset{w \in \Wd_{n,d}(e_{0 \cup \mu_{N}})}{\max} \big|f[w]\big|_{p}$. The inequality $\geqslant$ is obtained by choosing $w_{2} = \ldots = w_{l} = \emptyset$ in the left-hand side since $f[\emptyset] = 1$ ; the inequality $\leqslant$ follows from the shuffle equation for $f$ and from that the shuffle product restricts, for all $n_{1},n_{2},d_{1},d_{2} \in \mathbb{N}$, to a map ${}_{\mathbb{Z}}\mathcal{O}_{n_{1},d_{1}}^{\sh} \times {}_{\mathbb{Z}}\mathcal{O}_{n_{2},d_{2}}^{\sh} \rightarrow {}_{\mathbb{Z}}\mathcal{O}_{n_{1}+n_{2},d_{1}+d_{2}}^{\sh}$.
\newline\indent Now, for $f \in \Pi_{y,x}(K)$, we have $f^{-1} = \sum_{l \in \mathbb{N}} (1-f)^{l}$, where for each $w \in \mathcal{O}^{\sh,e_{0 \cup \mu_{N}}}$, the sum $\sum_{l \in \mathbb{N}} (1-f)^{l}[w]$ is finite. In particular, the ultrametric triangle inequality for $\mathcal{N}_{\Lambda,D}$ has a sense and remains true for this infinite sum, and we have $\mathcal{N}_{\Lambda,D}(f^{-1}) \leqslant \max_{l \in \mathbb{N}} \mathcal{N}_{\Lambda,D}((1-f)^{l}) \leqslant \max_{l\in\mathbb{N}} \mathcal{N}_{\Lambda,D}(f^{l}) = \mathcal{N}_{\Lambda,D}(f)$, where the last inequality follows from the binomial expansion of $(1-f)^{l}$ and from the ultrametric triangle inequality for $\mathcal{N}_{\Lambda,D}$. By symmetry of the roles of $f$ and $f^{-1}$, we deduce $\mathcal{N}_{\Lambda,D}(f^{-1}) = \mathcal{N}_{\Lambda,D}(f)$. This implies that $\Pi_{y,x}(K)_{<\infty}$ and $\Pi_{y,x}(K)_{o(1)}$ are stable by inversion.
\newline\indent In particular, $\Pi_{y,x}(K)_{<\infty}$ and $\Pi_{y,x}(K)_{o(1)}$ are subgroups of $\Pi_{y,x}(K)$. On the other hand, they are defined by the shuffle equation and $f[\emptyset]=1$ so they are closed subsets respectively of $K\langle \langle e_{0\cup \mu_{N}}\rangle\rangle_{<\infty}$ and $K\langle \langle e_{0\cup \mu_{N}}\rangle\rangle_{o(1)}$ which are complete by Proposition \ref{axioms}, so they are complete.
\newline\indent (ii) follows from (i) and the fact that the canonical isomorphisms $\Pi_{y,x} \simeq \Pi$ reviewed in \S1.1.1 are compatible with the groupoid structure.
\end{proof}

As a conclusion, we have two groupoids in ultrametric complete groups, defined respectively by the bounded and the summable points of $\pi_{1}^{\un,\DR}(\mathbb{P}^{1} \setminus \{0,\mu_{N},\infty\}(K)$, whose topology is defined by the uniform convergence in bounded depth.

\subsubsection{Compatibility with the byproducts of the motivic Galois action}

We prove that the groupoids constructed in \S1.3.2 are stable by the usual operations of $\pi_{1}^{\un,\DR}(\mathbb{P}^{1} \setminus \{0,\mu_{N},\infty\}(K)$ related to the motivic Galois action, and that these operations are continuous. We use Definition 1.2.1.

\begin{Proposition} \label{norm adjoint action}\label{prop norm ihara bis}\label{prop summable group ihara}  
(i) $(\tilde{\Pi}_{\xi,0},\circ^{\smallint_{\xi,0}})$, resp. 
$(\tilde{\Pi}_{\xi,0}(K)_{o(1)},\circ^{\smallint_{\xi,0}})$, $(\tilde{\Pi}_{\xi,0}(K)_{<\infty},\circ^{\smallint_{\xi,0}})$
are complete topological groups for the $\mathcal{N}_{\Lambda,D}$-topology, resp. the $\mathcal{N}_{D}$-topology ; $\Ad(e_{\xi})$ induces isomorphisms of complete topological groups between them and their images.
\newline (ii) $\circ^{\smallint_{0,0}} : \tilde{\Pi}_{1,0}(K) \times \Pi_{0,0}(K) \mapsto \Pi_{0,0}(K)$ is a continuous group action for the $\mathcal{N}_{D}$-topology ; $\Ad(e_{1})$ induces an isomorphism of continuous group actions between $\circ^{\smallint_{0,0}}$ restricted to $\tilde{\Pi}_{1,0}(K)$ and $\circ_{\Ad}^{\smallint_{0,0}}$ restricted to its image.
\newline (iii) The map $\tau : K^{\times} \times K\langle\langle e_{0\cup \mu_{N}}\rangle\rangle \rightarrow K\langle\langle e_{0\cup \mu_{N}}\rangle\rangle$, resp.  $\tau : \{\lambda \in K^{\ast} \text{ }|\text{ }|\lambda_{p}\leqslant 1\}\times  K\langle\langle e_{0\cup \mu_{N}}\rangle\rangle \times K\langle\langle e_{0\cup \mu_{N}}\rangle\rangle_{<\infty} \rightarrow K\langle\langle e_{0\cup \mu_{N}}\rangle\rangle_{<\infty}$ is continuous for the $\mathcal{N}_{\Lambda,D}$-topology resp. it is continuous for the $\mathcal{N}_{D}$-topology and stabilizes the groups of (i).
\end{Proposition}

\begin{proof} (i) 
It follows from Proposition \ref{prop dual composition} that we have 
for all $f,g \in \tilde{\Pi}_{1,0}(K)$, $\mathcal{N}_{\Lambda,D}( g \circ^{\smallint_{1,0}}f) \leqslant \mathcal{N}_{\Lambda,D}(g) \times \mathcal{N}_{\Lambda,D}(f)$ and
$\mathcal{N}_{\Lambda,D}(f^{-1_{\circ^{\smallint_{1,0}}}}) = \mathcal{N}_{\Lambda,D}(f)$, and similarly with $\mathcal{N}_{D}$ if $f$ and $g$ are bounded.
This proves $\tilde{\Pi}_{1,0}(K)_{<\infty}$, $\tilde{\Pi}_{1,0}(K)_{o(1)}$ are subgroups of $\Pi_{1,0}(K)$ for $\circ^{\smallint_{1,0}}$. On the other hand, by Proposition \ref{norm adjoint action},  $\tilde{\Pi}_{1,0}(K)_{<\infty}$, $\tilde{\Pi}_{1,0}(K)_{o(1)}$ are complete.
\newline\indent By the shuffle equation for $f$, we have $\mathcal{N}_{\Lambda,D}(\Ad_{f}(e_{\xi})) \leqslant \Lambda D \mathcal{N}_{\Lambda,D}(f)$, and $\mathcal{N}_{D}(\Ad_{f}(e_{\xi})) \leqslant D \mathcal{N}_{D}(f)$, whence $\Ad(e_{\xi})$ is continuous and we have $\Ad(e_{\xi})(\Pi_{y,x}(K)_{<\infty}) \subset K\langle\langle e_{0 \cup \mu_{N}}\rangle\rangle_{<\infty}$
and $\Ad(e_{\xi})(\Pi_{y,x}(K)_{o(1)}) \subset K\langle\langle e_{0 \cup \mu_{N}}\rangle\rangle_{o(1)}$. Moreover, we have seen that $\Ad(e_{\xi})$ restricted to $\tilde{\Pi}_{\xi,0}(K)$ is injective. The isomorphisms $(x \mapsto \xi x)_{\ast} : \Pi_{1,0}\simeq \Pi_{\xi,0}$ are homeomorphisms both for the $\mathcal{N}_{\Lambda,D}$-topology and the $\mathcal{N}_{D}$-topology.
\newline\indent (ii) is proved like (i) and (iii) is immediate by, for all $f \in \Pi_{y,x}(K)$, $\mathcal{N}_{\Lambda,D}(\tau(\lambda)(f))(\Lambda,D)
= \mathcal{N}_{\Lambda,D}(f)(\lambda\Lambda,D)$.
\end{proof}

In particular, by Proposition \ref{norm adjoint action}, combined to equation (\ref{eq: Frobenius at 10}) and (\ref{eq:first}), we deduce that the Frobenius $\phi$ is compatible with these topological structures (i.e. it is continuous and stabilizes the subgroupoids of bounded and summable elements, in the above sense). Our computation of the Frobenius will be also compatible with this structure.

\section{The pro-unipotent harmonic action of integrals}

We observe a simplification in the differential equation of the Frobenius (\S3.1), and we introduce the pro-unipotent harmonic action of integrals (\S3.2) which enables to express the simplified differential equation of the Frobenius and prove the "integral" part of the theorem (\S3.3).

\subsection{A simplification of the equation of the Frobenius}

\subsubsection{Suppressing a parameter in the equation of the Frobenius}

We reformulate the differential equation of the Frobenius (equation (\ref{eq:first})) in terms of the coefficients of its power series expansion at $0$.
 
\begin{Notation}
If $S\in K[[z]]$ and $m\in\mathbb{N}$, we denote the coefficient of $z^{m}$ in $S$ by $S[z^{m}]$.
\end{Notation}

\noindent In the next statement, the significative feature is that the right-hand side of the equation does not depend on $l$, whereas the left-hand side is a priori a complicated function of $l$.

\begin{Proposition} \label{previous formula} Let $d$ and $n_{i}$ $(1 \leqslant i \leqslant d)$ be positive integers, and let $\xi_{i}$ be $N$-th roots of unity $(1 \leqslant i \leqslant d)$. Let $w = \big( (n_{i})_{d};(\xi_{i})_{d+1}\big)$,
$w^{(p^{\alpha})} = \big( (n_{i})_{d};(\xi_{i}^{p^{\alpha}})_{d+1}\big)$, and $w_{l} = e_{0}^{l-1}e_{\xi_{d+1}} e_{0}^{n_{d}-1}e_{\xi_{d}} \ldots e_{0}^{n_{1}-1}e_{\xi_{1}}$ for all $l \in \mathbb{N}^{\ast}$. We have, for all $m \in \mathbb{N}^{\ast}$,
\begin{multline}
\label{eq:eq previous formula}
\tau(m)\bigg[ \Li_{p,X^{(p^{\alpha})}}^{\KZ}(z)\big(e_{0},
(\Ad_{(\Phi^{(\xi)}_{p,\alpha})}(e_{\xi}))_{\xi \in \mu_{N}(K)}\big)[w_{l}][z^{m}] +
 \Li^{\dagger}_{p,\alpha}[w_{l}][z^{p^{\alpha}m}] +
\\ \sum_{\substack{ \{ 1 \leqslant b \leqslant p^{\alpha}m-1 \text{ }|\text{ } p^{\alpha}|b \} \\ \{ (u_{l},v) \text{ }|  w_{l} = u_{l}v, \text{ }\depth(u_{l}) \geqslant 1 , \text{ }v \neq \emptyset \}}}
\Li_{p,\alpha}[u_{l}][z^{b}]\text{ }\text{ }\Li_{p,X^{(p^{\alpha})}}^{\KZ}(z^{p^{\alpha}}) \big(e_{0},\Ad_{(\Phi^{(\xi)}_{p,\alpha})}(e_{\xi})\big) [v_{r}][z^{m-\frac{b}{p^{\alpha}}}\log(z)^{0}]\bigg] \\ = (-1)^{d+1} \har_{p^{\alpha}m}(w).
\end{multline}
\end{Proposition}

\begin{proof} In equation (\ref{eq:first}), we take the coefficient of $w_{l}$, then, the coefficient of $z^{p^{\alpha}m}$ in the series expansion at $0$ with respect to $z$, then we apply $\tau(m)$ (equation (\ref{eq:tau})).
\newline\indent (a) By the definition of $p$-adic multiple polylogarithms in terms of the KZ equation \cite{Furusho 1} \cite{Yamashita}, we have, for all $m' \in \mathbb{N}^{\ast}$ :
$\tau(m')\Li_{p,X}^{\KZ}[w_{l}][z^{m'}] = (-1)^{d+1} \har_{m'} \big( w \big)$. This gives an expression of the right-hand side.
\newline\indent (b) The left-hand side, which is defined by a product involving $\Li_{p,\alpha}^{\dagger}(z)$ and $\Li_{p,X^{(p^{\alpha})}}^{\KZ}(z^{p^{\alpha}})$, is a sum over $b$ in the set $\{0,\ldots,p^{\alpha}m \}$, and over couples $(u_{l},v)$ such that $w_{l} = u_{l}v$. By
\label{lemma vanishing Li} \cite{I-1}, Lemma 4.2.1, we have, for all $n \in \mathbb{N}^{\ast}$, $\Li_{p,\alpha}^{\dagger}[e_{0}^{n}](z) = 0$. Moreover, by the definitions, we have $\Li_{p,\alpha}^{\dagger}[w][z^{0}] = \Li_{p,X^{(p^{\alpha})}}^{\KZ}[w^{(p^{\alpha})}][z^{0}] = \Li_{p,X}^{\KZ}[w][z^{0}] = 1$.
Thus the sum over $b$ can be restricted to $\alpha \in \{1,\ldots,p^{\alpha}m-1\}$, and the sum over $u_{l}$ can be restricted to terms such that $\depth(u_{l}) \geqslant 1$. The sum over $b$ can be reindexed by $\tilde{b}=\frac{b}{p^{\alpha}} \in \mathbb{N}$, since for any power series $S$, we have $S(z^{p^{\alpha}})[z^{p^{\alpha}\tilde{b}}] = S(z)[z^{\tilde{b}}]$.
\end{proof}

\subsubsection{Vanishing of a certain limit of the terms having an overconvergent factor}

We are going to exploit the suprising observation in Proposition \ref{previous formula} by computing the limit of the left-hand side of equation (\ref{eq:eq previous formula}) when $l \rightarrow \infty$. We are going to show that this limit is particularly simple. This is going to follow from the main result of \cite{I-1}.
\newline\indent Heuristically, since the functions $\Li_{p,\alpha}^{\dagger}[w']$ are overconvergent whereas the power series expansion at $0$ of the functions $\Li_{p}^{\KZ}[w'']$ converge only on $\{z \in K \text{ }|\text{ }|z|_{p}<1 \}$, the coefficients $\Li_{p,\alpha}^{\dagger}[w'][z^{m}]$ have a priori significantly smaller $p$-adic norms that the coefficients $\Li_{p}^{\KZ}[w''][z^{m}]$, at least for $m$ large.
\newline\indent The main result of \cite{I-1} can be reformulated as follows. Let $U^{\an}=(\mathbb{P}^{1,\an} \setminus \underset{\xi \in \mu_{N}(K)}{\cup} B(\xi,1))/K$, and let $A(U^{\an})$ be its $K$-algebra of rigid analytic functions, which is a Banach algebra over $K$. Extending in a natural way the notion of summable $K$-points of $\pi_{1}^{\un,\DR}(\mathbb{P}^{1} \setminus \{0,\mu_{N},\infty\})$ (Definition \ref{summable points}) to points having coefficients in any Banach algebra over $K$, the main theorem of \cite{I-1} is
$\Li_{p,\alpha}^{\dagger} \in \Pi_{0,0}(A(U^{\an}))_{o(1)}$, where, implicitly, we view the $\nabla_{\KZ}$ as a connection on a bundle trivialized at $\vec{1}_{0}$, as in \cite{I-1}.

\begin{Corollary} \label{cor tend to 0} The term
\begin{equation} \label{eq:tend to 0} \tau(m) \bigg[\Li^{\dagger}_{p,\alpha}[w_{l}][z^{p^{\alpha}m}] +
\\ \sum_{\substack{ \{ 1 \leqslant b \leqslant p^{\alpha}m-1 \text{ }|\text{ } p^{\alpha}|b \} \\ \{ (u_{l},v) \text{ }|  w_{l} = u_{l}v, \text{ }\depth(u_{l}) \geqslant 1 , \text{ }v \neq \emptyset \}}}
\Li_{p,\alpha}^{\dagger}[u_{l}][z^{b}].\Li_{p,X^{(p^{\alpha})}}^{\KZ}(z^{p^{m}})\big(e_{0},(\Ad_{\Phi^{(\xi)}_{p,\alpha}}(e_{\xi}))_{\xi \in \mu_{N}(K)}\big) [v][z^{m-\frac{b}{p^{\alpha}}}] \bigg]
\end{equation}
tends to $0$ when $l \rightarrow \infty$.
\end{Corollary}

\begin{proof} The set of $v$'s in the sum does not depend on $l$ ; thus, the factors depending on $v$ in the second line are contained in a bounded subset of $K$ depending only on 
$\big((n_{i})_{d},(\xi_{i})_{d+1}\big)$. Moreover, each $u_{l}$ is determined by the unique $v$ such that $w_{l}=u_{l}v$, and there are a finite number, bounded independently of $l$, of such $v$'s. Finally, we have $\limsup \depth u_{l} < +\infty$ and $\weight u_{l} \rightarrow +\infty$, and similarly for $w_{l}$. Whence the result by the theorem of \cite{I-1}.
\end{proof}

\subsection{The pro-unipotent harmonic action of integrals}

\subsubsection{Definition}

We now define a variant of $K \langle \langle e_{0 \cup \mu_{N}} \rangle\rangle$ which will contain in a natural way some non-commutative generating series of cyclotomic multiple harmonic sums. For convenience, in the rest of this text, we will restrict to words $w$ whose rightmost letter is an $e_{\xi}$, with $\xi \in \mu_{N}(K)$. This is sufficient for our purposes. The role of the other words will appear in a subsequent paper.

\begin{Definition} \label{harmonic generating series set}(i) Let $K \langle \langle e_{0 \cup \mu_{N}} \rangle\rangle_{\har}^{\smallint} \subset K \langle\langle e_{0 \cup \mu_{N}} \rangle\rangle$ be the vector subspace of the elements $f$ such that, for all words $w$ on $e_{0 \cup \mu_{N}}$, the sequence $(f[e_{0}^{l}w])_{l\in \mathbb{N}}$ is constant and $f[w'e_{0}]=0$ for all words $w'$.
\newline (ii) Let $K \langle\langle e_{0 \cup \mu_{N}} \rangle\rangle^{\lim} \subset K \langle\langle e_{0 \cup \mu_{N}} \rangle\rangle$ be the vector subspace consisting of the elements $f\in K \langle\langle e_{0 \cup \mu_{N}} \rangle\rangle$ such that, for all words $w$ on $e_{0 \cup \mu_{N}}$, the sequence $(f[e_{0}^{l}w])_{l\in \mathbb{N}}$ has a limit in $K$, and $f[w'e_{0}]=0$ for all words $w'$.
\newline (iii) \label{def limit map}Let 
$\lim : K \langle \langle e_{0 \cup \mu_{N}} \rangle\rangle^{\lim} \rightarrow K \langle \langle e_{0 \cup \mu_{N}} \rangle\rangle_{\har}^{\smallint}$ be the map defined by, for all words $w$ over $e_{0\cup\mu_{N}}$,
$$ (\lim f)[w] = \lim_{l\rightarrow \infty} f[e_{0}^{l}w]. $$
\end{Definition}

If $f \in K \langle \langle e_{0 \cup \mu_{N}} \rangle \rangle_{\har}^{\smallint}$, we denote by $f \big( (n_{i});(\xi_{i})\big) = f[e_{\xi_{d+1}} e_{0}^{n_{d}-1}e_{\xi_{d}} \ldots e_{0}^{n_{1}-1}e_{\xi_{1}}]$ for any positive integers $d$ and $n_{i}$ $(1 \leqslant i \leqslant d)$, and for any $N$-th roots of unity $\xi_{i}$ $(1 \leqslant i \leqslant d+1)$. 

\begin{Definition} \label{dR-rt harmonic Ihara action} The $p$-adic pro-unipotent harmonic action of integrals for $\mathbb{P}^{1} \setminus \{0,\mu_{N},\infty\}$ is the map 
$$ \circ^{\smallint}_{\har} :
\begin{array}{c}
\Ad_{\tilde{\Pi}_{1,0}(K)_{o(1)}}(e_{1}) \times
(K \langle\langle e_{0 \cup \mu_{N}} \rangle\rangle_{\har}^{\smallint})^{\mathbb{N}}
 \rightarrow 
(K \langle\langle e_{0 \cup \mu_{N}} \rangle\rangle_{\har}^{\smallint})^{\mathbb{N}} 
\\ \big( g, (h_{m})_{m\in\mathbb{N}} \big) \mapsto g \circ_{\har}^{\smallint} (h_{m})_{m\in\mathbb{N}} = \big( \lim \big( \tau(m)(g) \circ^{\smallint_{0,0}}_{\Ad} h_{m} \big)\big)_{m\in\mathbb{N}} \end{array}. $$
\end{Definition}

We will prove in Proposition \ref{harmonic Ihara properties} that it is well-defined.

\begin{Examples} \label{example circ har int} For $\mathbb{P}^{1} \setminus \{0,1,\infty\}$ and in depth 1 and 2, for all $n_{1},n_{2} \in \mathbb{N}^{\ast}$, $m \in \mathbb{N}$, for any $g$ and $h=(h_{m})_{m \in \mathbb{N}}$, we have
\begin{equation}
(g \circ^{\smallint}_{\har} h )_{m}(n_{1}) =  f_{m}(n_{1}) + \sum_{l \in \mathbb{N}} m^{n_{1}+l} g[e_{0}^{l}e_{1} e_{0}^{n_{1}-1}e_{1}], 
\end{equation}
\begin{multline}
(g \circ_{\har}^{\smallint} h )_{m}(n_{1},n_{2}) = f_{m}(n_{1},n_{2}) + \sum_{l\in \mathbb{N}} m^{l+n_{1}+n_{2}} g[e_{0}^{l}e_{1}e_{0}^{n_{2}-1}e_{1} e_{0}^{n_{1}-1}e_{1}] +
\\ \sum_{r_{2}=0}^{n_{2}-1} f_{m}(n_{2}-r_{2}) m^{r_{2}+n_{1}} g[e_{0}^{r_{2}}e_{1} e_{0}^{n_{1}-1}e_{1}] 
+ \sum_{r_{1}=0}^{n_{1}-1} f_{m}(n_{1}-r_{1}) \sum_{l \in \mathbb{N}} m^{l+n_{2}+r_{1}} g[e_{0}^{l}e_{1}e_{0}^{n_{2}-1}e_{1} e_{0}^{r_{1}}].
\end{multline}
\end{Examples}

Definition \ref{dR-rt harmonic Ihara action} involves only the summable elements of $\tilde{\Pi}_{1,0}(K)$ ; this restriction is removed below by replacing $m \in \mathbb{N}$ by a formal variable $\underline{m}$ ; we define $K[[\underline{m}]] \langle\langle e_{0 \cup \mu_{N}}\rangle\rangle_{\har}^{\smallint}$ like $K \langle\langle e_{0 \cup \mu_{N}} \rangle\rangle_{\har}^{\smallint}$ (Definition \ref{harmonic generating series set}). We will use most of the time the point of view of Definition \ref{dR-rt harmonic Ihara action} ; the point of view below is practical when we want to consider duals. We need first to write the dual of Definition \ref{harmonic generating series set}.

\begin{Definition} Let $\mathcal{O}_{\har}^{\sh,e_{0 \cup \mu_{N}}}$ be $\mathbb{Q}$-vector space generated by sequences of words of the form
\newline $(e_{0}^{l}e_{\xi_{d+1}}e_{0}^{n_{d}-1}e_{\xi_{d}}\ldots e_{0}^{n_{1}-1}e_{\xi_{1}})_{l \in \mathbb{N}}$.
\end{Definition}

In the next statement, the hats refer to the completions for the weight-adic topology.

\begin{Definition} \label{variant of dr RT} (i) Let the formal pro-unipotent harmonic action of integrals be the map :
$$ \tilde{\circ}^{\smallint}_{\har} : \begin{array}{c}
\Ad_{\tilde{\Pi}_{1,0}(K[m])}(e_{1}) \times
K \langle\langle e_{0 \cup \mu_{N}} \rangle\rangle_{\har}^{\smallint}
\rightarrow 
K[[\underline{m}]] \langle\langle e_{0 \cup \mu_{N}} \rangle\rangle_{\har}^{\smallint} 
\\ \big( g, h \big) \mapsto  \big( \lim \big( \tau(\underline{m})(g) \circ^{\smallint_{0,0}}_{\Ad} h \big) \end{array}. $$
(ii) Let the pro-unipotent harmonic coaction of integrals be the dual of $ \tilde{\circ}^{\smallint}_{\har}$ :
$$ (\tilde{\circ}_{\har}^{\smallint})^{\vee} : \mathcal{O}_{\har}^{\sh,e_{0 \cup \mu_{N}}} 
\longrightarrow \mathcal{O}_{\har}^{\sh,e_{0 \cup \mu_{N}}} \otimes \hat{\tau}(\underline{m})\widehat{\mathcal{O}^{\sh,e_{0 \cup \mu_{N}}}}. $$
(iii) Let the natural factorization, where $T$ denotes the tensor algebra,
$$ (\tilde{\circ}_{\har}^{\smallint})^{\vee} : \mathcal{O}_{\har}^{\sh,e_{0 \cup \mu_{N}}} 
\overset{(\tilde{\circ}_{\har}^{\smallint})^{\vee,T}}{\longrightarrow} \mathcal{O}_{\har}^{\sh,e_{0 \cup \mu_{N}}} \otimes T \big((\hat{\tau}(m)\widehat{\mathcal{O}^{\sh,e_{0 \cup \mu_{N}}}}) \big)^{\otimes N} \longrightarrow \mathcal{O}_{\har}^{\sh,e_{0 \cup \mu_{N}}} \otimes \hat{\tau}(m) \widehat{\mathcal{O}^{\sh,e_{0 \cup \mu_{N}}}} $$
where the $N$ tensor components refer to the coefficients of the $g^{(\xi)}$'s, $\xi \in \mu_{N}(K)$, and the tensor algebra encodes the products of coefficients of each $g^{(\xi)}$.
\end{Definition}

Indeed, the formula for $(\tilde{\circ}_{\har}^{\smallint})^{\vee}$ appears in a natural way in terms of products of coefficients of the $g^{(\xi)}$'s, $\xi \in \mu_{N}(K)$ (see Definition \ref{def Ad action 0,0} and Proposition \ref{prop dual composition}). In (ii) above, this is hidden behind the expressions of the $g^{(\xi^{i})}$'s in terms of $g$, and the shuffle equation for the unique $f$ such that $g = f^{-1}e_{1}f$. We recover it in (iii) above.

\subsubsection{Algebraic and topological properties}

\begin{Proposition} \label{harmonic Ihara properties}The $p$-adic pro-unipotent harmonic action of integrals is a well-defined group action of 
$(\Ad_{\tilde{\Pi}_{1,0}(K)_{o(1)}}(e_{1}),\circ^{\smallint_{1,0}}_{\Ad})$, continuous for the $\mathcal{N}_{D}$-topology on 
$\Ad_{\tilde{\Pi}_{1,0}(K)_{o(1)}}(e_{1})$, and the product topology on $\big(K \langle\langle e_{0 \cup \mu_{N}} \rangle\rangle_{\har}^{\smallint}\big)^{\mathbb{N}}$ of the $\mathcal{N}_{D}$-topologies on each factor $K \langle\langle e_{0 \cup \mu_{N}} \rangle\rangle_{\har}^{\smallint}$.
\end{Proposition}

\begin{proof} (a) Let $h \in K \langle \langle e_{0 \cup \mu_{N}} \rangle\rangle_{\har}^{\smallint}$ and $g \in \Ad_{\tilde{\Pi}_{1,0}(K)_{o(1)}}(e_{1})$ ; let a sequence of words $(w_{l})_{l\in \mathbb{N}}$ of the form $(e_{0}^{l}e_{\xi_{d+1}}e_{0}^{n_{d}-1}e_{\xi_{d}}\ldots e_{0}^{n_{1}-1}e_{\xi_{1}})_{l\in \mathbb{N}}$ ; we must show that $f(e_{0},(g^{(\xi)})_{\xi \in \mu_{N}(K)})[w_{l}]$ has a limit in $K$ when $l \rightarrow \infty$. Equation (\ref{eq:part c of proof}) gives a formula for $h(e_{0},(g^{(\xi)})_{\xi \in \mu_{N}(K)})[w_{l}]$ ; because of the assumption on $h$, that formula depends on $l$ in the following way, where the sum over $u$ is indexed by certain connected subsequences of $e_{\xi_{d+1}}e_{0}^{n_{d}-1}e_{\xi_{d}} \ldots e_{0}^{n_{1}-1}e_{\xi_{1}}$, thus does not depend on $l$,
\begin{equation} \label{eq:const lim} f(e_{0},(g^{(\xi)})_{\xi \in \mu_{N}(K)})[w_{l}] = \rho + \sum_{\xi \in \mu_{N}(K)} \sum_{u} \theta_{u,\xi} \sum_{b=1}^{l} g_{\xi}[e_{0}^{b-1}u],
\end{equation}
\noindent and where $\rho,\theta_{u,\xi} \in K$ do not depend on $l$. Because of the assumption that $g$ is summable, the right-hand side of equation (\ref{eq:const lim}) converges in $K$ when $l \rightarrow \infty$.
\newline\indent (b) By Proposition \ref{norm adjoint action}, $\circ_{\Ad}^{\smallint_{0,0}}$ and $\tau$ are continuous, and the map $\lim$ of Definition \ref{def limit map} is clearly continuous for restriction of the $\mathcal{N}_{D}$-topology to its source and target.
\newline\indent (c) Let $(h_{m})_{m \in\mathbb{N}} \in \big(K \langle \langle e_{0 \cup \mu_{N}} \rangle \rangle_{\har}^{\smallint}\big)^{\mathbb{N}}$, and let $g_{1},g_{2} \in \tilde{\Pi}_{1,0}(K)_{o(1)}$. By the associativity of the composition of formal power series in $K \langle \langle e_{0 \cup \mu_{N}} \rangle \rangle$, we have, for all $m \in \mathbb{N}$ :
\begin{equation} \label{eq:this equality} \tau(m)(g_{2}) \circ^{\smallint_{0,0}}_{\Ad} ( \tau(m)(g_{1}) \circ^{\smallint_{0,0}}_{\Ad} h_{m}) = ( \tau(m)(g_{2}) \circ^{\smallint_{1,0}}_{\Ad} \tau(m)(g_{1})) \circ^{\smallint_{0,0}}_{\Ad} h_{m}.
\end{equation}
\noindent By Definition \ref{dR-rt harmonic Ihara action}, by $\tau(m)(g_{2}) \circ^{\smallint_{1,0}}_{\Ad} \tau(m)(g_{1})= \tau(m)(g_{2} \circ^{\smallint_{1,0}}_{\Ad} g_{1})$ and by Proposition \ref{prop summable group ihara}, the right-hand side of (\ref{eq:this equality}) is in $K \langle \langle e_{0 \cup \mu_{N}} \rangle \rangle^{\lim}$ and its limit is the $m$-th term of the sequence $ (g_{2} \circ_{\Ad} g_{1}) \circ_{\har}^{\smallint} (h_{m})_{m \in \mathbb{N}}$. The Lemma \ref{following lemma} below shows that the expression $ g_{2} \circ_{\har}^{\smallint} \big( g_{1} \circ_{\har}^{\smallint} (h_{m})_{m \in \mathbb{N}} \big)$ is well-defined and equal to the sequence indexed by $m \in \mathbb{N}$ of limits of the left hand-side of (\ref{eq:this equality}).
\end{proof}

\begin{Lemma} \label{following lemma}Let $h' \in K \langle \langle e_{0 \cup \mu_{N}} \rangle \rangle_{\har}^{\smallint}$ and $g'_{1},g'_{2} \in \tilde{\Pi}_{1,0}(K)_{o(1)}$. Then, $g'_{2} \circ_{\Ad}^{\smallint_{0,0}} (g'_{1} \circ_{\Ad}^{\smallint_{0,0}} h')$ is in $K \langle \langle e_{0 \cup \mu_{N}} \rangle \rangle^{\lim}$ and we have
$$ \lim \big( g'_{2} \circ_{\Ad}^{\smallint_{0,0}} (g'_{1} \circ_{\Ad}^{\smallint_{0,0}} h') \big)
= \lim \big( g'_{2} \circ_{\Ad}^{\smallint_{0,0}} \lim (g'_{1} \circ_{\Ad}^{\smallint_{0,0}} h') \big). $$
\end{Lemma}

\begin{proof} The fact that $g'_{2} \circ_{\Ad}^{\smallint_{0,0}} (g'_{1} \circ_{\Ad}^{\smallint_{0,0}} h')$ is in $K \langle \langle e_{0 \cup \mu_{N}} \rangle \rangle^{\lim}$ follows from the previous proof. Let us prove the rest of the statement. Let $(w_{l})_{l \in \mathbb{N}}$ a sequence of the form 
$(e_{0}^{l}e_{\xi_{d+1}} e_{0}^{n_{d}-1}e_{\xi_{d}} \ldots e_{0}^{n_{1}-1}e_{\xi_{1}})_{l \in \mathbb{N}}$. Equation (\ref{eq:part c of proof}) applied two times gives a formula for $\big(g'_{2} \circ_{\Ad}^{\smallint_{0,0}} (g'_{1} \circ_{\Ad}^{\smallint_{0,0}} h')\big)[w_{l}]$. Since $h'$ is in $K \langle\langle e_{0 \cup \mu_{N}}\rangle\rangle_{\har}^{\smallint}$, we check that this depends on $l$ in the following way, where the sums over $u$ and $u',u''$ are over certain connected subsequences of $e_{\xi_{d+1}}e_{0}^{n_{d}-1}e_{\xi_{d}} \ldots e_{0}^{n_{1}-1}e_{\xi_{1}}$ :
\begin{multline} \label{eq:const lim lim} \big(g'_{2} \circ_{\Ad}^{\smallint_{0,0}} (g'_{1})^{(\xi)} \circ_{\Ad}^{\smallint_{0,0}} h')\big)[w_{l}] = \rho + 
\sum_{\xi \in \mu_{N}(K)} \sum_{u} \theta_{u,j} \sum_{b=1}^{l} {(g'_{1})}_{\xi}[e_{0}^{b-1}u] \\
+ \sum_{\xi',\xi''\in \mu_{N}(K)} \sum_{u',u''} \theta_{(u',u''),(\xi',\xi'')} \sum_{\substack{b',b''\geq 1\\ b'+b''=l}} {(g'_{1})}^{(\xi')}[e_{0}^{b'-1}u]{(g'_{2})}^{(\xi''}[e_{0}^{b''-1}u''],
\end{multline}
\noindent and where $\rho, \theta_{u,j},\theta_{(u',u''),(j',j'')} \in K$ do not depend on $l$. Equation (\ref{eq:const lim lim}) has a limit when $l \rightarrow \infty$ because $g'_{1}$ and $g'_{2}$ are summable. The limit when $l \rightarrow \infty$ of the third term of (\ref{eq:const lim lim}) is
\newline $\displaystyle \sum_{\xi' \in \mu_{N}(K)} \sum_{\xi''\in \mu_{N}(K)} \sum_{u',u''} \theta_{(u',u''),(\xi',\xi'')} \sum_{b'=1}^{\infty} {(g'_{1})}^{(\xi')}[e_{0}^{b'-1}u'] \sum_{b''=1}^{\infty} {(g'_{2})}_{(\xi'')}[e_{0}^{b''-1}u'']$. In particular, this formula separates $g'_{1}$ and $g'_{2}$ in the factors depending on $l$. This enables to check, first, that the limit when $l \rightarrow \infty$ of equation (\ref{eq:const lim lim}) is a function of $g'_{2}$ and $\lim \big( g'_{1} \circ_{\Ad}^{\smallint_{0,0}} h \big)$, and, then, that this function is exactly $\lim (g'_{2} \circ^{\smallint_{0,0}}_{\Ad} \lim (g'_{1} \circ^{\smallint_{0,0}}_{\Ad} h))$.
\end{proof}

\subsection{Application to the simplified equation of the Frobenius}

We combine \S2.1 and \S2.2 and we prove the integral part of the Theorem.

\subsubsection{Proof of equations (\ref{eq: property of ac 01}), (\ref{eq:formula for n=1}) and (\ref{eq:explicit inversion of series expansion N=1})}

We need some non-commutative generating series of weighted multiple harmonic sums ; in the next statement, we use the notation $w^{(p^{\alpha})}$ defined in Proposition \ref{previous formula}.

\begin{Definition} \label{def trois quatre un}(i) We define an element $\har_{m}$ of $K\langle\langle e_{0 \cup \mu_{N}} \rangle\rangle_{\har}^{\smallint}$ by, for all words : 
\newline  $\har_{m}[e_{0}^{l}e_{0}^{n_{d}-1}e_{\xi_{d}}\ldots e_{0}^{n_{1}-1}e_{\xi_{1}}] = \har_{m} \big((n_{i})_{d};(\xi_{i})_{d}\big)$.
\newline (ii) We define $\har_{m}^{(p^{\alpha})}\in K\langle\langle e_{0 \cup \mu_{N}} \rangle\rangle_{\har}^{\smallint}$ by $\har_{m}^{(p^{\alpha})}[w] = \har_{m}[w^{(p^{\alpha})}]$ for all words $w$.
\newline (iii) Let $\har_{p^{\alpha}\mathbb{N}},\har^{(p^{\alpha})}_{\mathbb{N}}\in (K\langle\langle e_{0 \cup \mu_{N}}\rangle\rangle_{\har}^{\smallint})^{\mathbb{N}}$ be respectively the sequences $(\har_{p^{\alpha}m})_{m \in \mathbb{N}}$, $(\har_{m}^{(p^{\alpha})})_{m\in\mathbb{N}}$.
\end{Definition}

We now prove equation (\ref{eq: property of ac 01}) which relates $p$-adic cyclotomic multiple zeta values and cyclotomic multiple harmonic sums by the pro-unipotent harmonic action of integrals. We note that the proof of the main theorem in \cite{I-1} also provides that $\Ad_{\Phi^{(\xi)}_{p,\alpha}}(e_{\xi}) \in K\langle \langle e_{0\cup \mu_{N}}\rangle\rangle_{o(1)}$ for any $\xi \in \mu_{N}(K)$ (\cite{I-1}, Corollary 4.3.2).

\begin{proof}
By Definition \ref{dR-rt harmonic Ihara action}, we have, with the notation $w_{l} = e_{0}^{l-1}e_{\xi_{d+1}} e_{0}^{n_{d}-1}e_{\xi_{d}} \ldots e_{0}^{n_{1}-1}e_{\xi_{1}}$,
\begin{multline} \label{eq:laleq}\big( \lim_{l \rightarrow \infty}  \tau(m) \Li_{p,X^{(p^{\alpha})}}^{\KZ}(z)\big(e_{0},
(\Ad_{\Phi^{(\xi)}_{p,\alpha}}(e_{\xi}))_{\xi \in \mu_{N}(K)}\big)[w_{l}][z^{m}]\big)_{m\in\mathbb{N}} = \big( \Phi_{p,\alpha} \circ_{\har}^{\smallint}
\har^{(p^{\alpha})}_{m} \big) \big((n_{i})_{d};(\xi_{i})_{d+1}\big).
\end{multline}
This combined to Proposition \ref{previous formula}, Corollary \ref{cor tend to 0} and Definition \ref{def trois quatre un} gives equation (\ref{eq: property of ac 01}), provided we can check that $\Ad_{\Phi_{p,\alpha}}(e_{1}) \in \Ad_{\tilde{\Pi}_{1,0}(K)_{o(1)}}(e_{1})$.
By \cite{I-1}, Corollary 4.3.2, we have $\Ad_{\Phi^{(\xi)}_{p,\alpha}}(e_{\xi}) \in K\langle \langle e_{0\cup \mu_{N}}\rangle\rangle_{o(1)} \cap \Ad_{\tilde{\Pi}_{1,0}(K)}(e_{1})$, and by \cite{J Assoc} we have $K\langle \langle e_{0\cup \mu_{N}}\rangle\rangle_{o(1)} \cap \Ad_{\tilde{\Pi}_{1,0}(K)}(e_{1})= \Ad_{\tilde{\Pi}_{1,0}(K)_{o(1)}}(e_{1})$.
\end{proof}

The simplest terms of equation (\ref{eq: property of ac 01}) are obtained in Example \ref{example circ har int}, in which we can take $g$, $(h_{m})_{m \in \mathbb{N}}$, $g \text{ }\circ_{\har}^{\smallint}\text{ } (h_{m})_{m \in \mathbb{N}}$ to be $\Phi_{p,\alpha}^{-1}e_{1}\Phi_{p,\alpha}$, $\har_{\mathbb{N}}^{(p^{\alpha})}$,  $\har_{p^{\alpha}\mathbb{N}}$ respectively.
\newline\indent We now prove the expansion of prime weighted multiple harmonic sums in terms of $p$-adic cyclotomic multiple zeta values mentioned in the theorem : equation (\ref{eq:formula for n=1}) and its $N=1$ case, equation (\ref{eq:explicit inversion of series expansion N=1}).

\begin{proof} In equation (\ref{eq:laleq}), the $m=1$ term is
$\displaystyle\sum_{b=0}^{\infty}\sum_{\xi \in \mu_{N}(K)} -\xi^{-p^{\alpha}}
\big(\Ad_{{\Phi^{(\xi)}_{p,\alpha}}}(e_{\xi})\big)[w_{b}]$.
\newline Indeed, $\har^{(p^{\alpha})}_{1}=\Li_{p,X^{(p^{\alpha})}}^{\KZ}(z^{p^{\alpha}})[z^{p^{\alpha}}]=\Li_{p,X^{(p^{\alpha})}}^{\KZ}(z)[z]$ is given in depth one by $\Li_{p,X^{(p^{\alpha})}}{\KZ}(z)[z][e_{0}^{l-b}e_{\xi}] = \xi^{-p^{\alpha}}$, and in zero in any depth $\geqslant 2$, because all the weighted multiple harmonic sums $\har_{1}$ have an empty domains of summation and are zero in all depths $\geq 2$. More details are in Lemma \ref{depth 0 harmonic Ihara} (iii).
\newline This and equation (\ref{eq: property of ac 01}) imply equation (\ref{eq:formula for n=1}), and, in the $N=1$ case, equation (\ref{eq:explicit inversion of series expansion N=1}).
\end{proof}

\begin{Remark} The remainder in the sum of series of equation (\ref{eq:formula for n=1}) has the following simple expression :
$(-1)^{d+1} \har_{p^{\alpha}}(\tilde{w}) - \sum\limits_{l'=0}^{l-1} \sum\limits_{\xi \in \mu_{N}(K)} -\xi^{-p^{\alpha}} \big( \Ad_{\Phi^{(\xi)}_{p,\alpha}}(e_{\xi})\big)[w_{l'}] =\Li_{p,\alpha}^{\dagger}[w_{l}][z^{p^{\alpha}}]$ and will find an interpretation of it in \cite{II-1}, \S4. Moreover, for all $r \in \{1,\ldots,p^{\alpha}-1\}$, we have $\Li_{p,\alpha}^{\dagger}[z^{r}][w_{l}^{(p^{\alpha})}] = p^{\weight(w)+l}r^{l}\har_{r}(w)$.
\end{Remark}

\subsubsection{Construction of a torsor containing $\har_{\mathbb{N}}^{(p^{\alpha})}$ for the pro-unipotent harmonic action of integrals}

\noindent We now prove that there exists a torsor for $\circ_{\har}^{\smallint}$ containing $\har_{\mathbb{N}}^{(p^{\alpha})}$ ; this will guarantee that equation (\ref{eq: property of ac 01}) characterizes $p$-adic cyclotomic multiple zeta values in terms of weighted multiple harmonic sums.
\newline\indent We need to prove that $\circ_{\har}^{\smallint}$ is compatible with the depth filtration, and to write explicitly the terms of extremal depth.

\begin{Definition} We denote $e_{\xi_{d+1}}e_{0}^{n_{d}-1}e_{\xi_{d}} \ldots e_{0}^{n_{1}-1}e_{\xi_{1}} \in \mathcal{O}_{\har}^{\sh,e_{0 \cup \mu_{N}}}$ by $((n_{i});(\xi_{i}))_{d,d+1}$, with $(n_{i})=\emptyset$ if $d=0$, and we say that such words have depth $d$ and weight $\sum\limits_{i=1}^{d}n_{i}$ ; we denote by ${\mathcal{O}_{\har,\ast,d}^{\sh,e_{0 \cup \mu_{N}}}}$, the subspace generated by such words of depth $d$.
\end{Definition}

\begin{Lemma} \label{depth 0 harmonic Ihara}(i) For any $d \in \mathbb{N}$, the map $(\tilde{\circ}_{\har}^{\smallint})^{\vee,T}$ sends 
$$ \mathcal{O}^{\sh,e_{0 \cup \mu_{N}}}_{\har,\ast,d}
\rightarrow \bigoplus_{d'=0}^{d} {\mathcal{O}_{\har,\ast,d'}^{\sh,e_{0 \cup \mu_{N}}}} \otimes T \big( \tau(m_{f})\widehat{\mathcal{O}^{\sh,e_{0 \cup \mu_{N}}}_{\ast,d-d'}} \otimes ( \underset{\xi \in \mu_{N}(K)}{\oplus} \mathbb{Q}\xi ) \big). $$
\noindent (ii) In depth $0$, the map $(\circ_{\har}^{\smallint})^{\vee}$ is $\big( \emptyset;\xi \big) \mapsto \big(\emptyset;\xi\big) \otimes 1$ for all $\xi \in \mu_{N}(K)$.
\newline (iii) The term in ${\mathcal{O}_{\har}^{\sh,e_{0 \cup \mu_{N}}}}_{\ast,0} \otimes T \big(  \tau(m)\widehat{\mathcal{O}^{\sh,e_{0 \cup \mu_{N}}}_{\ast,d}} \otimes ( \oplus_{\xi \in \mu_{N}(K)} \mathbb{Q}\xi) \big)$ of any $\displaystyle(\circ_{\har}^{\smallint})^{\vee,T}\big((n_{i})_{d};(\xi_{i}))_{d+1}\big)$, in the sense of (i), is $\displaystyle \sum_{\xi \in \mu_{N}(K)} \big(\emptyset;\xi\big) \otimes \sum_{b=0}^{+\infty} m_{f}^{b+n_{d}+\ldots+n_{1}} e_{0}^{b}e_{\xi_{d+1}}\ldots e_{0}^{n_{1}-1}e_{\xi_{1}} \otimes \xi$.
\end{Lemma}

\begin{proof} This follows from Proposition \ref{prop dual composition}, equation (\ref{eq:part c of proof}) and the definition of $\circ_{\har}^{\smallint}$ (Definition \ref{dR-rt harmonic Ihara action}).
\end{proof}

We now prove the torsor structure mentioned in (i) of the theorem, i.e. that the orbit of $\har_{\mathbb{N}}^{(p^{\alpha})}$ is a torsor containing $\har_{\mathbb{N}}^{(p^{\alpha})}$ for the $p$-adic pro-unipotent harmonic action of integrals.

\begin{proof} Let $\mathcal{E}^{\smallint}_{\har} \subset \big( K \langle\langle e_{0 \cup \mu_{N}}\rangle\rangle_{\har}^{\smallint}\big)^{\mathbb{N}}$ be the subset of elements $h=(h_{m})_{m\in\mathbb{N}}$ such that the maps $m\in \mathbb{N}^{\ast} \mapsto h_{m}\big(\emptyset;\xi\big)$, $\xi \in \mu_{N}(K)$, are linearly independent over the ring $A(\mathbb{Z}_{p})$ of rigid analytic functions of $m\in \mathbb{Z}_{p}$. Then 
\newline (a) $\har_{\mathbb{N}}^{(p^{\alpha})} \in \mathcal{E}^{\smallint}_{\har}$. This is because for all $\xi \in \mu_{N}(K)$, we have $\har_{m}(\emptyset ;\xi) = \xi^{-m}$ and the result follows from the invertibility of a Vandermonde matrix. 
\newline (b) $\mathcal{E}^{\smallint}_{\har}$ is stable by $\circ_{\har}^{\smallint}$. This follows from part (ii) of Lemma \ref{depth 0 harmonic Ihara}.
\newline (c) $\circ_{\har}^{\smallint}$ restricted to $\mathcal{E}^{\smallint}_{\har}$ is free. One proves by induction on $d$ that $\circ_{\har}^{\smallint}$ truncated to depths at most $d$ is free, by (iii) of Lemma \ref{depth 0 harmonic Ihara}.
\newline This implies that the orbit of $\har_{\mathbb{N}}^{(p^{\alpha})}$ is included in  $\mathcal{E}^{\smallint}_{\har}$ and is a torsor.
\end{proof}

\subsubsection{The harmonic Frobenius of integrals}

\begin{Definition} \label{def har Frob}Let the harmonic Frobenius of integrals, iterated $\alpha$ times, be the map
\begin{center}
$(\tau(p^{\alpha})\phi^{\alpha})^{\smallint}_{\har} :
\begin{array}{c} \big( K\langle\langle e_{0 \cup \mu_{N}} \rangle\rangle_
{\har}^{\smallint} \big)^{\mathbb{N}}\rightarrow \big( K\langle\langle e_{0 \cup \mu_{N}} \rangle\rangle_{\har}^{\smallint}\big)^{\mathbb{N}}
\\ f \mapsto \Phi_{p,\alpha}^{-1}e_{1}\Phi_{p,\alpha} \circ_{\har}^{\smallint} \sigma^{\alpha}(f) \end{array}.$
\end{center}
\end{Definition}

Indeed, the passage from the Frobenius to the harmonic Frobenius commutes with the iteration ; see \cite{I-3}.

\begin{Proposition} The harmonic Frobenius of integrals is continuous for the product indexed by $\mathbb{N}$ of the $\mathcal{N}_{D}$-topology on $K\langle\langle e_{0 \cup \mu_{N}} \rangle\rangle_{\har}^{\smallint}$.
\end{Proposition}

\begin{proof} Follows from the continuity of $\circ_{\har}^{\smallint}$ (Proposition \ref{harmonic Ihara properties}).
\end{proof}

With Definition \ref{def har Frob}, equation (\ref{eq: property of ac 01}) is restated as

\begin{equation} (\phi^{\alpha})^{\smallint}_{\har} (\har_{\mathbb{N}}) = \har_{p^{\alpha}\mathbb{N}}.
\end{equation}

In Definition \ref{dR-rt harmonic Ihara action} and Definition \ref{def har Frob}, the adjective ``harmonic'' means ``adapted to weighted multiple harmonic sums'' : we will check in the next sections that these objects are indeed natural as operations on weighted multiple harmonic sums.

\section{Setting for the pro-unipotent harmonic action of series}

We define (\S3.1) and study (\S3.2,\S3.3) a generalization of cyclotomic multiple harmonic sums which we call localized cyclotomic multiple harmonic sums. The term localized refers to the inversion of a differential operator which is implicit behind the definition. This is a preliminary to \S4.

\subsection{Localized cyclotomic multiple harmonic sums}

\subsubsection{Cyclotomic multiple harmonic sums}

The cyclotomic multiple harmonic sums are the following numbers, with the notations of equation (\ref{eq:multiple harmonic sum}) and $m_{0} \in \mathbb{N}$,
$$ \frak{h}_{m_{0},m} \big( (n_{i})_{d};(\xi_{i})_{d+1}\big) = \sum_{(m_{1},\ldots,m_{d}) \in \Delta_{m_{0},m}^{\mathbb{N}^{d}}} \frac{
\big( \xi_{1}\big)^{m_{0}}\big( \frac{\xi_{2}}{\xi_{1}}\big)^{m_{1}} \ldots \big( \frac{\xi_{d+1}}{\xi_{d}}\big)^{m_{d}}
\big(\frac{1}{\xi_{d+1}}\big)^{m}}{m_{1}^{n_{1}} \ldots m_{d}^{n_{d}}} $$ 
where, for $d \in \mathbb{N}^{\ast}$, $m_{0},m \in \mathbb{N}^{\ast}$, 
$$\Delta^{\mathbb{N}^{d}}_{m_{0},m}=\{(m_{1},\ldots,m_{d}) \in \mathbb{N}^{d} \text{ }|\text{ }m_{0}<m_{1}<\ldots<m_{d}<m\}$$ 
and $\Delta^{\mathbb{N}^{d}}_{m}=\Delta^{\mathbb{N}^{d}}_{0,m}$ ; the weighted cyclotomic multiple harmonic sums are the numbers
$$ \har_{m_{0},m}\big( (n_{i})_{d};(\xi_{i})_{d+1}\big) = (m-m_{0})^{n_{1}+\ldots+n_{d}} \frak{h}_{m_{0},m}\big( (n_{i});(\xi_{i})\big)_{d}. $$
The prime weighted multiple harmonic sums are the numbers $\har_{p^{\alpha}}\big( (n_{i})_{d};(\xi_{i})_{d+1}\big)$ (\cite{I-1}, Definition B.0.1).
\newline\indent We call harmonic word a sequence  $((n_{i})_{d};(\xi_{i})_{d+1})$ where $d$ and the $n_{i}$'s (for $1\leqslant i \leqslant d$) are positive integers and the $\xi_{i}$'s are $N$-th roots of unity (for $1\leqslant i \leqslant d+1$). Let $\Wd_{\har}(e_{0 \cup \mu_{N}})$ be the set of harmonic words. We define a natural series counterpart of $K \langle \langle e_{0 \cup \mu_{N}}\rangle\rangle_{\har}^{\smallint}$ from Definition \ref{harmonic generating series set}, which is isomorphic as a $K$-vector space.

\begin{Definition} Let $K \langle \langle e_{0 \cup \mu_{N}}\rangle\rangle_{\har}^{\Sigma} = \bigg\{ \sum\limits_{w \in \Wd_{\har}(e_{0 \cup \mu_{N}})} \lambda_{w} w \text{ }\bigg|\text{ }\forall w, \lambda_{w} \in K \bigg\}$.
\end{Definition}

We will view $\har_{m}$, $\har_{m}^{(p^{\alpha})}$ from Definition \ref{def trois quatre un} as elements of $K \langle \langle e_{0 \cup \mu_{N}}\rangle\rangle_{\har}^{\Sigma}$, and  $\har_{\mathbb{N}}$, $\har_{\mathbb{N}}^{(p^{\alpha})}$ from Definition \ref{def trois quatre un} as elements of $\big( K \langle\langle e_{0 \cup \mu_{N}}\rangle\rangle_{\har}^{\Sigma}\big)^{\mathbb{N}}$, via this isomorphism. Let, more generally, with the notation $w^{(p^{\alpha})}$ introduced in Proposition \ref{previous formula}.

\begin{Definition} \label{def non comm gen seri wmhs} For $m_{0},m \in \mathbb{N}^{\ast}$, let :
\newline (i) $\har_{m_{0},m}= \sum\limits_{w \in \Wd_{\har}(e_{0 \cup \mu_{N}})} \har_{m_{0},m}(w)w \in K \langle\langle e_{0 \cup \mu_{N}}\rangle\rangle_{\har}^{\Sigma}$.
\newline (ii) For $I,J\subset \mathbb{N}$ such that $I\times J\simeq \mathbb{N}$, let $\har_{I,J} =( \har_{m_{0},m})_{(m_{0},m) \in I \times J} \in \big( K \langle\langle e_{0 \cup \mu_{N}}\rangle\rangle_{\har}^{\Sigma} \big)^{\mathbb{N}}$.
\newline (iii) We define similarly $\har_{m_{0},m}^{(p^{\alpha})}$ and  $\har_{I,J}^{(p^{\alpha})}$, by replacing $w$ by $w^{(p^{\alpha})}$ in (i).
\end{Definition}

\subsubsection{Localized cyclotomic multiple harmonic sums}

This is the central object of this \S3.

\begin{Definition} \label{quatre quatre un}
(i) A localized harmonic word is a sequence $((n_{i})_{d};(\xi_{i}))_{d+1}$ as above except that we allow the $n_{i}$'s to be any elements of $\mathbb{Z}$. Let $\Wd_{\har}(e_{0\cup\mu_{N}})_{\loc}$ be the set of localized harmonic words.
\newline (ii) Let $K \langle\langle e_{0 \cup \mu_{N}}\rangle\rangle_{\har,\loc}^{\Sigma} = \bigg\{ \sum\limits_{w \in \Wd_{\har}(e_{0 \cup \mu_{N}})_{\loc}} \lambda_{w} w \text{ }|\text{ }\forall w, \lambda_{w} \in K \bigg\}$.
\end{Definition}

\begin{Definition} \label{def localized mult har sums} Let $m_{0},m \in \mathbb{N}^{\ast}$, $((n_{i})_{d};(\xi_{i}))_{d+1}$ a harmonic word. Let $i_{1},\ldots,i_{r} \in \{1,\ldots,d\}$ be the elements such that $n_{i}>0$ and $n_{i-1}<0$, or $n_{i}>0$ and $i=1$. We call localized cyclotomic harmonic sums the numbers
\begin{equation}
\label{eq:localized cyclotomic multiple harmonic sums}
\frak{h}_{m_{0},m}\big( (n_{i});(\xi_{i})\big)_{d}= \sum_{(m_{1},\ldots,m_{d}) \in \Delta_{(i_{1},\ldots,i_{r}),m_{0},m}} \frac{\big( \xi_{1}\big)^{m_{0}}\big( \frac{\xi_{2}}{\xi_{1}}\big)^{m_{1}} \ldots \big( \frac{\xi_{d+1}}{\xi_{d}}\big)^{m_{d}}\big(\frac{1}{\xi_{d+1}}\big)^{m}}{m_{1}^{n_{1}} \ldots m_{d}^{n_{d}}},
\end{equation}
where 
$$ \Delta^{\mathbb{N}^{d}}_{(i_{1},\ldots,i_{r}),m_{0},m} = \{ (m_{1},\ldots,m_{d}) \in \mathbb{N}^{d}\text{ }|\text{ }m_{0}< \ldots < m_{i_{1}-1} \leq m_{i_{1}} < \ldots < m_{i_{r}-1} \leq m_{i_{r}} < \ldots < m \}.$$
We call weighted localized multiple harmonic sums the numbers :
$$\har_{m_{0},m}(w) = (m-m_{0})^{n_{1}+\ldots+n_{d}} \frak{h}_{m_{0},m}(w)$$
For all $m \in \mathbb{N}^{\ast}$, we denote by $\frak{h}_{m}=\frak{h}_{0,m}$, $\har_{m}=\har_{0,m}$.
\end{Definition}

\begin{Definition} \label{gen series harm}For $m \in \mathbb{N}^{\ast}$, let  :
\newline (i) $\har_{m_{0},m,\loc}= \sum\limits_{w \in \Wd_{\har}(e_{0 \cup \mu_{N}})_{\loc}} \har_{m_{0},m}(w)w \in  K \langle\langle e_{0 \cup \mu_{N}}\rangle\rangle_{\har,\loc}^{\Sigma}$, and $\har_{m,\loc}=\har_{0,m,\loc}$.
\newline (ii) For $I,J\subset \mathbb{N}$ such that $I\times J\simeq \mathbb{N}$, let $\har_{I,J,\loc} =( \har_{m_{0},m,\loc})_{(m_{0},m) \in I \times J} \in \big( K \langle\langle e_{0 \cup \mu_{N}}\rangle\rangle_{\har,\loc}^{\Sigma} \big)^{\mathbb{N}}$, and $\har_{J,\loc}=\har_{\{0\},J,\loc}$
\newline (iii) We define similarly $\har_{m_{0},m,\loc}^{(p^{\alpha})}$, $\har^{(p^{\alpha})}_{m,\loc}$,
$\har_{I,J,\loc}^{(p^{\alpha})}$, $\har^{(p^\alpha)}_{J,\loc}$ by replacing $w$ by $w^{(p^{\alpha})}$ in (i).
\end{Definition}

\begin{Remark} Most of the computations in the rest of this paper can be immediately extended to the generalizations of localized multiple harmonic sums obtained as follows : replacing the factors $m_{i}^{n_{i}}$ ($i=1,\ldots,d$) in equation (\ref{eq:localized cyclotomic multiple harmonic sums}), by $\chi_{i}(m_{i})$, with $\chi_{i}$ group morphisms $(K^{\ast},\times) \rightarrow (K^{\ast},\times)$ which are analytic on $\{z \in K^{\ast}\text{ }|\text{ }|z-1|_{p}\leq \frac{1}{p^{\alpha}}\}$, and are thus locally analytic on $K^{\ast}$ ; replacing the weight $\sum\limits_{i=1}^{d} n_{i}$ of a sequence $(n_{1},\ldots,n_{d})$ by $\sum\limits_{i=1}^{d}-\frac{\log_{p}(\chi_{i}(p))}{\log_{p}(p)}$ ; replacing the factor $\frac{\xi_{1}^{m_{0}}}{\xi_{d+1}^{m}}$in Definition \ref{def localized mult har sums} by any $\frac{\tilde{\xi}_{1}^{m_{0}}}{\tilde{\xi}_{d+1}^{m}}$, with $\tilde{\xi}_{1},\tilde{\xi}_{d+1} \in \mu_{N}(K)$, such that $\tilde{\xi}_{1}\tilde{\xi}_{d+1}^{-1} = \xi_{1}\xi^{-1}_{d+1}$ ; replacing 
$m_{0},m \in \mathbb{N}^{\ast}$ by elements of $\mathbb{Z}$ such that the right-hand side of equation (\ref{eq:localized cyclotomic multiple harmonic sums}) is well-defined.
\end{Remark}

\subsubsection{Operations on the indices of localized cyclotomic multiple harmonic sums}

The next definition is an analogue of the notion of subword from Definition \ref{def subword} (i).

\begin{Definition} Let $w=(n_{i};\xi_{i})_{d} \in  \Wd_{\har}^{\loc}(e_{0 \cup \mu_{N}})$ be a localized harmonic word and $S =[a,b] \subset \{1,\ldots,d\}$. We denote 
$(n_{a},\ldots,n_{b};\xi_{a},\ldots,\xi_{b+1}) \in \Wd_{\har}^{\loc}(e_{0 \cup \mu_{N}})$ by $w|_{S}$.
\end{Definition}

The next definition is an analogue of the notion of quotient word from Definition \ref{def subword} (v), which will appear implicitly afterwards.

\begin{Definition} \label{definition of connected partition}Let $S$ be a subset of $\mathbb{N}$.
\noindent\newline (i) A connected partition of $S$ is a partition of $S$ into segments.
\newline (ii) An increasing connected partition of $S$ is a connected partition of $S$ with an order on the corresponding set of parts of $S$, such that if a part $C$ is inferior to a part $C'$ for this order, we have $j<j'$ in $\mathbb{N}$ for all $j \in C$ and $j' \in C'$.
\newline (iii) The canonical increasing connected partition of $S$ is the increasing connected partition of $S$ defined by the segments included in $S$ and maximal for the inclusion, which we call the connected components of $S$.
\newline (iv) Let $S$ be a subset of $\mathbb{N}$. We call the boundary of $S$ and denote by $\partial S$ the subset of $S$ made of the elements $x$ such that $x-1 \not\in S$ or $x+1 \not\in S$.
\end{Definition}

\subsection{Computation of totally negative cyclotomic multiple harmonic sums}

\begin{Definition} We say that a localized harmonic word $w = \big( (n_{i});(\xi_{i}) \big)_{d}$ is totally negative if, for all $i$, $n_{i}<0$ ; in that case we also say, for all $m_{0},m$, that $\frak{h}_{m_{0},m}(w)$ is totally negative. Let $\Wd_{\har}^{-}(e_{0 \cup \mu_{N}})$ be the set of totally negative harmonic words.
\end{Definition}

\begin{Proposition-Definition} \label{numbers mathcal B}For any $w=\big((n_{i});(\xi_{i}) \big)_{d}\in\Wd_{\har}^{-}(e_{0 \cup \mu_{N}})$, there exists a unique sequence $(\mathcal{B}_{\delta_{0},\delta,\xi_{0},\xi}^{w})_{\substack{\delta_{0},\delta\in \{0,\ldots,\sum_{i=1}^{d}|n_{i}|+d+1\} \\ \xi_{0},\xi \in \mu_{N}(K)}}$ of elements of the $N$-th cyclotomic field, such that, for all $m_{0},m$ we have
\begin{equation} \label{eq:recurrence} \har_{m_{0},m}(w) =
\sum_{\xi_{0},\xi \in \mu_{N}(K)} \sum_{\delta_{0},\delta=0}^{l_{1}+\ldots+l_{d}+d+1}\mathcal{B}_{\delta_{0},\delta,\xi_{0},\xi}^{w} m_{0}^{\delta_{0}}m^{\delta} \xi_{0}^{m_{0}}\xi^{m}.
\end{equation}
Moreover, for all $\delta_{0},\delta,\xi_{0},\xi$, we have $\displaystyle v_{p}(\mathcal{B}_{\delta_{0},\delta,\xi_{0},\xi}^{w})\geq - d - \frac{\log(|n_{1}|+\ldots+|n_{d}|+d+1)}{\log(p)}$.
\end{Proposition-Definition}

\begin{proof} The existence of these numbers is proved by induction on $d$, using that, for $m \in \mathbb{N}^{\ast}$, $l \in \mathbb{N}^{\ast}$ we have : $\sum\limits_{m_{1}=0}^{m-1} m_{1}^{l} =\sum\limits_{\delta=0}^{l+1} \frac{1}{l+1}{l+1 \choose \delta} B_{l+1-\delta}T^{\delta}$ and $\sum\limits_{m_{1}=0}^{m-1} m_{1}^{l}T^{m_{1}}= (T\frac{d}{dT})^{l}(\sum\limits_{m_{1}=0}^{m} T^{m_{1}}) = (T\frac{d}{dT})^{l}\big( \frac{T^{m}-1}{T-1}\big)$, where $T$ is a formal variable, to which we can substitute an element of $\mu_{N}(K) \setminus \{1\}$.
\newline\indent The uniqueness follows from the uniqueness of the coefficients of a polynomial and the invertibility of a Vandermonde matrix.
\newline\indent The bound of valuations is proved by induction on $d$ by Von-Staudt Clausen's theorem, as well as $v_{p}(\frac{1}{l})\geq - \frac{\log(l)}{\log(p)}$ for all $l \in \mathbb{N}^{\ast}$, and $|\xi-1|_{p}=1$ for all $\xi \in \mu_{N}(K) \setminus \{1\}$.
\end{proof}

\begin{Notation} (i) For all $\delta_{0},\delta \geqslant l_{1}+\ldots+l_{d}+d$, $\xi_{0},\xi \in \mu_{N}(K)$, $w \in \Wd_{\har}^{-}(e_{0 \cup \mu_{N}})$, let $\mathcal{B}_{\delta_{0},\delta,\xi_{0},\xi}^{w}=0$ 
\newline (ii) For all $\delta \in \mathbb{N}$, $w \in \Wd_{\har}^{-}(e_{0 \cup \mu_{N}})$, $\xi_{0},\xi \in \mu_{N}(K)$, let  $\mathcal{B}_{\delta,\xi_{0},\xi}^{w}=\mathcal{B}_{0,\delta,\xi_{0},\xi}^{w}$
\newline (iii) For all $\delta \in \mathbb{N}$, $l,l_{1},l_{2} \in \mathbb{N}^{\ast}$, let $\mathcal{B}_{\delta}^{l}=\mathcal{B}_{0,\delta,1,1}^{(l;1,1)}$, $\mathcal{B}_{\delta}^{l_{1},l_{2}}=\mathcal{B}_{0,\delta,1,1}^{(l_{1},l_{2};1,N,1)}$
\newline (iv) For $l\in \mathbb{N}^{\ast}$, $\delta \in \mathbb{N}$, $\tilde{\xi} \mu_{N}(K)$, let $\mathcal{B}_{\delta,\xi}^{l}(\tilde{\xi})=\mathcal{B}_{0,\delta,N,\xi}^{(l;\tilde{\xi},1)}$.
\newline (v) We omit $\xi_{0},\xi$ in all the notations if $N=1$.
\end{Notation}

We note that if $N \not= 1$, the coefficients $\mathcal{B}$ depend on the $N$-th roots of unity via rational functions in
$\mathbb{Z}[T_{1},\ldots,T_{N-1},\frac{1}{T_{1}},\ldots,\frac{1}{T_{N-1}},\frac{1}{T_{1}-1},\ldots,\frac{1}{T_{N-1}-1}]$. This type of expression already appeared in \cite{I-1}, \S3.

\subsection{Formulas on adding and multiplying upper bounds of the domain of summations\label{multiplication rational}}

We write some analogues for cyclotomic multiple harmonic sums of some basic rules of computation on iterated integrals (more details on this analogy will appear in \cite{II-3}). 

\begin{Notation} In the next statements, the abbreviation i.c.p. stands for increasing connected partition, in the sense of Definition \ref{definition of connected partition}.
\end{Notation}

\subsubsection{Addition of upper bounds of domains of summation}

We want to relate $\har_{m+m'}$ to $\har_{m}$ and $\har_{m'}$, for any $m,m' \in \mathbb{N}^{\ast}$. If we stay in the $N$-th cyclotomic field, what we obtain is a formula for the "splitting" of the domain of summation of localized cyclotomic multiple harmonic sums.

\begin{Proposition}\label{addition}
Let $m,m_{0} \in \mathbb{N}^{\ast}$, such that $m_{0} < m$. Let 
$ \tilde{m}_{1}<\ldots<\tilde{m}_{r}\in \{m_{0},\ldots,m-1\}$. We also denote by $\tilde{m}_{0} = m_{0}$ and $\tilde{m}_{r+1} = m$. Then we have, for all harmonic words $w=\big((n_{i})_{d};(\xi_{i})_{d+1}\big)$ :
\begin{equation} \label{eq:formula splitting} \frak{h}_{m_{0},m}(w) = \sum_{\substack{0 \leq \tilde{r} \leq r \\ 1\leq i_{1}<\ldots<i_{\tilde{r}} \leq d
\\ 1 \leq \tilde{i}_{1}<\ldots<\tilde{i}_{\tilde{r}} \leq r
\\ \{1,\ldots,d\} - \{i_{1},\ldots,i_{\tilde{r}}\} = S_{0} \amalg \ldots \amalg S_{r} \text{ i.c.p.}}}
\prod_{\tilde{i}=1}^{\tilde{r}} \frac{1}{\tilde{m}_{\tilde{i}}^{n_{i_{\tilde{i}}}}} \prod_{a=0}^{m-1} \frak{h}_{\tilde{m}_{a},\tilde{m}_{a+1}}(w|_{S_{a}}).
\end{equation}
\end{Proposition}

\begin{proof} For each $(m_{1},\ldots,m_{d})$ in the domain of summation $\Delta^{\mathbb{N}^{d}}_{m_{0},m}$ of $\frak{h}_{m_{0},m}$, we let $\{i_{1},\ldots,i_{\tilde{r}}\}=\{i \in \{1,\ldots,d\} \text{ }|\text{ } m_{i} \in \{\tilde{m}_{1},\ldots,\tilde{m}_{d}\}\}$, and $\{\tilde{i}_{1},\ldots,\tilde{i}_{\tilde{r}}\}=\{\tilde{i} \in \{1,\ldots,r\} \text{ }|\text{ } \tilde{m}_{\tilde{i}} \in \{m_{1},\ldots,m_{d}\}\}$, with $i_{1}<\ldots<i_{r}$ and $\tilde{i}_{1}<\ldots< \tilde{i}_{\tilde{r}}$. In particular,  $m_{i_{\tilde{i}}}=\tilde{m}_{\tilde{i}}$.
\end{proof}

\begin{Example} Equation (\ref{eq:formula splitting}) in the case $r=1$ is
\begin{multline} \frak{h}_{m_{0},m} \big((n_{i});(\xi_{i})\big)_{d} = \sum_{i_{1}=1}^{d}  \frak{h}_{m_{0},\tilde{m}_{1}} \big((n_{i});(\xi_{i}) \big)_{i_{1}}\frak{h}_{\tilde{m}_{1},m}
\big((n_{i+i_{1}});(\xi^{j_{i+i_{1}}}) \big)_{d-i_{1}} 
\\ +
\sum_{i_{1}=1}^{d} \frac{1}{\tilde{m}_{1}^{n_{i_{1}}}} \frak{h}_{m_{0},\tilde{m}_{1}} \big((n_{i});(\xi_{i})\big)_{i_{1}-1} \frak{h}_{\tilde{m}_{1},m}
\big((n_{i+i_{1}});(\xi^{j_{i+i_{1}}}) \big)_{d-i_{1}}.
\end{multline}
\end{Example}

\subsubsection{Multiplication of upper bounds of domains of summation}

We now want to relate $\har_{mm'}$ to $\har_{m}$ and $\har_{m'}$, for any $m,m' \in \mathbb{N}^{\ast}$. If we stay in the $N$-th cyclotomic field, what we obtain is the following formula, which express the Euclidean division by $m$ of the coordinates of elements of the domain of summation of $\har_{mm'}$. In the next statement, we use the convention that $\har_{m_{0},m}(\emptyset)=1$.

\begin{Proposition} \label{multiplication} For all harmonic words $w=\big((n_{i})_{d};(\xi_{i})_{d+1} \big)$, we have :
\begin{equation} \label{eq:equation multiplication}
\har_{\mu m_{0},\mu m}(w)
= \sum_{\substack{0 \leq \tilde{r} \leq m-1 \\ 1\leq i_{1}<\ldots<i_{\tilde{r}} \leq d
		\\ 1 \leq \tilde{i}_{1}<\ldots<\tilde{i}_{\tilde{r}} \leq m-1
		\\ \{1,\ldots,d\} - \{i_{1},\ldots,i_{\tilde{r}}\} = S_{0} \amalg \ldots \amalg S_{m-1} \text{ i.c.p.}}}
\prod_{\tilde{i}=1}^{\tilde{r}} \frac{1}{\tilde{i}^{n_{i_{\tilde{i}}}}} \prod_{a=0}^{m-1} \har_{\mu a,\mu(a+1)}(w|_{S_{a}}).
\end{equation}
\end{Proposition}

\begin{proof} By applying Proposition \ref{addition} to $\frak{h}_{\mu m_{0},\mu m}$ and $\{\tilde{m}_{1},\ldots,\tilde{m}_{r}\} = \{\mu,2\mu,\ldots,(m-1)\mu\}$.
\end{proof}

\begin{Example} (i) If $w$ has depth one, the right-hand side of (\ref{eq:equation multiplication}) has two terms : this corresponds the partition  $\Delta^{\mathbb{N}^{1}}_{\mu m_{0},\mu m} = \{\mu|m_{1}\} \amalg \{\mu \nmid m_{1}\}$.
\newline (ii) If $w$ has depth two, the right-hand side of (\ref{eq:equation multiplication}) has five terms ; this corresponds to the partition $ \Delta^{\mathbb{N}^{2}}_{\mu m_{0},\mu m}=\big\{\mu|m_{1}, \mu|m_{2} \big\} 
\text{ }\amalg\text{ } 
\big\{\mu\nmid m_{1},\mu\nmid m_{2}, \big[\frac{m_{1}}{\mu}\big] = \big[ \frac{m_{2}}{\mu}\big] \big\}
\text{ }\amalg\text{ } 
\big\{ \mu\nmid m_{1},\mu\nmid m_{2}, \big[\frac{m_{1}}{\mu}\big] < \big[\frac{m_{2}}{\mu}\big] \big\}
\amalg\text{ }  \big\{\mu|m_{1},\mu\nmid m_{2}\big\} \text{ }\amalg\text{ } 	
\big\{\mu\nmid m_{1},\mu|m_{2} \big\}$.
\end{Example}

\section{The pro-unipotent harmonic action of series}

We construct a ``localized pro-unipotent harmonic action of series'' $\circ_{\har,\loc}^{\Sigma}$ (Proposition-Definition \ref{loc series harmonic action}) and a ``map of delocalization'' $\deloc$ (Proposition-Definition \ref{loc for har n}). We give explicit formulas for these two maps ; the most significative one combinatorially is the formula for $\deloc$ (Proposition \ref{formula for loc}). Composing these two maps gives the pro-unipotent harmonic action of series (Proposition-Definition \ref{series harmonic action}) and proves the ``series'' part of the theorem.

\subsection{The localized pro-unipotent harmonic action of series}

We need first a $p$-adic formula for shifting the bounds of the domain of summation of a cyclotomic multiple harmonic sums.

\begin{Definition} For $w=((n_{i});(\xi_{i}))_{d+1}$ a harmonic word and $l_{1},\ldots,l_{d} \in \mathbb{N}$, let $\shft_{l_{1},\ldots,l_{d}}(w)=((n_{i}+l_{i})_{d},(\xi_{i})_{d+1})$.
\end{Definition}

\begin{Lemma} \label{lemma p-adic shifting} Let $m_{0},m,\delta\in \mathbb{N}$. Assume that $|\delta|_{p}<|m'-\delta|_{p}$ for all $m' \in \{m_{0},\ldots,m\}$ ; then we have, for all harmonic words $w$,
$$ \har_{m_{0}+\delta,m+\delta}(w)
= \sum_{l_{1},\ldots,l_{d} \geqslant 0} \bigg( \prod_{i=1}^{d} \delta^{l_{i}} {-n_{i} \choose l_{i}} \bigg) \text{ } \har_{m_{0},m}( \shft_{l_{1},\ldots,l_{d}}(w)). $$
\end{Lemma}

\begin{proof} We make the change of variable $(m_{1},\ldots,m_{d})=(m'_{1}+\delta,\ldots,m'_{d}+\delta)$ in the domain of summation $\Delta_{m_{0},m}^{\mathbb{N}^{d}}$,
and we write the power series expansion 
$(m'_{i}+\delta)^{-n_{i}} = {m'_{i}}^{-n_{i}} \sum_{l_{i}\geq 0} {-n_{i} \choose l_{i}} \big( \frac{\delta}{m'_{i}}\big)^{l_{i}}$ for all $i \in \{1,\ldots,d\}$.
\end{proof}

The next proposition continues in $K$ the computation of Proposition \ref{multiplication}, assuming $\mu=p^{\alpha}$.

\begin{Definition} Let $K\langle \langle e_{0 \cup \mu_{N}} \rangle \rangle_{\har,o(1)}^{\Sigma} \subset K\langle \langle e_{0 \cup \mu_{N}} \rangle \rangle_{\har}^{\Sigma}$ be the subset of elements $f$ such that for all sequence $(w_{l})_{l\in\mathbb{N}}$ of words of bounded depth and such that $\displaystyle \weight(w_{l}) \underset{l\rightarrow \infty}{\rightarrow} \infty$, we have $\sum\limits_{l \geqslant 0} |f[w_{l}]|_{p}<\infty$.
\end{Definition}

Below we use the notation $g \circ_{\har}^{\Sigma} f = \circ_{\har}^{\Sigma} (g,f)$.

\begin{Proposition-Definition} \label{loc series harmonic action} Let the localized $p$-adic pro-unipotent harmonic action of series of $\mathbb{P}^{1} \setminus \{0,\mu_{N},\infty\}$ be the map
\begin{multline}
\circ_{\har,\loc}^{\Sigma} :
K\langle \langle e_{0 \cup \mu_{N}} \rangle \rangle_{\har,o(1)}^{\Sigma}  
\times \big( K\langle\langle e_{0 \cup \mu_{N}}^{\loc} \rangle\rangle \big)^{\mathbb{N}}
\rightarrow \big( K\langle \langle e_{0 \cup \mu_{N}} \rangle \rangle \big)^{\mathbb{N}} 
\\ (g,h) \mapsto 
g \circ_{\har,\loc}^{\Sigma} h =
 \sum_{\substack{0 \leq \tilde{r} \leq m-1 \\ 1 \leq i_{1}< \ldots < i_{\tilde{r}} \leq d
\\ 1 \leq  \tilde{i}_{1} < \ldots < \tilde{i}_{\tilde{r}} \leqslant m-1
\\ \{1,\ldots,d\} - \{i_{1},\ldots,i_{\tilde{r}}\} = S_{0} \amalg \ldots \amalg S_{m-1} \text{ i.c.p.}}} \sum_{l_{I_{1}},\ldots,i_{I_{d-\tilde{r}}} \geqslant 0} \bigg(
\\ \prod_{t=1}^{d-\tilde{r}}
 {-n_{i} \choose l_{i}} h(-\sum_{i \in S_{0}} l_{i},\ldots,-\sum_{i \in S_{M_{1}-1}} l_{i},n_{i_{1}}, -\sum_{i \in M_{1}} l_{i},\ldots,-\sum_{i \in M_{2}-1} l_{i},\ldots, n_{i_{\tilde{r}}},-\sum_{i \in M_{r}} l_{i},\ldots,-\sum_{i \in M_{r+1}-1} l_{i},)
\\ \times 
\prod_{a=0}^{m-1} g(\shft_{l_{1},\ldots,l_{d}}(w)|_{S_{a}}) \bigg),
\end{multline}
where $M_{1},\ldots,M_{r}$ are such that $S_{M_{i}} \amalg S_{M_{i}+1} \amalg \ldots \amalg S_{M_{i+1}-1}=]i_{l},i_{l+1}[$. Then we have
\begin{equation} \label{eq:equation of localized action}\har_{p^{\alpha}\mathbb{N}} = \har_{p^{\alpha}} \text{ } \circ_{\har,\loc}^{\Sigma} \text{ } \frak{h}^{(p^{\alpha})}_{\mathbb{N},\loc}
\end{equation}
\end{Proposition-Definition}

\begin{proof} We write Proposition \ref{multiplication} in the particular case $\mu=p^{\alpha}$, and we apply Lemma \ref{lemma p-adic shifting} to the factors 
$\har_{\mu a,\mu(a+1)}(w|_{S_{a}})$ in the right-hand side of equation (4.4.3), with $\delta=p^{\alpha}a$.
\end{proof}

\begin{Example} In depth one and two and if $N=1$, for all $m,n,n_{1},n_{2} \in \mathbb{N}^{\ast}$, for any $g$ and $h=(h_{m})_{m\in \mathbb{N}^{\ast}}$,
\begin{equation} (g \circ_{\har,\loc}^{\Sigma}h)_{m}(n) = h_{m}(n) + \sum_{l_{1} \in \mathbb{N}} m^{n} h_{l_{1}}(n) {-n \choose l_{1}} g(n+l_{1}),
\end{equation}
\begin{multline}
(g \circ_{\har,\loc}^{\Sigma}h)_{m}(n_{1},n_{2}) = h_{m}(n_{1})\sum_{l_{2} \geqslant 0} h_{m}(l_{2}) {-n_{2} \choose l_{2}} g(n_{2}+l_{2}) + \\ h_{m}(n_{1},n_{2}) +
 \sum_{l_{1},l_{2}\geq 0}  
\prod_{i=1}^{2} {-n_{i} \choose l_{i}} m^{n_{i}}  \times
\bigg[ h_{m}(-l_{1}-l_{2}) g(n_{1}+l_{1},n_{2}+l_{2})
+ h_{m}(-l_{1},-l_{2}) \prod_{i=1}^{2} g(n_{i}+l_{i}) \bigg] 
\\ + m^{n_{1}+n_{2}} \bigg[ \sum_{l_{2}\geq 0} g(n_{2}+l_{2}) {-n_{1} \choose l_{1}} h_{m}(-l_{1},n_{2})
- \sum_{l_{2}\geq 0} g(n_{2}+l_{2}) {-n_{2} \choose l_{2}} \frak{h}_{m}(n_{1},-l_{2}) \bigg].
\end{multline}
\end{Example}

\subsection{The delocalization of localized cyclotomic multiple harmonic sums}

We show that localized cyclotomic multiple harmonic sums $\frak{h}_{m_{0},m}(w)$ can be expressed as linear combinations of cyclotomic multiple harmonic sums over a ring of explicit polynomial-exponential functions of $(m_{0},m)$. This is a series analogue of the fact that an iterated integral of any differential forms on $\mathbb{P}^{1} \setminus \{0,\mu_{N},\infty\}$ can be related to iterated integrals of $\frac{dz}{z}$, $\frac{dz}{z-\xi}$, $\xi \in \mu_{N}(K)$.

\subsubsection{Definition of the localization map and recursive formula}

\begin{Definition} \label{la definition du localise}For $w =\big((n_{i})_{d};(\xi_{i})_{d+1}\big)$ a localized harmonic word, let $\Sign^{-}(w)= \{i \in \{1,\ldots,d\} \text{ | } n_{i}<0\}$, and $\Sign^{+}(w) = \{ i \in \{1,\ldots,d\} \text{ | } n_{i} \geqslant 0\}$.
\end{Definition}

Below we use the notations of Definition \ref{def trois quatre un} and Definition \ref{gen series harm}

\begin{Proposition-Definition} \label{loc for har n} of the following linear map $\deloc$, defined by induction on the depth as follows. Let $\frak{h}_{m_{0},m}(w) = \sum\limits_{w'} \har_{m_{0},m}(w') P_{w'}(m_{0},m)$ is the equality obtained by applying equation (\ref{eq:recurrence}) to $\frak{h}_{m'_{0},m'}(w|_{[i_{C},j_{C}]})$ for all $[i_{C},j_{C}]$ connected components of $\Sign^{-}(w)$ and summing over all the appropriate $(m'_{0},m')$ ($P_{w'}$ is a polynomial-exponential function of $(m_{0},m)$). We let, for any $w \in \Wd_{\har}^{\loc}(e_{0\cup \mu_{N}})$, 
$$\deloc(w) = (\sum_{w'} \loc(w') P_{w'}(m_{0},m))_{(m_{0},m)\in \mathbb{N}^{2},m_{0}<m}. $$
Then, $\deloc$ is well-defined and its dual restricted to the terms $m_{0}=0$, 
$\deloc^{\vee} : (K\langle\langle e_{0\cup \mu_{N}}\rangle\rangle_{\har}^{\Sigma})^{\mathbb{N}} \longrightarrow (K\langle\langle e_{0\cup \mu_{N}}\rangle\rangle_{\har,\loc}^{\Sigma})^{\mathbb{N}}$, satisfies :
\begin{equation} \label{eq:equation of loc} \deloc^{\vee}\har_{\mathbb{N}}^{(p^{\alpha})} = \frak{h}_{\mathbb{N},\loc}^{(p^{\alpha})}.
\end{equation}
\end{Proposition-Definition}

\begin{proof} This follows from the computation of totally negative cyclotomic multiple harmonic sums (Proposition-Definition \ref{numbers mathcal B}) and from the fact that, with the notations of the statement we always have $\depth(w')<\depth(w)$.
\end{proof}

\begin{Example} In depth one and two and if $N=1$, for all $l_{1},l_{2} \in \mathbb{N}$, $n_{1},n_{2} \in \mathbb{N}^{\ast}$, $m \in \mathbb{N}^{\ast}$, we have
\newline 
$\displaystyle \deloc(-l_{1},n_{2})_{0,m} = \left\{ \begin{array}{ll}\sum\limits_{\delta_{1}=1}^{l_{1}+1} \mathcal{B}_{\delta_{1}}^{l_{1}} \big(n_{2}-\delta_{1}\big) & \text{ if }l_{1}+1 \leqslant n_{2} \\ \sum\limits_{\delta_{1}=1}^{n_{2}-1} \mathcal{B}_{\delta_{1}}^{l_{1}} \big(n_{2}-\delta_{1}\big) 
+ \sum\limits_{\tilde{\delta}_{1}=0}^{l_{1}-n_{2}+1}
\sum\limits_{\delta_{2}=1}^{\delta_{1}-n_{2}+1}
\mathcal{B}_{\delta_{1}}^{l_{1}}\mathcal{B}_{\delta_{2}}^{\delta_{1}-n_{2}} m^{\delta_{2}} &\text{ if }l_{1}+1 > n_{2} 
\end{array} \right.$
\newline
$\displaystyle\deloc(n_{1},-l_{2})_{0,m} = \left\{ \begin{array}{ll}
\sum\limits_{\delta=1}^{l_{2}+1}\mathcal{B}_{\delta}^{l_{2}} m^{l_{2}} \big(n_{1}\big) - \sum\limits_{\delta_{2}=1}^{l_{2}+1} \mathcal{B}_{\delta}^{l_{2}} \big( n_{1}-\delta_{2} \big) & \text{ if } l_{2}+1 < n_{1} 
\\ \sum\limits_{\delta=1}^{l_{2}+1}\mathcal{B}_{\delta}^{l_{2}} m^{l_{2}} \big(n_{1}\big)
- \sum\limits_{\delta_{2}=1}^{n_{1}-1} \mathcal{B}_{\delta_{2}}^{l_{1}} \big(n_{1}-\delta_{2}\big) -
\sum\limits_{\tilde{\delta}_{2}=0}^{l_{2}-n_{1}+1}\sum\limits_{\delta_{1}=1}^{\tilde{\delta}_{2}+1} \mathcal{B}_{\tilde{\delta}_{2}+n_{1}}^{(l_{2}-n_{1})+n_{1}} \mathcal{B}_{\delta_{1}}^{\tilde{\delta}_{2}} m^{\delta_{1}} & \text{ if }l_{2}+1 \geqslant n_{1}
\end{array} \right.$
\end{Example}

\subsubsection{Closed formula for the localization map}

\begin{Definition} For $w \in \Wd_{\har}^{\loc}(e_{0 \cup \mu_{N}})$, let $\mathcal{T}(w)$ be the finite tree built inductively as follows : the root of the tree is labeled by $(\Sign^{-}(w),\Sign^{+}(w))$ and, for each vertex $V$ of the tree labeled by a couple of parts $(S^{-},S^{+})$ of $\{1,\ldots,d\}$, if $S^{-} \not= \emptyset$ then, for each $P \subset \partial S^{+}(w)$, we draw an arrow starting from $V$ to a new vertex $V'$, and we label $V'$ by the couple $(P,S^{+} - P)$.
\end{Definition}

Clearly $\mathcal{T}(w)$ depends only on the couple $(\Sign^{-}(w),\Sign^{+}(w))$.

\begin{Example} The trees $\mathcal{T}(w)$ with $w$ of depth 1 are $(1)^{-}$ and $(1)^{+}$. The trees $\mathcal{T}(w)$ with $w$ of depth 2 are $(12)^{-}$ $(12)^{+}$, and the two following ones :
$$		
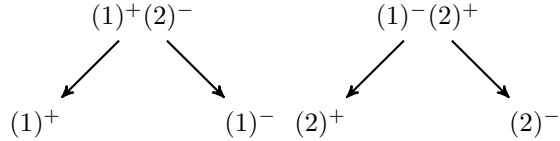
\begin{tikzpicture}[->,>=stealth',shorten >=1pt,auto,node distance=2cm,
		thick,main node/.style={font=\sffamily}]
		\node[main node] (1) {$(1)^{+}(2)^{-}$};
		\node[main node] (2) [below left of=1] {$(1)^{+}$};
		\node[main node] (3) [below right of=1] {$(1)^{-}$};
		\path[every node/.style={font=\sffamily\small}]
		(1) edge node [left] {} (2)
		edge [right] node[left] {} (3)
		;
		\end{tikzpicture} 
\begin{tikzpicture}[->,>=stealth',shorten >=1pt,auto,node distance=2cm,
thick,main node/.style={font=\sffamily}]
\node[main node] (1) {$(1)^{-}(2)^{+}$};
\node[main node] (2) [below left of=1] {$(2)^{+}$};
\node[main node] (3) [below right of=1] {$(2)^{-}$};
\path[every node/.style={font=\sffamily\small}]
(1) edge node [left] {} (2)
edge [right] node[left] {} (3) ;
\end{tikzpicture} $$
The trees $\mathcal{T}(w)$ with $w$ of depth 3 are $(123)^{+}$, $(123)^{-}$, and the six following ones :
$$ 
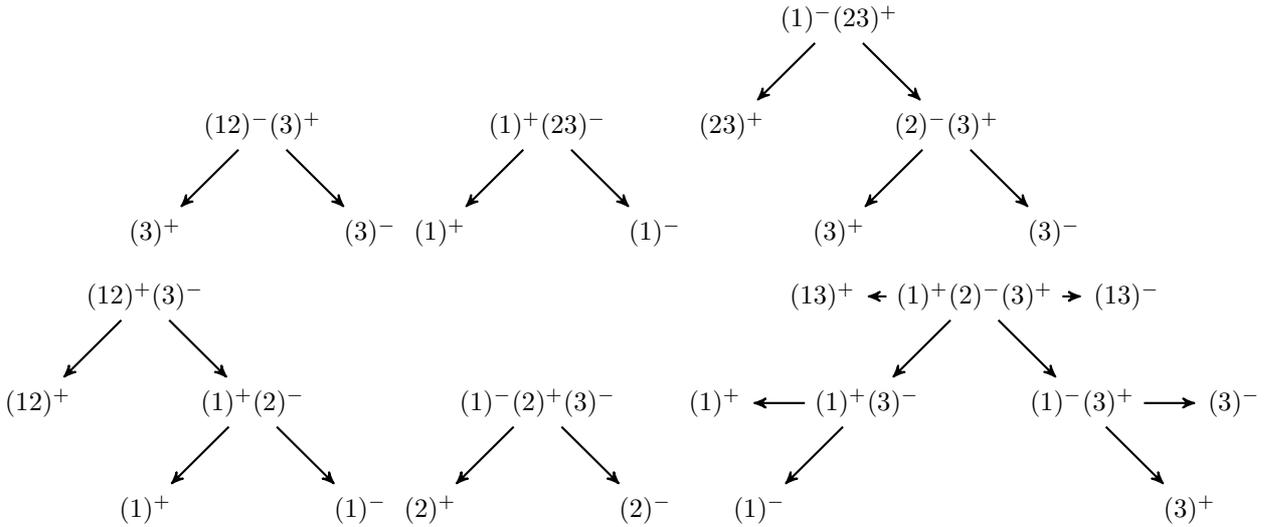
\begin{tikzpicture}[->,>=stealth',shorten >=1pt,auto,node distance=2cm,
thick,main node/.style={font=\sffamily}]
\node[main node] (1) {$(12)^{-}(3)^{+}$};
\node[main node] (2) [below left of=1] {$(3)^{+}$};
\node[main node] (3) [below right of=1] {$(3)^{-}$};
\path[every node/.style={font=\sffamily\small}]
(1) edge node [left] {} (2)
edge [right] node[left] {} (3)
;
\end{tikzpicture}
\begin{tikzpicture}[->,>=stealth',shorten >=1pt,auto,node distance=2cm,
thick,main node/.style={font=\sffamily}]
\node[main node] (1) {$(1)^{+}(23)^{-}$};
\node[main node] (2) [below left of=1] {$(1)^{+}$};
\node[main node] (3) [below right of=1] {$(1)^{-}$};
\path[every node/.style={font=\sffamily\small}]
(1) edge node [left] {} (2)
edge [right] node[left] {} (3);
\end{tikzpicture}
\begin{tikzpicture}[->,>=stealth',shorten >=1pt,auto,node distance=2cm,
thick,main node/.style={font=\sffamily}]
\node[main node] (1) {$(1)^{-}(23)^{+}$};
\node[main node] (2) [below left of=1] {$(23)^{+}$};
\node[main node] (3) [below right of=1] {$(2)^{-}(3)^{+}$};
\node[main node] (4) [below left of=3] {$(3)^{+}$};
\node[main node] (5) [below right of=3] {$(3)^{-}$};
\path[every node/.style={font=\sffamily\small}]
(1) edge node [left] {} (2)
edge [right] node[left] {} (3)
(3) edge node [left] {} (4)
edge [right] node[left] {} (5) ;
\end{tikzpicture} $$
$$ \begin{tikzpicture}[->,>=stealth',shorten >=1pt,auto,node distance=2cm,
thick,main node/.style={font=\sffamily}]
\node[main node] (1) {$(12)^{+}(3)^{-}$};
\node[main node] (2) [below left of=1] {$(12)^{+}$};
\node[main node] (3) [below right of=1] {$(1)^{+}(2)^{-}$};
		\node[main node] (4) [below left of=3] {$(1)^{+}$};
		\node[main node] (5) [below right of=3] {$(1)^{-}$};
		\path[every node/.style={font=\sffamily\small}]
		(1) edge node [left] {} (2)
		edge [right] node[left] {} (3)
		(3) edge node [left] {} (4)
		edge [right] node[left] {} (5) ;
		\end{tikzpicture}
		\begin{tikzpicture}[->,>=stealth',shorten >=1pt,auto,node distance=2cm,
		thick,main node/.style={font=\sffamily}]
		\node[main node] (1) {$(1)^{-}(2)^{+}(3)^{-}$};
		\node[main node] (4) [below left of=1] {$(2)^{+}$};
		\node[main node] (5) [below right of=1] {$(2)^{-}$};
		\path[every node/.style={font=\sffamily\small}]
		(1) 
		edge node [left] {} (4)
		(1) edge [right] node[left] {} (5) ;
		\end{tikzpicture}
  \begin{tikzpicture}[->,>=stealth',shorten >=1pt,auto,node distance=2cm,
 thick,main node/.style={font=\sffamily}]
 \node[main node] (1) {$(1)^{+}(2)^{-}(3)^{+}$};
 \node[main node] (2) [left of=1] {$(13)^{+}$};
 \node[main node] (3) [right of=1] {$(13)^{-}$};
 \node[main node] (4) [below left of=1] {$(1)^{+}(3)^{-}$};
 \node[main node] (5) [below right of=1] {$(1)^{-}(3)^{+}$};
 \node[main node] (6) [left of=4] {$(1)^{+}$};
 \node[main node] (7) [below left of=4] {$(1)^{-}$};		
 \node[main node] (8) [below right of=5] {$(3)^{+}$};
 \node[main node] (9) [right of=5] {$(3)^{-}$};
 ;
 \path[every node/.style={font=\sffamily\small}]
 (1) edge node [left] {} (2)
 edge [right] node[left] {} (3)
 edge node [left] {} (4)
 edge [right] node[left] {} (5)
 (4) edge node [left] {} (6)
 edge [right] node[left] {} (7)
 (5)	edge node [left] {} (8)
 edge [right] node[left] {} (9) ;
 \end{tikzpicture} $$
\end{Example}

We now consider paths from the root to the leaves of a $\mathcal{T}(w)$.

\begin{Definition} For $w \in \Wd_{\har}^{\loc}(e_{0 \cup \mu_{N}})$, let $\mathcal{P}(w)$ be the set of sequences of nodes  $\big((S_{i}^{-},S_{i}^{+})_{i=0,\ldots,u}$ in $\mathcal{T}(w)$ whose first element is the root, whose last element is the leaf, and such that for all $i$, the $i$-th node in the sequence is the son of the $(i-1)$-th node in the sequence.
\newline For all $\big((S_{i}^{-},S_{i}^{+})_{i=0,\ldots,u} \in \mathcal{P}(w)$, let us denote by
\newline (i) $[a^{-}_{1,i},b^{-}_{1,i}]\amalg \ldots \amalg[a^{-}_{r^{-}_{i},i},b^{-}_{r^{-}_{i},i}]$ resp. $[a^{+}_{1,i},b^{+}_{1,i}]\amalg \ldots \amalg[a^{+}_{r^{+}_{i},i},b^{+}_{r^{+}_{i},i}]$ the canonical increasing connected partition of each $S_{i}^{-}$ resp. $S_{i}^{+}$
\newline (ii) $ \{A_{1,i+1}^{-},\ldots,A^{-}_{t^{A,-}_{i+1},i+1}\} = \{a_{1,i}^{-}-1,\ldots,a_{r_{i}^{-},i}^{-}-1\} \cap S_{i+1}^{-}$, $ \{A_{1,i+1}^{+},\ldots,A^{+}_{t^{A,+}_{i+1},i+1}\} = \{a_{1,i}^{-}-1,\ldots,a_{r_{i}^{-},i}^{-}-1\} \cap S_{i+1}^{+}$,  $ \{B_{1,i+1}^{-},\ldots,B^{-}_{t^{B,-}_{i+1},i+1}\} = \{b_{1,i}^{-}+1,\ldots,b_{r_{i}^{-},i}+1\} \cap S_{i+1}^{-}$, $ \{B_{1,i+1}^{+},\ldots,B^{+}_{t^{B,+}_{i+1},i+1}\} = \{b_{1,i}^{-}+1,\ldots,b_{r_{i}^{-},i}+1\} \cap S_{i+1}^{+}$
\newline (iii) $x_{1,i},\ldots,x_{y_{i},i}$ the connected components of $\partial S_{i}^{+}$ which are singletons, and  $\{x_{1,i},\ldots,x_{i,t_{i}}\} = \partial S_{i}^{+} - S_{i+1}^{-}$, with $x_{1,i}<\ldots<x_{i,t_{i}}$.
\end{Definition}

\begin{Definition} For each $w=((n_{i})_{d},(\xi_{i})_{d+1})$ localized harmonic word, for each element of $\mathcal{P}(w)$ as above, and $\underline{\delta}$, sequence of variables in $\mathbb{N}$, and $\underline{j_{0}},\underline{j}$ sequences of functions with values in $\{1,\ldots,N\}$, let
\newline $w^{[i]}(\underline{\delta},\underline{j_{0}},\underline{j})=(n^{[i]}_{a_{1,i}^{+}},\ldots,n^{[i]}_{b_{1,i}^{+}},\ldots\ldots,n^{[i]}_{a_{r_{i}^{+},i}^{+}},\ldots,n^{[i]}_{b_{r_{i}^{+},i}^{+}} ; \xi^{[i]}_{a_{1,i}^{+}},\ldots,\xi^{[i]}_{b_{1,i}^{+}},\ldots\ldots,\xi^{[i]}_{a_{r_{i}^{+},i}^{+}},\ldots,\xi^{[i]}_{b_{r_{i}^{+},i}^{+}})$ where
\newline (i) $(n^{[i]}_{a_{1,i}^{+}},\ldots,n^{[i]}_{b_{1,i}^{+}},\ldots\ldots,n^{[i]}_{a_{r_{i}^{+},i}^{+}},\ldots,n^{[i]}_{b_{r_{i}^{+},i}^{+}})(\underline{\delta}) =
\\ \displaystyle\bigg( n_{a_{1,u}^{+}}-\sum_{\substack{x\in \cup_{i'=1}^{i}(\partial S_{i'}^{+} - S_{i'+1}^{-}) \\ \text{s.t. }x<a^{+}_{1,i}}} \delta_{x},n_{a_{1,i}^{+}+1},\ldots,n_{b_{1,i}^{+}-1},n_{b_{1,i}^{+}}-\sum_{\substack{x\in \cup_{i'=1}^{i}(\partial S_{i'}^{+} - S_{i'+1}^{-}) \\ \text{s.t. }b^{+}_{1,i}<x<a^{+}_{2,i}}}\delta_{0,x},\ldots\ldots
\\ \ldots\ldots,
n_{a_{r^{+}_{i},i}^{+}}-\sum_{\substack{x\in \cup_{i'=1}^{i}(\partial S_{i'}^{+} - S_{i'+1}^{-}) \\ \text{s.t. }b^{+}_{r^{+}_{i}-1,i}<x<a^{+}_{r^{+}_{i},i}}}\delta_{x},n_{a_{r^{+}_{i},i}^{+}+1},\ldots,n_{b_{r^{+}_{i},i}^{+}-1},n_{b_{r^{+}_{i}},i}^{+}-\sum_{\substack{x\in \cup_{i'=1}^{i}(\partial S_{i'}^{+} - S_{i'+1}^{-}) \\ \text{s.t. }b^{+}_{r^{+}_{i},i}<x}}\delta_{x}\bigg)$.
\newline (ii) $(\xi^{[i]}_{a_{1,i}^{+}},\ldots,\xi^{[i]}_{b_{1,i}^{+}},\ldots\ldots,\xi^{[i]}_{a_{r_{i}^{+},i}^{+}},\ldots,\xi^{[i]}_{b_{r_{i}^{+},i}^{+}})(\underline{\delta},\underline{j_{0}},\underline{j}) = 
\\
\displaystyle\bigg( \xi_{a_{1,u}^{+}} \prod_{\substack{1\leq i'\leq i-1\\b_{f,i'}^{-}<a_{1,i}^{+} }} \xi^{j(w^{[i']}|_{[a_{f,i'}^{-},b_{f,i'}^{-}]})} ,\xi_{a_{1,i}^{+}+1},\ldots,\xi_{b_{1,i}^{+}-1},\xi_{b_{1,i}^{+}}\prod_{\substack{1\leq i'\leq i-1\\b_{1,i}^{+}<a_{f,i'}^{-} <b_{f,i'}^{-}<a_{2,i}^{+}}} \xi^{j_{0}(w^{[i']}|_{[a_{f,i'}^{-},b_{f,i'}^{-}]})} ,
\ldots\ldots
\\ \ldots\ldots,
\xi_{a_{r^{+}_{i},i}^{+}}\prod_{\substack{1\leq i'\leq i-1\\a_{r_{i}^{+},i}<a_{f,i'}^{-} <b_{f,i'}^{-}<a_{r_{i}^{+},i}^{+} }} \xi^{j(w^{[i']}|_{[a_{f,i'}^{-},b_{f,i'}^{-}]})} , \xi_{a_{r^{+}_{i},i}^{+}+1},\ldots,\xi_{b_{r^{+}_{i},i}^{+}-1},\xi_{b_{r^{+}_{i}},i}^{+} \prod_{\substack{1\leq i'\leq i-1\\ a_{r_{i}^{+},i}^{+}<a_{f,i'}^{-}}} \xi^{j_{0}(w^{[i']}|_{[a_{f,i'}^{-},b_{f,i'}^{-}]})} \bigg)$.
\end{Definition}

\begin{Proposition} \label{formula for loc}For any $w=\big((n_{i})_{d},(\xi_{i})_{d+1}\big)$, and $m_{0},m \in \mathbb{N}$ with $m_{0}<m$, we have :
\begin{multline} \deloc(w)_{m_{0},m} =
\sum_{(S_{i}^{-},S_{i}^{+})_{i=0,\ldots,u} \in \mathcal{P}(w)}
\sum_{\substack{j_{0}(w^{[i]}|_[a_{f,i}^{-},b_{f,i}^{-}]),j(w^{[i]}|_[a_{f,i}^{-},b_{f,i}^{-}]) \in \{1,\ldots,N\}
\\  
\delta_{0,A_{1,i+1}^{+}} \leqslant n_{A_{1,i+1}^{+}}^{[i]}-1, \ldots, \delta_{0,A^{+}_{t^{A,+}_{i+1},i+1}} \leqslant n_{A^{+}_{t^{A,+}_{i+1},i+1}}^{[i]}-1 
\\ \delta_{B_{1,i+1}^{+}} \leqslant n_{B_{1,i+1}^{+}}^{[i]}-1, \ldots, \delta_{B^{+}_{t^{B,+}_{i+1},i+1}} \leqslant n_{B^{+}_{t^{B,+}_{i+1},i+1}}^{[i]}-1 \\ \delta_{0,A_{1,i+1}^{-}} \geqslant n_{A_{1,i+1}^{-}}^{[i]}, \ldots, \delta_{0,A^{-}_{t^{A,-}_{i+1},i+1}} \geqslant n_{A^{-}_{t^{A,-}_{i+1},i+1}}^{[i]} 
\\ \delta_{B_{1,i+1}^{-}} \geqslant n_{B_{1,i+1}^{-}}^{[i]}, \ldots, \delta_{B^{-}_{t^{B,-}_{i+1},i+1}} \geqslant n_{B^{-}_{t^{B,-}_{i+1},i+1}}^{[i]}}} 
\\
\bigg( \prod_{\substack{0 \leqslant i \leqslant u-1 \\ 1 \leqslant f \leqslant r_{i}^{-}}} \mathcal{B}_{\delta_{0},\delta,j_{0}(w^{[i]}|_[a_{f,i}^{-},b_{f,i}^{-}]),j(w^{[i]}|_[a_{f,i}^{-},b_{f,i}^{-}])}^{w^{[i]}|_{[a_{f,i}^{-},b_{f,i}^{-}]}} \bigg) m_{0}^{n^{[u]}_{a_{1,u}^{+}}}m^{n^{[u]}_{a_{r^{+}_{u},u}^{+}}} (\xi^{[u]}_{a_{1,u}^{+}})^{m_{0}} 
(\xi^{[u]}_{a_{r_{u}^{+}},u})^{m} w_{+}^{[u]}(\underline{\delta},\underline{j_{0}},\underline{j}).
\end{multline}
\end{Proposition}

\begin{proof} By induction using Proposition \ref{loc for har n}.
\end{proof}

\subsection{The pro-unipotent harmonic action of series}

We combine \S4.1 and \S4.2 and we prove equation (\ref{eq: property ac sigma}). In the next statement, we use the notation $g \circ_{\har}^{\Sigma} f = \circ_{\har}^{\Sigma} (g,f)$.

\begin{Proposition-Definition} \label{series harmonic action}Let the $p$-adic pro-unipotent harmonic action of series for $\mathbb{P}^{1} \setminus \{0,\mu_{N},\infty\}$ be the map
\begin{equation} \circ_{\har}^{\Sigma} =  \circ_{\har,\loc}^{\Sigma} \circ (\id \times \loc^{\vee}) : K\langle \langle e_{0 \cup \mu_{N}} \rangle \rangle^{\Sigma}_{\har,o(1)} \times \big( K\langle \langle e_{0 \cup \mu_{N}} \rangle \rangle ^{\Sigma}_{\har}\big)^{\mathbb{N}} \rightarrow 
\big( K\langle \langle e_{0 \cup \mu_{N}} \rangle \rangle^{\Sigma}_{\har} \big)^{\mathbb{N}}.
\end{equation}
Then $\circ_{\har}^{\Sigma}$ is continuous for the $\mathcal{N}_{D}$-topology and satisfies equation (\ref{eq: property ac sigma}).
\end{Proposition-Definition}

\begin{proof} This follows from equation (\ref{eq:equation of localized action}) and equation (\ref{eq:equation of loc}). The convergence of the series involved follows from the bounds on the $p$-adic valuations of the coefficients $\mathcal{B}$ (Proposition-Definition 4.2.2). The continuity is clear.
\end{proof}

Joining the formula for $\circ_{\har,\loc}^{\Sigma}$ (Proposition-Definition \ref{loc series harmonic action}) and the formula for $\deloc$ (Proposition \ref{formula for loc}) we have a formula for $\circ_{\har}^{\Sigma}$.

\begin{Example} \label{example theorem part Sigma} In depth one and two and if $N=1$, for all $n \in \mathbb{N}^{\ast}$, $n_{1},n_{2} \in \mathbb{N}^{\ast}$, for any $g$ and $h=(h_{m})_{m\in\mathbb{N}}$,
\begin{equation} (g \circ_{\har}^{\Sigma} h)(n) = h(n) + \sum_{l \geqslant 1} m^{n+l} \sum_{l_{1}\geq l-1} \mathcal{B}^{l_{1}}_{l} {-n \choose l_{1}} g(n+l_{1}),
\end{equation}
\begin{multline} (g \circ_{\har}^{\Sigma} h)_{m}(n_{1},n_{2}) =
h(n_{1},n_{2}) \text{ }+ 
\\ \sum_{t \geqslant 1} m^{n_{1}+n_{2}+t}
\big[
\sum_{\substack{l_{1},l_{2} \geqslant 0 \\ l_{1}+l_{2} \geqslant t-1}} 
\mathcal{B}_{t}^{l_{1}+l_{2}} 
\prod_{i=1}^{2} {-n_{i} \choose l_{i}} g(n_{1}+l_{1},n_{2}+l_{2}) 
+ \sum_{\substack{l_{1},l_{2} \geqslant 0 \\ l_{1}+l_{2} \geqslant t-2}} 
\mathcal{B}_{t}^{l_{1},l_{2}} 
\prod_{i=1}^{2} {-n_{i} \choose l_{i}} g(n_{i}+l_{i}) \big]
\\ 
+ \sum_{\substack{ 1 \leqslant t \leqslant n_{2}-1 \\ l_{1} \geqslant t-1}} 
m^{n_{1}+t} h_{m}(n_{2}-t) \mathcal{B}_{t}^{l_{1}} {-n_{1} \choose l_{1}} g(n_{1}+l_{1}) - \sum_{\substack{1 \leqslant t \leqslant n_{1}-1 \\ l_{2} \geqslant t-1}}
m^{n_{2}+t} h_{m}(n_{1}-t) \mathcal{B}_{t}^{l_{2}} {-n_{2} \choose l_{2}} g(n_{2}+l_{2})
\\ - m^{n_{2}+n_{1}} \bigg[ \sum_{l_{1} \geqslant n_{2}-1} \mathcal{B}_{n_{2}}^{l_{1}} {-n_{1} \choose l_{1}} g(n_{1}+l_{1})
- \sum_{l_{2} \geqslant n_{1}-1} \mathcal{B}_{n_{1}}^{l_{2}} {-n_{2} \choose l_{2}} g(n_{2}+l_{2})  \bigg] 
\\ 
+ \sum_{t'\geq 1} n^{t'}
\bigg[ \sum_{\substack{t \geqslant n_{2}+t'-1 \\ l_{1} \geqslant t-1}}
\mathcal{B}_{t'}^{t-n_{2}} \mathcal{B}_{t}^{l_{1}} {-n_{1} \choose l_{1}}g(n_{1}+l_{1}) - \sum_{\substack{t \geqslant n_{1}+t'-1 \\ l_{2} \geqslant t-1}}
\mathcal{B}_{t'}^{t-n_{1}} \mathcal{B}_{t}^{l_{2}} {-n_{2} \choose l_{2}}g(n_{2}+l_{2}) \bigg].
\end{multline}
\end{Example}

\begin{Definition} \label{def harmonic Frobenius of series}Let the harmonic Frobenius of series, iterated $\alpha$ times, be the map
\newline $(\phi^{\alpha}_{\har})^{\Sigma} : \begin{array}{cc} \big( K\langle\langle e_{0 \cup \mu_{N}} \rangle\rangle_
{\har}^{\Sigma} \big)^{\mathbb{N}}\rightarrow \big( K\langle\langle e_{0 \cup \mu_{N}} \rangle\rangle_{\har}^{\Sigma}\big)^{\mathbb{N}}
\\  f \mapsto \har_{p^{\alpha}} \text{ }\circ_{\har}^{\Sigma}\text{ } \sigma^{\alpha}(f)\end{array}.$
\end{Definition}

With Definition \ref{def harmonic Frobenius of series}, equation (\ref{eq: property ac sigma}) is restated as 
\begin{equation} (\phi^{\alpha})_{\har}^{\Sigma} (\har_{\mathbb{N}}) = \har_{p^{\alpha}\mathbb{N}}.
\end{equation}

\section{Comparison between results on integrals and on series}

We relate the computations on integrals (\S1,\S2) and the computations on series (\S3,\S4). We prove the part ``comparison between integrals and series'' of the theorem.

\subsection{Maps of comparisons from integrals to series and from series to integrals}

In order to relate the pro-unipotent harmonic action of integrals $\circ_{\har}^{\smallint}$ (Proposition-Definition \ref{dR-rt harmonic Ihara action}) and the pro-unipotent harmonic action of series $\circ_{\har}^{\Sigma}$ (Proposition-Definition \ref{series harmonic action}), we need firstly to extend the definition of $\circ_{\har}^{\smallint}$.

\begin{Proposition-Definition} \label{extension harmonic Ihara action}(i) $K \langle\langle e_{0 \cup \mu_{N}} \rangle\rangle_{o(1)}^{N}$ equipped with 
$(g_{\xi'})_{\xi'\in \mu_{N}(K)} \circ^{\smallint_{0,0}} (f_{\xi})_{\xi'\in \mu_{N}(K)} =$
\newline $(f_{\xi}(e_{0},(g_{\xi'})_{\xi'\in \mu_{N}(K)}))_{\xi \in \mu_{N}(K)}$ is a topological group for the $\mathcal{N}_{D}$-topology.
\newline (ii) The map 
$K \langle\langle e_{0 \cup \mu_{N}} \rangle\rangle_{o(1)}^{N} \times K \langle\langle e_{0 \cup \mu_{N}} \rangle\rangle \rightarrow K \langle\langle e_{0 \cup \mu_{N}} \rangle\rangle$, $(h_{\xi})_{\xi \in \mu_{N}(K)} \circ_{\Ad}^{\smallint_{0,0}} f = f(e_{0},(h_{\xi})_{\xi \in \mu_{N}(K)})$ is a continuous action of the group $K \langle\langle e_{0 \cup \mu_{N}} \rangle\rangle_{o(1)}^{N}$, for the $\mathcal{N}_{D}$-topology.
\newline (iii) Let the extended pro-unipotent harmonic action of integrals be the following map
$$ \circ^{\smallint}_{\har,U} :
\begin{array}{c}
K \langle\langle e_{0 \cup \mu_{N}} \rangle\rangle^{N}_{o(1)} \times
(K \langle\langle e_{0 \cup \mu_{N}} \rangle\rangle_{\har}^{\smallint})^{\mathbb{N}}
\rightarrow 
(K \langle\langle e_{0 \cup \mu_{N}} \rangle\rangle_{\har}^{\smallint})^{\mathbb{N}} 
\\ \big( (g_{\xi})_{\xi \in \mu_{N}(K)}, (h_{m})_{m\in\mathbb{N}} \big) \mapsto g \circ_{\har}^{\smallint} (h_{m})_{m\in\mathbb{N}} = \big( \lim \big( (\tau(m)(g_{\xi}))_{\xi \in \mu_{N}(K)} \circ^{\smallint_{0,0}}_{\Ad} h_{m} \big)\big)_{m\in\mathbb{N}} \end{array}. $$
$\circ^{\smallint}_{\har,U}$ is well-defined and is a continuous group action of the topological group $(K \langle\langle e_{0 \cup \mu_{N}} \rangle\rangle^{N}_{o(1)},\circ^{\smallint_{0,0}})$.
\end{Proposition-Definition}

\begin{proof} (i) and (ii) : the algebraic properties follow from the associativity of the composition of non-commutative formal power series, and the continuity follows as in \S1 ; the topological properties follow from equation (\ref{eq:part c of proof}).
\newline (iii) Same with the proof of Proposition \ref{harmonic Ihara properties}.
\end{proof}

By considering (\ref{eq: property of ac 01}) and (\ref{eq: property ac sigma}), we can now define the maps of comparison between series and integrals.

\begin{Definition} \label{def comp from series to integrals} Let $\comp^{\smallint\Sigma}=(\comp_{\xi}^{\smallint\Sigma})_{\xi \in \mu_{N}(K)} : K\langle\langle e_{0\cup \mu_{N}}\rangle\rangle_{\har,o(1)} \rightarrow K\langle\langle e_{0\cup \mu_{N}}\rangle\rangle_{o(1)}^{N}$ be the map defined as follows : for all $g_{\Sigma}$, and $w_{\har}=\big((n_{i})_{d},(\xi_{i})_{d+1})$, and $w=e_{\xi_{d+1}}e_{0}^{n_{d}-1}e_{\xi_{d}} \ldots e_{0}^{n_{1}-1}e_{\xi_{1}}$,
\newline $(\comp_{\xi}^{\smallint\Sigma}g_{\Sigma})[ e_{0}^{l}e_{\xi_{d+1}}e_{0}^{n_{d}-1}e_{\xi_{d}} \ldots e_{0}^{n_{1}-1}e_{\xi_{1}}]$ is the coefficient of $h_{m_{0},m}(\emptyset)\xi^{-jm}m^{l}$ in the formula for $g_{\Sigma}  \circ_{\har}^{\Sigma}h$, $\comp_{\xi}^{\smallint\Sigma}g_{\Sigma}[e_{0}]=0$ and all other coefficients of $\comp_{\xi}^{\smallint\Sigma}g_{\Sigma}$ are deduced by applying the relation of shuffle modulo products.
\end{Definition}

One can read a formula for $\comp^{\smallint\Sigma}$ via the expression of $\circ_{\har}^{\Sigma}$ explained in \S5.3. We note that writing that formula requires to make a distinction between the words as above for which $l> 0$ and those for which $l=0$.

\begin{Definition} \label{def comp from integrals to series} Let $\comp^{\Sigma\smallint} : K\langle\langle e_{0\cup \mu_{N}}\rangle\rangle_{o(1)}^{N} \rightarrow K\langle\langle e_{0\cup \mu_{N}}\rangle\rangle_{\har,o(1)}$ be defined by
\newline 
$(\comp^{\Sigma\smallint}((g_{\xi})_{\xi \in \mu_{N}(K)}))[e_{\xi_{d+1}}e_{0}^{n_{d}-1}e_{\xi_{d}} \ldots e_{0}^{n_{1}-1}e_{\xi_{1}}] = (-1)^{d} \sum\limits_{\xi \in \mu_{N}(K)} \xi^{-p^{\alpha}}  g_{\xi}[\frac{1}{1-e_{0}}e_{\xi_{d+1}}e_{0}^{n_{d}-1}e_{\xi_{d}}\ldots e_{0}^{n_{1}-1}e_{\xi_{1}}]$.
\end{Definition}

We can now prove equation (\ref{eq:comparison 1}) which relates $\circ_{\har}^{\Sigma}$ and $\circ_{\har}^{\smallint}$.

\begin{proof} The proof is by induction on the depth. Let us mention the two main ingredients of the proof and leave the details to the reader :
\newline (a) $\circ^{\Sigma}_{\har}$ is compatible with restrictions on the domain of summation; namely, the term of $\circ^{\Sigma}_{\har}$ corresponding to a domain of summation of bounds $(m_{0},m)$ and depth $d$ can be computed by computing the term $(m'_{0},m')$ for any $m_{0}<m'_{0}<m'<m$ and depth $d' \leqslant d$, and summing over $(m'_{0},m')$'s and $d'$.
\newline (b) For any $\xi \in \mu_{N}(K)$, $g \in \tilde{\Pi}_{1,0}(K)$,
$\Ad_{g^{(\xi)}}(e_{\xi})$ satisfies the shuffle equation modulo products and $\Ad_{g^{(\xi)}}(e_{\xi})[e_{0}]=0$ ; this implies a formula for all its coefficients in terms of those at words whose rightmost letter is not $e_{0}$ : 
$$\Ad_{g^{(\xi)}}(e_{\xi})
[e_{0}^{n_{d}-1}e_{\xi_{d}}\ldots e_{0}^{n_{1}-1}e_{\xi_{1}}e_{0}^{r}] = \sum\limits_{\substack{l_{1},\ldots,l_{d}\geqslant 0 \\
l_{1}+\ldots+l_{d}=r}} \prod_{i=1}^{d}{-n_{i} \choose l_{i}}
\Ad_{g^{(\xi)}}(e_{\xi})
[e_{0}^{n_{d}+l_{d}-1}e_{\xi_{d}}\ldots e_{0}^{n_{1}+l_{1}-1}e_{\xi_{1}}]. $$
\end{proof}

We now prove equation (\ref{eq:comparison 2}), which relates the two comparison maps.

\begin{proof} For all non-empty totally negative harmonic words $w$,  by Proposition-Definition \ref{numbers mathcal B} and by $\har_{0,1}(w)=0$ (an iterated sum on an empty domain of summation is zero), we have $\sum\limits_{\delta \in\mathbb{N}} \mathcal{B}_{\delta}^{w}=0$. This implies equation (\ref{eq:comparison 2}).
\end{proof}

We now deduce equation (\ref{eq:series expansion}) and (\ref{eq:inversion of series expansion}) which relate $p$-adic cyclotomic multiple zeta values and prime weighted cyclotomic multiple harmonic sums.

\begin{proof}
(a) Equation (\ref{eq:series expansion}) is a consequence of equations (\ref{eq: property of ac 01}), (\ref{eq: property ac sigma}), (\ref{eq:comparison 2})) and the following property. \newline Let $\mathcal{E}^{\smallint}_{\har} \subset \big( K \langle\langle e_{0 \cup \mu_{N}}\rangle\rangle_{\har}^{\smallint}\big)^{\mathbb{N}}$ be the subset introduced in the proof of equation (\ref{eq: property of ac 01}) in \S2.3.2. We have proved in \S2.3.2 that the action $\circ_{\har}^{\smallint}$ restricted to $\mathcal{E}^{\smallint}_{\har}$ is free. This property remains true for the action $\circ_{\har,U}^{\smallint}$ of $K\langle \langle e_{0\cup \mu_{N}}\rangle\rangle_{o(1)}^{N}$ introduced in Proposition-Definition \ref{extension harmonic Ihara action} : indeed, the proof of that property in \S2.3.2 relies on Lemma \ref{depth 0 harmonic Ihara}, which remains true for the extension of $\circ_{\har}^{\smallint}$ introduced in Proposition-Definition \ref{extension harmonic Ihara action}.
\newline (b) Equation (\ref{eq:inversion of series expansion}) is a direct consequence of equation (\ref{eq:comparison 2}) and (\ref{eq:series expansion}).
\end{proof}

By equations (\ref{eq:comparison 1}) and (\ref{eq:series expansion}), the harmonic Frobenius of integrals (Definition \ref{def har Frob}) and the harmonic Frobenius of series (Definition \ref{def harmonic Frobenius of series}) are equal, with the canonical identification $K\langle\langle e_{0\cup \mu_{N}}\rangle\rangle_{\har}^{\smallint}=K\langle\langle e_{0\cup \mu_{N}}\rangle\rangle_{\har}^{\Sigma}$, and can be called ``the harmonic Frobenius'', without ambiguity.

\begin{Remark} The formulas of the theorem can be extended to a formula for the Frobenius itself :
\newline\indent (i) A formula for $\Li_{p,\alpha}^{\dagger}$ in terms of series can be obtained by injecting equation (\ref{eq:series expansion}) in equation (\ref{eq:first}). This enables to interpret in terms of series the parameter $l \in \mathbb{N}^{\ast}$ of the words $e_{0}^{l-1}e_{\xi_{d+1}}e_{0}^{n_{d}-1}e_{\xi_{d}}\ldots e_{0}^{n_{1}-1}e_{\xi_{1}}$, which we have suppressed when we have passed from the Frobenius to the harmonic Frobenius in \S2.
\newline\indent (ii) Let $r \in \{1,\ldots,p^{\alpha}-1\}$. Then, for all $w$, $\har_{r+p^{\alpha}m}[w]$ is a polynomial of values of $\har_{p^{\alpha}m}$ and of analytic functions of $p^{\alpha}m$ whose coefficients are expressed in terms of $\har_{r}$.
\newline We apply the formula of splitting at $p^{\alpha}m$ (\S4.2.2) to express $\har_{p^{\alpha}m+r}$ in terms of $\har_{p^{\alpha}m}$ and $\har_{p^{\alpha}m,p^{\alpha}m+r}$ ; then, the formula of shifting (\S4.2.3) to express $\har_{p^{\alpha}m,p^{\alpha}m+r}$ as an analytic function of $p^{\alpha}m$ with coefficients expressed in terms of $\har_{r}$.
\end{Remark}

\subsection{An adelic interpretation}

Let us now consider all possible values of $p$ and $\alpha$ at the same time : we denote the field $K$ of the previous paragraphs by $K_{p}$, and we let $\mathcal{P}_{N}$ be the set of prime numbers that are prime to $N$. Let also $C_{N}$ be the $N$-th cyclotomic field, embedded diagonally in $\underset{(p,\alpha) \in \mathcal{P}_{N} \times \mathbb{N}^{\ast}}{\prod} K_{p}$. In \cite{I-1}, Definition B.0.3, we have defined, for any positive integer $d$, a $\mathbb{Z}$-module 
$\widehat{\Har}_{\mathcal{P}_{N}^{\mathbb{N}^{\ast}},d}$
as the image of the map $\displaystyle\widehat{\mathcal{O}_{\text{Bound}(d)}^{\sh,e_{0\cup \mu_{N},d}}} \rightarrow \prod_{(p,\alpha)\in \mathcal{P}_{N}\times \mathbb{N}^{\ast}} K_{p}$ which sends $\sum\limits_{n \geqslant 0} w_{n} \mapsto \big( \sum\limits_{n \geqslant 0} \har_{p^{\alpha}}(w_{n}) \big)_{(p,\alpha) \in\mathcal{P}_{N}\times \mathbb{N}^{\ast}}$ ;
here, $\widehat{\mathcal{O}_{\text{Bound}(d)}^{\sh,e_{0\cup \mu_{N},d}}}$ the set of formal infinite sums $\sum\limits_{n \in \mathbb{N}} w_{n}$ where $w_{n}$ is a $C_{N}$-linear combination of words of weight $n$ and depth $\leqslant d$ with coefficients in $\{x  \in C_{N}\text{ }|\text{ }\forall p \in \mathcal{P}_{N}, v_{p}(x) \geqslant - \kappa_{d} - \kappa'_{d}\log(n+\kappa''_{d}) \}$, 
and $\kappa_{d},\kappa'_{d},\kappa''_{d} \in \mathbb{R}^{+\ast}$ are constants defined by the computations of \cite{I-1}. For any positive integer $d$, the rational coefficients in ths sums of series in depth $\leqslant d$ which appear in \S3, \S4, \S5 clearly satisfy the same bounds with those of \cite{I-1}, so we can keep the same constants $\kappa_{d},\kappa'_{d},\kappa''_{d}$. 
\newline\indent We now have not only a formula for $p$-adic cyclotomic multiple zeta values as a sum of series involving prime weighted multiple harmonic sums (equation (\ref{eq:series expansion})), but also a converse formula of the same type : equation (\ref{eq:inversion of series expansion}).

\begin{Definition-Notation} (i) Let us denote by 
$\widehat{\mathcal{Z}}_{\mathcal{P}_{N}^{\mathbb{N}^{\ast}},d}^{\Sigma}=\widehat{\Har}_{\mathcal{P}_{N}^{\mathbb{N}^{\ast}},d}$.
\newline (ii) Let $\widehat{\mathcal{Z}}^{\smallint}_{\mathcal{P}_{N}^{\mathbb{N}^{\ast}},d}$
be the image of the map $\displaystyle\widehat{\mathcal{O}_{\text{Bound}(d)}^{\sh,e_{0\cup \mu_{N},d}}} \rightarrow \prod_{(p,\alpha)\in \mathcal{P}_{N}\times \mathbb{N}^{\ast}} K_{p}$ which sends
\newline $\sum\limits_{n \geqslant 0} w_{n} \mapsto \big( \sum\limits_{n \geqslant 0} \zeta_{p,\alpha}(w_{n}) \big)_{(p,\alpha) \in\mathcal{P}_{N}\times \mathbb{N}^{\ast}}$.
\end{Definition-Notation}

We deduce a last result of comparison between integrals and series :

\begin{Corollary} \label{equality of algebras} We have $\widehat{\mathcal{Z}}^{\smallint}_{\mathcal{P}_{N}^{\mathbb{N}^{\ast}},d} = \widehat{\mathcal{Z}}^{\Sigma}_{\mathcal{P}_{N}^{\mathbb{N}^{\ast}},d}$.
\end{Corollary}

\begin{proof} The inclusion $\subset$ is proved by \cite{I-1} or equation (\ref{eq:series expansion}), combined to the relations between the coefficients of $\Phi_{p,\alpha}$ and $\Phi_{p,\alpha}^{-1}e_{1}\Phi_{p,\alpha}$ explained in \cite{J Assoc}. The inclusion $\supset$ follows from equation (\ref{eq:inversion of series expansion}).
\end{proof}

\section{Application : bounds for the dimension of the spaces of cyclotomic finite multiple zeta values}

The following definition generalizes the notion of finite multiple zeta values introduced by Kaneko and Zagier to the cyclotomic case. Several variants of this definition have appeared in the literature, including in \cite{II-1}.

Let $\mathcal{P}_{N}$ be the set of prime numbers which do not divide $N$.

\begin{Definition} Let $\displaystyle \overline{\mathbb{F}}^{(N)}_{p\rightarrow \infty} = \bigg( \prod_{p\in \mathcal{P}_{N}} \overline{\mathbb{F}}_{p}\bigg) / \bigg( \bigoplus_{p\in \mathcal{P}_{N}} \overline{\mathbb{F}}_{p}\bigg)$ 

Let cyclotomic finite multiple zeta values be the following numbers : for $d \in \mathbb{N}_{\geq 1}$,  $n_{i} \in \mathbb{N}_{\geq 1}$, ($1 \leqslant i \leqslant d$) and $\xi_{i}$ $N$-th roots of unity ($1 \leqslant i \leqslant d$),

$$ \zeta_{f} ((n_{i})_{d};(\xi_{i})_{d+1}) = \bigg( \sum_{0<m_{1}<\cdots < m_{d}<p} \frac{ \big( \frac{\xi_{2}}{\xi_{1}} \big)^{m_{1}} \ldots \big( \frac{\xi_{d+1}}{\xi_{d}} \big)^{m_{d}}  
\big(\frac{1}{\xi_{d+1}} \big)^{m}}{m_{1}^{n_{1}}\ldots m_{d}^{n_{d}}} \bigg)_{p \in \mathcal{P}_{N}} \in \overline{\mathbb{F}}^{(N)}_{p\rightarrow \infty} . $$

\end{Definition}

For any $n \in \mathbb{N}$, we let $Z_{n,f}$, resp. $Z_{n,p}$ be the $K$-vector space generated by finite cyclotomic multiple zeta values, resp. $p$-adic cyclotomic multiple zeta values $\zeta_{p,1}$ of weight $n$. By convention $Z_{0,f}=Z_{0,p}=K$.

For any word $w$, denote by $\zeta_{p}(w) = p^{-\weight(w)} \zeta_{p,1}(w)$.

The following application has been derived by Agaki-Hirose-Yasuda in the $N=1$ case (apparently unpublished). We generalize it to the cyclotomic case.

\begin{Corollary} For all $n \in \mathbb{N}$, we have $\dim Z_{n,f}\leq \dim Z_{n,p}$.
\end{Corollary}

\begin{proof}

By Chatzistamatiou's integrality result \cite{Chatzistamatiou}, for any word $w$, we have, for $p$ large enough, $v_{p}(\zeta_{p}^{\mathrm{KZ}}(w)) \geq \mathrm{weight}(w)$, where $\zeta_{p}^{\mathrm{KZ}}$ means the $p$CMZVs in the sense of Furusho as defined in \cite{Yamashita}. As a consequence, we also have, for any word $w$, for $p$ large enough, $v_{p}(\zeta_{p}(w)) \geq 0$. This is deduced by the formula for the Frobenius of $\Pi_{1,0}(K)$ (equation (1.1.5)) and the fact that the numbers $\zeta_{p}^{KZ}(w)$ are (up to a sign) coefficients of the Frobenius-invariant path in $\Pi_{1,0}(K)$.

Thus, by taking reduction modulo large $p$ in equation (\ref{eq:explicit inversion of series expansion N=1}), and dividing by $p^{n_{1}+\cdots+n_{d}}$, we obtain, for large $p$,

\begin{multline} \sum_{d'=0}^{d}  \xi_{d-d'+1}^{p^{\alpha}} \bigg( \prod_{i=d'+1}^{d} (-1)^{n_{i}}
\bigg) \zeta^{(\xi_{d'+1})}_{p} 
\big((n_{d'+i})_{d-d'} ; (\xi_{d'+1+i})_{d-d'} \big)
\text{ }\zeta^{(\xi_{d'+1})}_{p} \big((n_{i})_{d'}\big)
\\ \equiv \sum_{0<m_{1}<\cdots<m_{d}<p}\frac{ \big( \frac{\xi_{2}}{\xi_{1}} \big)^{m_{1}} \ldots \big( \frac{\xi_{d+1}}{\xi_{d}} \big)^{m_{d}}  
\big(\frac{1}{\xi_{d+1}} \big)^{m}}{m_{1}^{n_{1}}\ldots m_{d}^{n_{d}}} \mod p. \end{multline}

Moreover, we can deduce from Anzawa's theorem \cite{A} that the numbers appearing in the left-hand side of equation (6.0.1) generate the $K$-vector space $Z_{n,p}$ with $n=n_{1}+\cdots+n_{d}$.

Thus the image of map $(\zeta_{p}(w)) \in \prod_{p} K_{p} \mapsto 
(\zeta_{p}(w) \mod p)_{p \in \mathcal{P}_{N}} \in 
\overline{\mathbb{F}}_{p\rightarrow \infty}^{(N)}$ (where $K_{p}$ is the extension of $\mathbb{Q}_{p}$ generated by $N$-th roots of unity) is contained in the $K$-vector space of finite CMZVs. This map is surjective by its definition. Thus we deduce the result.

\end{proof}

Combining this corollary with the upper bounds for $\dim(Z_{n,p})$ obtained from the crystalline realization of mixed Tate motives \cite{Yamashita}, we obtain a motivic upper bound for the dimension of $Z_{n,f}$. An analogue of the conjecture of periosd would be that this upper bound is an equality. Thus we can consider finite cyclotomic multiple zeta values as analogue of periods in the unusual ring $\mathbb{F}^{(N)}_{p\rightarrow \infty}$.

\end{document}